\documentclass[10pt]{amsart}
\usepackage{amscd}
\usepackage{amssymb}
\usepackage[all]{xy}
\usepackage{graphicx}
\author{Amnon Yekutieli}
\title{Deformation Quantization in Algebraic Geometry}
\date{8 July 2005}
\address{Department of  Mathematics,
Ben Gurion University, Be'er Sheva 84105, Israel}
\email{amyekut@math.bgu.ac.il}
\dedicatory{Dedicated to Professor Michael Artin on the Occasion
of his Seventieth Birthday}
\thanks{
Primary: 53D55; Secondary: 14D15, 13D10, 16S80.}
\thanks{ This work was partially supported by the 
US - Israel Binational Science Foundation}

\newtheorem{thm}[equation]{Theorem}
\newtheorem{cor}[equation]{Corollary}
\newtheorem{prop}[equation]{Proposition}
\newtheorem{lem}[equation]{Lemma}
\theoremstyle{definition}
\newtheorem{dfn}[equation]{Definition}
\newtheorem{rem}[equation]{Remark}
\newtheorem{exa}[equation]{Example}

\newtheorem{que}[equation]{Question}

\numberwithin{equation}{section}

\newcommand{\iso}{\xrightarrow{\simeq}}
\newcommand{\inj}{\hookrightarrow}

\newcommand{\xar}{\xrightarrow}
\newcommand{\opn}{\operatorname}
\newcommand{\cat}[1]{\operatorname{\mathsf{#1}}}

\newcommand{\rmitem}[1]{\item[\text{\textup{(#1)}}]}
\newcommand{\mfrak}[1]{\mathfrak{#1}}
\newcommand{\mcal}[1]{\mathcal{#1}}
\newcommand{\msf}[1]{\mathsf{#1}}
\newcommand{\mbf}[1]{\mathbf{#1}}
\newcommand{\mrm}[1]{\mathrm{#1}}
\newcommand{\mbb}[1]{\mathbb{#1}}

\newcommand{\smfrac}[2]{{\textstyle \frac{#1}{#2}}}
\newcommand{\tup}[1]{\textup{#1}}
\newcommand{\bsym}[1]{\boldsymbol{#1}}
\newcommand{\boplus}{\bigoplus\nolimits}

\newcommand{\wtil}[1]{\widetilde{#1}}
\newcommand{\til}[1]{\tilde{#1}}
\newcommand{\what}[1]{\widehat{#1}}
\newcommand{\bra}[1]{\langle #1 \rangle}
\newcommand{\set}[1]{\{ #1 \}}

\newcommand{\K}{\mbb{K} \hspace{0.05em}}
\renewcommand{\d}{\mrm{d}}
\newcommand{\bwedge}{{\textstyle \bigwedge}}
\newcommand{\sprod}{{\textstyle \prod}}
\newcommand{\hatotimes}[1]{\, \what{\otimes}_{#1} \,}

\begin{document}

\begin{abstract}
We study deformation quantizations of the structure sheaf 
$\mcal{O}_X$ of a smooth algebraic variety $X$ in 
characteristic $0$. Our main result is that when 
$X$ is $\mcal{D}$-affine, any formal Poisson structure
on $X$ determines a deformation quantization of $\mcal{O}_X$ 
(canonically, up to gauge equivalence). 
This is an algebro-geometric analogue of Kontsevich's 
celebrated result.
\end{abstract}

\maketitle

\setcounter{section}{-1}
\section{Introduction}

This article began with an attempt to understand the work of 
Kontsevich \cite{Ko1, Ko3}, 
Cattaneo-Felder-To\-ma\-ssini \cite{CFT} and
Nest-Tsygan \cite{NT} on deformation 
quantization of Poisson manifolds. Moreover we tried to see to 
what extent the methods applied in the case of $\mrm{C}^{\infty}$ 
manifolds can be carried over to the algebro-geometric case. 

If $X$ is a $\mrm{C}^{\infty}$ manifold with Poisson structure 
$\alpha$ then there is always a deformation quantization
of the algebra of functions $\mrm{C}^{\infty}(X)$ with 
first order term $\alpha$. 
This was proved by Kontsevich in \cite{Ko1}. Furthermore, 
Kontsevich proved that such a deformation quantization is unique 
in a suitable sense. 

If $X$ is either a complex 
analytic manifold or a smooth algebraic variety then one wants to 
deform the sheaf of functions $\mcal{O}_X$. As might be
expected there are potential obstructions, due to the lack of 
global (analytic or algebraic) functions and 
sections of bundles. The case of a complex analytic manifold with 
holomorphic symplectic structure was treated in \cite{NT}. The 
algebraic case was studied in \cite{Ko3}, where several 
approaches were discussed. In the present paper we take a somewhat 
different direction than \cite{Ko3}.

First let us explain what we mean by deformation quantization in 
the context of algebraic geometry. Let $\K$ be a field of 
characteristic $0$, and let $X$ be a smooth algebraic 
variety over $\K$. The tangent sheaf of $X$ is denoted by 
$\mcal{T}_X$. Given an element 
$\alpha \in \Gamma(X, \bwedge^2_{\mcal{O}_{X}} \mcal{T}_X)$
let $\{ -,- \}_{\alpha}$ be the $\K$-bilinear sheaf morphism 
$\mcal{O}_{X} \times \mcal{O}_{X} \to \mcal{O}_{X}$
defined by 
$\{ f, g \}_{\alpha} := \bra{\alpha, \d(f) \wedge \d(g)}$
for local sections $f, g \in \mcal{O}_{X}$. 
If $\{ -,- \}_{\alpha}$ is a Lie bracket on $\mcal{O}_{X}$
then it is called a {\em Poisson bracket}, $\alpha$ is called a 
{\em Poisson structure on $X$}, and the pair 
$(X, \alpha)$ is called a {\em Poisson variety}. It is known that 
$\alpha$ is a Poisson structure if and only if 
$[\alpha, \alpha]= 0$ for the Schouten-Nijenhuis bracket.

Let $\hbar$ be an indeterminate (the ``Planck constant''). 
A {\em star product} on $\mcal{O}_{X}[[\hbar]]$ is a 
$\K[[\hbar]]$-bilinear sheaf morphism
\[ \star :  \mcal{O}_{X}[[\hbar]] \times \mcal{O}_{X}[[\hbar]]
\to \mcal{O}_{X}[[\hbar]] \]
which makes $\mcal{O}_{X}[[\hbar]]$ 
into a sheaf of associative unital $\K[[\hbar]]$-algebras.
The unit element for $\star$ has to be $1 \in \mcal{O}_{X}$, 
and for any local sections $f, g \in \mcal{O}_X$ their product 
should satisfy $f \star g \equiv f g \ \opn{mod} \hbar$.
Furthermore there is a differential condition: there should be a 
sequence of bi-differential operators
$\beta_{j} : \mcal{O}_{X} \times \mcal{O}_{X} 
\to \mcal{O}_{X}$, such that 
\[ f \star g = f g + \sum_{j = 1}^{\infty} 
\beta_{j}(f, g) \hbar^j \in \mcal{O}_{X}[[\hbar]] . \]
A {\em deformation quantization} of $\mcal{O}_{X}$ is by 
definition a star product on $\mcal{O}_{X}[[\hbar]]$.

Actually there is a more refined notion of deformation 
quantization, which has a local nature; see Section \ref{sec1}. 
In the body of the paper the deformation defined in the previous 
paragraph is referred to as a {\em globally trivialized 
deformation quantization}. However, according to 
Theorem \ref{thm4.5}, if
$\mrm{H}^1(X, \mcal{D}_X) = 0$ then any deformation quantization 
is equivalent to a globally trivialized one. So for the purpose of 
the introduction (cf.\ Theorem \ref{thm0.2} below)
we might as well consider only globally trivialized 
deformation quantizations.

Suppose $\star$ is some star product on $\mcal{O}_{X}[[\hbar]]$.
Given two local sections $f, g \in \mcal{O}_X$ define
$\{ f, g \}_{\star} \in \mcal{O}_X$
to be the unique local section satisfying
\[ \hbar \{ f, g \}_{\star} \equiv f \star g - g \star f
\ \opn{mod} \hbar^2 . \]
This is a Poisson bracket on $\mcal{O}_X$. Note that
$\{ f, g \}_{\star} = \beta_{1} (f, g) - \beta_1(g, f)$.

Let $(X, \alpha)$ be a Poisson variety. A 
{\em deformation quantization of $(X, \alpha)$} is a 
deformation quantization $\star$ such that 
$\{ f, g \}_{\star} = 2 \{ f, g \}_{\alpha}$.

There is an obvious notion of {\em gauge equivalence} for deformation 
quantizations. First we need to define what is a gauge 
equivalence of $\mcal{O}_{X}[[\hbar]]$. This is a 
$\K[[\hbar]]$-linear sheaf automorphism 
$\gamma : \mcal{O}_{X}[[\hbar]] \iso \mcal{O}_{X}[[\hbar]]$ 
of the following form: there is a sequence 
$\gamma_{j} : \mcal{O}_{X} \to \mcal{O}_{X}$ 
of differential operators, such that 
\[ \gamma(f) = f + \sum_{j = 1}^{\infty} \gamma_{j}(f) \hbar^j \]
for all $f \in \mcal{O}_{X}$; and also 
$\gamma(1) = 1$. 
Two star products $\star$ and $\star'$ on $\mcal{O}_{X}[[\hbar]]$
are said to be gauge equivalent if there is some 
gauge equivalence $\gamma$ such that
\[ f \star' g = \gamma^{-1} \big( \gamma(f) \star \gamma(g) \big) 
\]
for all $f, g \in \mcal{O}_{X}$.  

To state the main result of our paper we need 
the notion of formal Poisson structure on $X$.
This is a series
$\alpha = \sum_{k = 1}^{\infty} \alpha_k \hbar^k
\in \Gamma(X, \bwedge^2_{\mcal{O}_{X}} \mcal{T}_{X})[[\hbar]]$
satisfying $[\alpha, \alpha] = 0$.
For instance, if $\alpha_1$ is a Poisson structure then 
$\alpha := \alpha_1 \hbar$ is a formal Poisson structure.
Two formal Poisson structure $\alpha$ and
$\alpha'$ are called gauge equivalent if there is some
$\gamma =  \sum_{k = 1}^{\infty} \gamma_k \hbar^k \in 
\Gamma(X, \mcal{T}_X)[[\hbar]]$
such that 
$\alpha' = \opn{exp}(\opn{ad}(\gamma))(\alpha)$. 

Recall that the variety $X$ is said to be {\em $\mcal{D}$-affine} 
if $\mrm{H}^i(X, \mcal{M}) = 0$ for all quasi-coherent left 
$\mcal{D}_X$-modules $\mcal{M}$ and all $i > 0$. Here $\mcal{D}_X$ 
is the sheaf of differential operators on $X$. 

\begin{thm} \label{thm0.2}
Let $X$ be a smooth algebraic variety over the field $\K$. 
Assume $X$ is $\mcal{D}$-affine and $\mbb{R} \subset \K$. 
Then there is a canonical function
\[ Q : \frac{ \{ \tup{formal Poisson structures on $X$} \} }
{\tup{gauge equivalence}} \iso
\frac{ \{ \tup{deformation quantizations of $\mcal{O}_X$} \} }
{\tup{gauge equivalence}} \]
called the {\em quantization map}. The map $Q$
preserves first order terms, and
commutes with \'etale morphisms $X' \to X$.
If $X$ is affine then $Q$ is bijective. 
There is an explicit formula for $Q$. 
\end{thm}

This is an algebraic analogue of \cite[Theorem 1.3]{Ko1}.
Theorem \ref{thm0.2} is repeated as Corollaries \ref{cor9.1} and 
\ref{cor7.9} in the body of the 
paper. Full details, including the explicit 
formula for the quantization map $Q$, are in Theorem \ref{thm5.2}.
By ``preserving first order terms'' we mean that given a formal 
Poisson structure 
$\alpha = \sum_{j = 1}^{\infty} \alpha_j \hbar^j$
and associated deformation quantization 
$Q(\alpha) = \star$, then 
$\{ -, - \}_{\star} =  2 \{ -, - \}_{\alpha_1}$. 
If $f : X' \to X$ is an \'etale morphism then any 
formal Poisson structure $\alpha$ on $X$ can be pulled back 
to a formal Poisson structure $f^*(\alpha)$ on $X'$; and likewise
any deformation quantization $\star$ on $X$ can be pulled back to 
a deformation quantization $f^*(\star)$ on $X'$. 
The third assertion in Theorem \ref{thm0.2} 
says that if $X'$ is also $\mcal{D}$-affine 
then $Q(f^*(\alpha)) = f^*(Q(\alpha))$.

There are two important classes of varieties satisfying the 
conditions of Theorem \ref{thm0.2}. 
The first consists of all 
smooth affine varieties. Note that even if $X$ is affine, yet does 
not admit an \'etale morphism $X \to \mbf{A}^n_{\K}$, the 
result is not trivial -- since changes of coordinates have 
to be accounted for (cf.\ Corollary \ref{cor6.1}).

The second class of examples is that of the flag varieties 
$X = G / P$, where $G$ is a connected reductive algebraic group 
and $P$ is a parabolic subgroup. By the Beilinson-Bernstein 
Theorem the variety $X$ is $\mcal{D}$-affine.
This class of varieties includes the projective spaces 
$\mbf{P}_{\K}^n$. 

Here is an outline of the paper (with some of the features 
simplified). There are two important
sheaves of DG Lie algebras on $X$: the sheaf of poly vector fields
$\mcal{T}_{\mrm{poly}, X}$, and the sheaf of poly differential 
operators $\mcal{D}_{\mrm{poly}, X}$ (see Section \ref{sec3}). 
Their global sections
$\mcal{T}_{\mrm{poly}}(X) := \Gamma(X, \mcal{T}_{\mrm{poly}, X})$
and
$\mcal{D}_{\mrm{poly}}(X) := \Gamma(X, \mcal{D}_{\mrm{poly}, X})$
control Poisson structures and deformation quantizations 
respectively. If one could find an $\mrm{L}_{\infty}$ 
quasi-isomorphism 
$\mcal{T}_{\mrm{poly}}(X) \to \mcal{D}_{\mrm{poly}}(X)$
this would imply Theorem \ref{thm0.2}. However, unless $X$ is 
affine and admits an \'etale morphism to $\mbf{A}^n_{\K}$, 
there is no reason why such a quasi-isomorphism should exist. 

Imitating Fedosov \cite{Fe} and Kontsevich \cite{Ko1}, 
we use formal geometry to solve the global 
problem. The adaptation of this theory to algebraic geometry is 
done in Section \ref{sec4}. 
There is an infinite dimensional bundle 
$\pi : \opn{LCC} X \to X$, which parameterizes formal coordinate 
systems on $X$ modulo linear change of 
coordinates. (In \cite{Ko1} the notation for this bundle
is $X^{\mrm{aff}}$.) Let $\mcal{P}_X$ be the sheaf of principal 
parts on $X$. The complete pullbacks 
$\pi^{\what{*}} (\mcal{P}_X \otimes_{\mcal{O}_X} 
\mcal{T}_{\mrm{poly}, X})$
and
$\pi^{\what{*}} (\mcal{P}_X \otimes_{\mcal{O}_X} 
\mcal{D}_{\mrm{poly}, X})$
are sheaves of DG Lie algebras on $\opn{LCC} X$
(see Section \ref{sec5}). 
The universal deformation formulas of Kontsevich give rise to an 
$\mrm{L}_{\infty}$ quasi-isomorphism 
\[ \pi^{\what{*}} (\mcal{P}_X \otimes_{\mcal{O}_X} 
\mcal{T}_{\mrm{poly}, X}) \to 
\pi^{\what{*}} (\mcal{P}_X \otimes_{\mcal{O}_X} 
\mcal{D}_{\mrm{poly}, X}) . \]

When $X$ is a $\mrm{C}^{\infty}$ manifold the bundle $\opn{LCC} X$ 
has contractible fibers, and thus it has global
$\mrm{C}^{\infty}$ sections. 
This fact is crucial for Kontsevich's proof. 
However, in our algebraic setup there is no reason to assume that 
$\pi : \opn{LCC} X \to X$ has any global sections.

We discovered a way to get around the absence of global 
sections in the case of an algebraic variety: the idea is to use 
{\em simplicial sections}. This idea is inspired by a construction 
of Bott; see \cite{Bo, HY}.
A simplicial section $\bsym{\sigma}$ of 
$\pi : \opn{LCC} X \to X$, based on an open  
covering $X = \bigcup U_{(i)}$, consists of a family of morphisms 
$\sigma_{\bsym{i}} : \bsym{\Delta}^q_{\K} \times 
U_{\bsym{i}} \to \opn{LCC} X$, 
where $\bsym{i} = (i_0, \ldots, i_q)$ is a multi-index; 
$\bsym{\Delta}^q_{\K}$ is the $q$-dimensional geometric simplex; 
and 
$U_{\bsym{i}} := U_{(i_{0})} \cap \cdots \cap U_{(i_{q})}$. 
The morphisms $\sigma_{\bsym{i}}$ are required to be compatible 
with $\pi$ and to satisfy simplicial relations. 

It is easy to show that sections of $\pi$ exist locally. 
Because of the particular geometry of the 
bundle $\opn{LCC} X$, if we take a sufficiently fine affine open 
covering $X = \bigcup U_{(i)}$, and choose a section 
$\sigma_{(i)} : U_{(i)} \to \opn{LCC} X$ for each $i$, 
then these sections can be extended to a 
simplicial section $\bsym{\sigma}$. 
(See Figure \ref{fig1} for an illustration. 
The details of this 
construction are worked out in the companion paper \cite{Ye4}.)

In order to make use of the simplicial section we need {\em mixed 
resolutions}. The mixed resolution of a quasi-coherent 
$\mcal{O}_X$-module $\mcal{M}$ is a complex
$\opn{Mix}_{\bsym{U}}(\mcal{M})$, which combines a de Rham type 
differential related to $\mcal{P}_X$, called the Grothendieck 
connection, together with a \v{C}ech-simplicial 
type differential related to the covering 
$\bsym{U} = \{ U_{(i)} \}$. 
(See Section \ref{sec6} for a review of mixed resolutions.)
We show that the inclusions 
$\mcal{T}_{\mrm{poly}, X} \to 
\opn{Mix}_{\bsym{U}}(\mcal{T}_{\mrm{poly}, X})$
and
$\mcal{D}_{\mrm{poly}, X} \to 
\opn{Mix}_{\bsym{U}}(\mcal{D}_{\mrm{poly}, X})$
are quasi-isomorphisms of sheaves of DG Lie algebras. 
We then prove the following result (which is Theorem \ref{thm5.1} 
in the body of the paper). 

\begin{thm}
Let $\K$ be a field containing $\mbb{R}$, and 
let $X$ be a smooth $n$-dimensional algebraic variety over $\K$. 
Suppose $\bsym{U} = \{ U_{(i)} \}$ is an open covering of $X$,
where each $U_{(i)}$ is affine and admits an \'{e}tale morphism to 
$\mbf{A}^n_{\K}$. Let $\bsym{\sigma}$ be the corresponding
simplicial section of $\pi : \opn{LCC} X \to X$. 
Then there is an induced $\mrm{L}_{\infty}$ quasi-isomorphism 
\[ \Psi_{\bsym{\sigma}} : 
\opn{Mix}_{\bsym{U}}(\mcal{T}_{\mrm{poly}, X}) \to 
\opn{Mix}_{\bsym{U}}(\mcal{D}_{\mrm{poly}, X})  \]
between sheaves of DG Lie algebras.
\end{thm}

We should point out that the construction of the $\mrm{L}_{\infty}$ 
morphism $\Psi_{\bsym{\sigma}}$ involves twisting, due to the 
presence of the Grothendieck connection in the mixed resolution
$\opn{Mix}_{\bsym{U}}(-)$. This sort of twisting is 
discussed in detail in the companion paper \cite{Ye2}. 

Passing to global sections we obtain an $\mrm{L}_{\infty}$ 
quasi-isomorphism 
\[ \Gamma(X, \Psi_{\bsym{\sigma}}) : 
\Gamma \bigl( X, \opn{Mix}_{\bsym{U}}(\mcal{T}_{\mrm{poly}, X})
\bigr) \to 
\Gamma \bigl( X, \opn{Mix}_{\bsym{U}}(\mcal{D}_{\mrm{poly}, X})
\bigr) . \]
There are DG Lie algebra homomorphisms
\begin{equation} \label{eqn0.2}
\mcal{T}_{\mrm{poly}}(X) \to 
\Gamma \bigl( X, \opn{Mix}_{\bsym{U}}(\mcal{T}_{\mrm{poly}, X})
\bigr)
\end{equation}
and
\begin{equation} \label{eqn0.1}
\mcal{D}_{\mrm{poly}}(X) \to 
\Gamma \bigl( X, \opn{Mix}_{\bsym{U}}(\mcal{D}_{\mrm{poly}, X})
\bigr) .
\end{equation}
Each sheaf $\mcal{D}^p_{\mrm{poly}, X}$ is a quasi-coherent left 
$\mcal{D}_{X}$-module. Hence if $X$ is $\mcal{D}$-affine the 
homomorphism (\ref{eqn0.1}) is a quasi-isomorphism. 
By standard results of deformation theory (that are reviewed in 
Section \ref{sec3}) this implies the existence of the quantization 
map $Q$ in Theorem \ref{thm0.2}.
In case $X$ is affine the homomorphism (\ref{eqn0.2}) 
is also a quasi-isomorphism, and thus $Q$ is bijective. 

An earlier version of this paper was much longer. The current 
version contains only the main results; auxiliary results were 
moved to the companion papers \cite{Ye2}, \cite{Ye3} and 
\cite{Ye4}. 
  
Finally let us mention several recent papers and surveys
dealing with deformation quantization: \cite{BK1, BK2}, 
\cite{CDH}, \cite{CF}, \cite{CI}, \cite{Do1, Do2}, and \cite{Ke}. 

\medskip \noindent
\textbf{Acknowledgments.}
The author wishes to thank Michael Artin, 
Joseph Bernstein, Paul Bressler, Vasiliy Dolgushev, 
Murray Gerstenhaber, Vladimir Hinich, 
David Kazhdan, Bernhard Keller, 
Maxim Kontsevich, Monique Lejeune-Jalabert, Zinovy Reichstein,
Pierre Schapira, James Stasheff, Michel Van den Bergh 
and James Zhang for helpful conversations. Thanks also to the 
referee for reading the paper carefully and offering suggestions 
and corrections.
Some of the work was done during visits to MIT, 
the University of Washington, Institut Mittag-Leffler 
and Universit\'e Denis Diderot 
(Paris 7), and their hospitality is gratefully acknowledged.

\section{Deformation Quantizations of $\mcal{O}_X$}
\label{sec1}

Throughout the paper $\K$ is a field of characteristic $0$.
By default all algebras and schemes in the paper are over $\K$,
and so are all morphisms. The symbol $\otimes$ denotes 
$\otimes_{\K}$. The letter $\hbar$ denotes an 
indeterminate, and $\K[[\hbar]]$ is the power series algebra. 

Let $X$ be a smooth separated irreducible $n$-dimensional scheme 
over $\K$. 

\begin{dfn} \label{dfn6.3}
Let $U \subset X$ be an open set. 
A {\em star product} on $\mcal{O}_{U}[[\hbar]]$ is a 
$\K[[\hbar]]$-bilinear sheaf morphism
\[ \star : \mcal{O}_{U}[[\hbar]] \times \mcal{O}_{U}[[\hbar]]
\to \mcal{O}_{U}[[\hbar]] \]
satisfying the following conditions:
\begin{enumerate}
\rmitem{i} The product $\star$ makes $\mcal{O}_{U}[[\hbar]]$ into 
a sheaf of associative unital $\K[[\hbar]]$-algebras with unit
$1 \in \mcal{O}_{U}$. 
\rmitem{ii} There is a sequence
$\beta_j : \mcal{O}_U \times \mcal{O}_U \to \mcal{O}_U$
of bi-differential operators, such that for any two local 
sections $f, g \in \mcal{O}_U$ one has
\[ f \star g = f g + 
\sum_{j = 1}^{\infty} \beta_j(f, g) \hbar^j . \]
\end{enumerate}
\end{dfn}

Note that $f \star g \equiv f g \, \opn{mod} \hbar$, and also
$\beta_j(f, 1) = \beta_j(1, f) = 0$ for all $f$ and $j$.

\begin{dfn} \label{dfn4.3}
Let $\mcal{A}$ be a sheaf of $\hbar$-adically complete flat 
$\K[[\hbar]]$-algebras on $X$, and let
$\psi : \mcal{A} / \hbar \mcal{A} \iso \mcal{O}_X$
be an isomorphism of sheaves of $\K$-algebras.
Let $U \subset X$ be an open set. A {\em differential 
trivialization of $(\mcal{A}, \psi)$ on $U$} 
is an isomorphism 
\[ \tau : \mcal{O}_U [[\hbar]] \iso \mcal{A}|_U \]
of sheaves of $\K[[\hbar]]$-modules satisfying the conditions 
below. \begin{enumerate}
\rmitem{i} Let $\star$ denote the product of $\mcal{A}$. Then the 
$\K[[\hbar]]$-bilinear 
product $\star_{\tau}$ on $\mcal{O}_U [[\hbar]]$, defined by
\[ f \star_{\tau} g := \tau^{-1} \big( \tau(f) \star \tau(g) \big) 
\]
for local sections $f, g \in \mcal{O}_U$, is a star product.
\rmitem{ii} For any local section $f \in \mcal{O}_U$ one 
has $(\psi \circ \tau)(f) = f$.
\end{enumerate}
\end{dfn} 

Condition (i) implies that $\tau(1_{\mcal{O}_X}) = 1_{\mcal{A}}$,
where $1_{\mcal{O}_X}$ and $1_{\mcal{A}}$ are the unit elements of
$\mcal{O}_X$ and $\mcal{A}$ respectively. 

\begin{dfn} \label{dfn4.2}
Let $U \subset X$ be an open set. A {\em gauge equivalence} 
of $\mcal{O}_U [[\hbar]]$ is a $\K[[\hbar]]$-linear automorphism
of sheaves
\[ \gamma : \mcal{O}_U [[\hbar]] \iso \mcal{O}_U [[\hbar]] \]
satisfying these conditions:
\begin{enumerate}
\rmitem{i} There is a sequence of differential operators
$\gamma_k : \mcal{O}_U \to \mcal{O}_U$, such that for 
any local section $f \in \mcal{O}_U$
\[ \gamma(f) = f+ \sum_{k = 1}^{\infty} \gamma_k(f) \hbar^k . \]
\rmitem{ii} 
$\gamma(1) = 1$.
\end{enumerate}
\end{dfn}

Condition (ii) is equivalent to 
$\gamma_k(1) = 0$ for all $k$.
The gauge equivalences of $\mcal{O}_U [[\hbar]]$ form a group 
under composition.

\begin{dfn} \label{dfn4.4}
Let $(\mcal{A}, \psi)$ be as in Definition \ref{dfn4.3}.
A {\em differential structure} 
$\bsym{\tau} = \{ \tau_i \}$ on $(\mcal{A}, \psi)$
consists of an open covering
$X = \bigcup_{i} U_i$, and 
for every $i$ a differential trivialization
$\tau_i : \mcal{O}_{U_i} [[\hbar]] \iso 
\mcal{A}|_{U_i}$
of  $(\mcal{A}, \psi)$ on $U_i$.
The condition is that for any two indices $i, j$ the transition 
automorphism $\tau_j^{-1} \circ \tau_i$ of 
$\mcal{O}_{U_i \cap U_j} [[\hbar]]$
is a gauge equivalence.
\end{dfn}

\begin{exa}
If $\mcal{A}$ is commutative then automatically it has a 
differential structure $\bsym{\tau} = \{ \tau_i \}$, 
with the additional property that each differential trivialization  
$\tau_i : \mcal{O}_{U_i} [[\hbar]] \iso \mcal{A}|_{U_i}$
is an isomorphism of algebras. Here $\mcal{O}_{U_i} [[\hbar]]$
is the usual power series algebra. Let us explain how this is done.
Choose an affine open covering $X = \bigcup U_i$. For any $i$ let 
$C_i := \Gamma(U_i, \mcal{O}_{X})$. 
By formal smoothness of $\K \to C_i$ the isomorphism
$\psi^{-1} : C_i \iso \Gamma(U_i, \mcal{A}) / (\hbar)$
lifts to an isomorphism of algebras 
$\tau_i : C_i[[\hbar]] \iso \Gamma(U_i, \mcal{A})$.
Due to commutativity the isomorphism 
$\tau_i$ sheafifies to a 
differential trivialization on $U_i$. Commutativity also implies 
that the transitions $\tau_j^{-1} \circ \tau_i$
are gauge equivalences. 
The differential structure $\bsym{\tau}$ is 
unique up to gauge equivalence (see Definition \ref{dfn4.6}
below). The first order terms of the gauge 
equivalences $\tau_j^{-1} \circ \tau_i$ are derivations, 
and they give the deformation class of $\mcal{A}$ in
$\mrm{H}^1(X, \mcal{T}_X)$. 
\end{exa}

For a noncommutative algebra 
$\mcal{A}$ it seems that {\em we must stipulate the existence of a 
differential structure}. Furthermore a given algebra 
$\mcal{A}$ might have distinct differential structures. 
Thus we are led to the next definition.

\begin{dfn} \label{dfn4.5}
A {\em deformation quantization of $\mcal{O}_X$} is the data
$(\mcal{A}, \psi, \bsym{\tau})$, where
$\mcal{A}$ is a sheaf of $\hbar$-adically complete flat 
$\K[[\hbar]]$-algebras on $X$; 
$\psi : \mcal{A} / \hbar \mcal{A} \iso \mcal{O}_X$
is an isomorphism of sheaves of $\K$-algebras;
and $\bsym{\tau}$ is a differential structure on
$(\mcal{A}, \psi)$.
\end{dfn}

If there is no danger of confusion we shall sometimes just say 
that $\mcal{A}$ is a deformation quantization, keeping the rest of 
the data implicit.

\begin{exa}
Let $Y$ be a smooth variety and $X := \mrm{T}^* Y$ the cotangent 
bundle, with projection $\pi : X \to Y$. 
$X$ is a symplectic variety, so it has a non-degenerate 
Poisson structure $\alpha$. Let
$\mcal{B} := \boplus_{i = 0}^{\infty} 
(\mrm{F}_i \mcal{D}_Y) \hbar^i
\subset \mcal{D}_Y[\hbar]$
be the Rees algebra of $\mcal{D}_Y$ w.r.t.\ the order filtration
$\{ \mrm{F}_i \mcal{D}_Y \}$. 
So $\mcal{B} / \hbar \mcal{B} \cong \pi_* \mcal{O}_{X}$
and $\mcal{B} / (\hbar - 1) \mcal{B} \cong \mcal{D}_{Y}$.
Define $\mcal{B}_m := \mcal{B} / \hbar^{m+1} \mcal{B}$.
Consider the sheaf of $\K[\hbar]$-algebras 
$\pi^{-1} \mcal{B}_m$ on $X$.
It can be localized to a sheaf of algebras $\mcal{A}_m$
on $X$ such that 
$\mcal{A}_m \cong \pi^{*} \mcal{B}_m = 
\mcal{O}_{X} \otimes_{\pi^{-1} \mcal{O}_Y} 
\pi^{-1} \mcal{B}_m$
as left $\mcal{O}_{X}$-modules. In particular 
$\mcal{A}_0 \cong \mcal{O}_{X}$. 
Let $\mcal{A} := \lim_{\leftarrow m} \mcal{A}_m$.
Then $\mcal{A}$ is a deformation quantization of $(X, \alpha)$. 
Note the similarity to microlocal differential operators \cite{Sch}.
\end{exa}

\begin{dfn} \label{dfn4.6}
Suppose $(\mcal{A}, \psi, \bsym{\tau})$ and
$(\mcal{A}', \psi', \bsym{\tau}')$
are two deformation quantizations of $\mcal{O}_X$. A 
{\em gauge equivalence}
\[ \gamma : (\mcal{A}, \psi, \bsym{\tau}) \to
(\mcal{A}', \psi', \bsym{\tau}') \]
is a isomorphism
$\gamma : \mcal{A} \iso \mcal{A}'$
of sheaves of $\K[[\hbar]]$-algebras satisfying the following two 
conditions:
\begin{enumerate}
\rmitem{i} One has 
$\psi =  \psi' \circ \gamma : \mcal{A} \to \mcal{O}_X$.
\rmitem{ii} Let $\{ U_i \}$ and $\{ U'_j \}$ be the open coverings 
associated to $\bsym{\tau}$ and $\bsym{\tau}'$ respectively. Then 
for any two indices $i, j$ the automorphism
$\tau_j'^{-1} \circ \gamma \circ \tau_i$
of $\mcal{O}_{U_i \cap U'_j}[[\hbar]]$
is a gauge equivalence, in the sense of Definition \ref{dfn4.2}.
\end{enumerate}
\end{dfn}

Let $\Omega^{p}_X = \Omega^{p}_{ X / \K}$ be the sheaf of 
differentials of degree $p$, and let 
$\mcal{T}_{X} = \mcal{T}_{X / \K}$ be the tangent sheaf of $X$.
For every $p \geq 0$ there is a canonical pairing
\[ \bra{-,-} : (\bwedge^p_{\mcal{O}_{X}} \mcal{T}_{X})
\times \Omega^{p}_X \to \mcal{O}_{X}  . \]

\begin{dfn} \label{dfn7.2}
\begin{enumerate}
\item A {\em Poisson bracket} on $\mcal{O}_{X}$ is a biderivation
\[ \set{-, -} : \mcal{O}_{X} \times \mcal{O}_{X} \to 
\mcal{O}_{X} \]
which makes $\mcal{O}_{X}$ into a sheaf of Lie algebras. 
\item Let 
$\alpha \in \Gamma(X, \bwedge^2_{\mcal{O}_{X}} \mcal{T}_{X})$.
Define a biderivation $\set{-, -}_{\alpha}$ by the formula
\[ \set{f, g}_{\alpha} := \bra{\alpha, \d(f) \wedge \d{g} } \]
for local sections $f, g \in \mcal{O}_X$.
If $\set{-, -}_{\alpha}$ is a Poisson bracket then 
$\alpha$ is called a {\em Poisson structure}, and
$(X, \alpha)$ is called a {\em Poisson variety}. 
\end{enumerate}
\end{dfn}

The next result is an easy calculation.

\begin{prop} \label{prop7.4}  
Let $(\mcal{A}, \psi, \bsym{\tau})$ be a deformation 
quantization of $\mcal{O}_X$, and denote by $\star$ the 
multiplication of $\mcal{A}$. Given two local sections 
$f, g \in \mcal{O}_X$ choose liftings
$\til{f}, \til{g} \in \mcal{A}$. Then the formula
\[ \{ f, g \}_{\mcal{A}} := \psi \bigl( \hbar^{-1}
(\til{f} \star \til{g} - \til{g} \star \til{f}) \bigr) 
\in \mcal{O}_X \]
defines a Poisson bracket on $\mcal{O}_X$. 
\end{prop}

Suppose $\{ U_i \}$ is the open covering associated with the 
differential trivialization $\bsym{\tau}$, and for each $i$ the
collection of bi-differential operators on $\mcal{O}_{U_i}$ 
occurring in Definition \tup{\ref{dfn6.3}} is
$\{ \beta_{i, j} \}_{j = 1}^{\infty}$. Then for local 
sections $f, g \in \mcal{O}_{U_i}$ one has
\[ \{ f, g \}_{\mcal{A}} = \beta_{i, 1} (f, g) - \beta_{i, 1} (g, f)
. \]

\begin{dfn} \label{dfn0.1}
Let $\alpha$ be a Poisson structure on $X$. 
A {\em deformation quantization} of 
the Poisson variety $(X, \alpha)$ is a 
deformation quantization $(\mcal{A}, \psi, \bsym{\tau})$ of 
$\mcal{O}_X$ such that the Poisson brackets satisfy
\[ \{ -,- \}_{\mcal{A}} = 2 \{ -,- \}_{\alpha} . \]
\end{dfn}

\begin{dfn}
A {\em globally trivialized deformation quantization of 
$\mcal{O}_X$} is a deformation quantization 
$(\mcal{A}, \psi, \bsym{\tau})$ 
in which the differential structure $\bsym{\tau}$ consists of a 
single differential trivialization 
$\tau : \mcal{O}_X [[\hbar]] \iso \mcal{A}$.
\end{dfn}

In effect a globally trivialized deformation quantization of 
$\mcal{O}_X$ is the same as a star product on 
$\mcal{O}_X [[\hbar]]$; the correspondence is 
$\star \mapsto \star_{\tau}$ in the notation of Definition 
\ref{dfn4.3}.  

Let $\mcal{D}_X$ be the sheaf of $\K$-linear
differential operators on $X$.

\begin{thm} \label{thm4.5}
Assume $\mrm{H}^1(X, \mcal{D}_X) = 0$. Then any
deformation quantization $(\mcal{A}, \psi, \bsym{\tau})$ 
of $\mcal{O}_X$ can be globally trivialized.
Namely there is a globally trivialized deformation quantization
$(\mcal{A}', \psi', \{ \tau' \})$ 
of $\mcal{O}_X$, and a gauge equivalence
$(\mcal{A}, \psi, \bsym{\tau}) \to 
(\mcal{A}', \psi', \{ \tau' \})$.
\end{thm}

\begin{proof}
We will take $(\mcal{A}' ,\psi') := (\mcal{A}, \psi)$,
and produce a global differential trivialization $\tau'$. 

By refining the open covering 
$\bsym{U} = \{ U_i \}$ associated with the differential 
structure $\bsym{\tau} = \{ \tau_i \}$
we may assume each of the open sets $U_i$ is affine.
We may also assume that $\bsym{U}$ is finite, say
$\bsym{U} = \{ U_0, \ldots, U_m \}$ .
For any pair of indices $(i, j)$ let 
$\rho_{(i, j)} := \tau_j^{-1} \circ \tau_i$, 
which is a gauge equivalence of
$\mcal{O}_{U_i \cap U_j}[[\hbar]]$. 
We are going to construct a gauge equivalence
$\rho_i$ of $\mcal{O}_{U_i}[[\hbar]]$, for every
$i \in \{ 0, \ldots, m \}$, such that 
$\rho_{(i, j)} = \rho_j \circ \rho_i^{-1}$.
Then the new differential structure $\bsym{\tau}'$, defined
by $\tau'_i := \tau_i \circ \rho_i$, will satisfy
${\tau'}_j^{-1} \circ \tau'_i = 
\bsym{1}_{\mcal{O}_{U_i \cap U_j}[[\hbar]]}$,
the identity automorphism of 
$\mcal{O}_{U_i \cap U_j}[[\hbar]]$. 
Therefore the various $\tau_i'$ can be glued 
to a global differential 
trivialization $\tau' : \mcal{O}_X [[\hbar]] \iso \mcal{A}$ 
as required.

Let $\mcal{D}^{\mrm{nor}}_X$ be the subsheaf of 
$\mcal{D}_X$ consisting 
of operators that vanish on $1_{\mcal{O}_X}$. This is the left 
ideal of $\mcal{D}_X$ generated by the sheaf of derivation
$\mcal{T}_X$. There is a direct sum decomposition
$\mcal{D}_X = \mcal{O}_X \oplus \mcal{D}^{\mrm{nor}}_X$
(as sheaves of left $\mcal{O}_X$-modules), and therefore
$\mrm{H}^1(X, \mcal{D}^{\mrm{nor}}_X) = 0$. 

Consider the sheaf of nonabelian groups $G$ on $X$ 
whose sections on an open set $U$ is the group of gauge 
equivalences of $\mcal{O}_{U}[[\hbar]]$. Let
$\mcal{D}^{\mrm{nor}}_X[[\hbar]]^+ :=
\hbar \mcal{D}^{\mrm{nor}}_X[[\hbar]]$.
As sheaves of sets there is a canonical isomorphism
$\mcal{D}^{\mrm{nor}}_X[[\hbar]]^+ \iso G$,
whose formula is
$\sum_{k = 1}^{\infty} D_k \hbar^k \mapsto 
\bsym{1}_{\mcal{O}_X} + \sum_{k = 1}^{\infty} D_k \hbar^k$. 
Define $G^k$ to be the subgroup of $G$ 
consisting of all equivalences 
congruent to $\bsym{1}_{\mcal{O}_{X}}$ modulo $\hbar^{k+1}$. 
Then each $G^k$ is a normal subgroup, and the map
$\mcal{D}^{\mrm{nor}}_X \iso G^k / G^{k+1}$,
$D \mapsto \bsym{1}_{\mcal{O}_X} + D \hbar^{k+1}$,
is an isomorphism of sheaves of abelian groups. 
Moreover the conjugation 
action of $G$ on $G^k / G^{k+1}$ is trivial, so that
$\gamma_1 \gamma_2 = \gamma_2 \gamma_1 \in 
G^k / G^{k+1}$
for every $\gamma_1 \in G^k$ and $\gamma_2 \in G$. 

The gauge equivalences
$\rho_i = \sum_{k = 0}^{\infty} \rho_{i, k} \hbar^k$
will be defined by successive approximations; namely the 
differential operators 
$\rho_{i, k} \in \Gamma(U_i, \mcal{D}^{\mrm{nor}}_X)$ 
shall be defined by recursion on $k$, simultaneously for all 
$i \in \{ 0, \ldots, m \}$. For $k = 0$ we take 
$\rho_{i, 0} := \bsym{1}_{\mcal{O}_{U_i}}$ of course.
Now assume at the $k$-th stage we have operators
$\rho_i^{(k)} := \sum_{l = 0}^k \rho_{i, l} \hbar^l$
which satisfy
\[ \rho_j^{(k)} \circ (\rho_i^{(k)})^{-1} \equiv \rho_{(i, j)}
\opn{mod} \hbar^{k + 1} . \]
This means that 
\[ \rho_j^{(k)} \circ (\rho_i^{(k)})^{-1} \circ 
(\rho_{(i, j)})^{-1} \in G^k . \]
By the properties of the group $G$ mentioned above the
function
$\{ 0, \ldots, m \}^2 \to G^k / G^{k+1}$,
\[ (i, j) \mapsto \rho_j^{(k)} \circ (\rho_i^{(k)})^{-1} \circ 
(\rho_{(i, j)})^{-1} \in 
\Gamma(U_i \cap U_j, G^k / G^{k+1}), \]
is a \v{C}ech $1$-cocycle for the affine covering $\bsym{U}$. 
Since $G^k / G^{k+1} \cong \mcal{D}^{\mrm{nor}}_X$, 
and we are given that 
$\mrm{H}^1(X, \mcal{D}^{\mrm{nor}}_X) = 0$, 
it follows that there exists a $0$-cochain 
$i \mapsto \rho_{i, k+1} \in \Gamma(U_i, \mcal{D}^{\mrm{nor}}_X)$
such that 
\[ (\bsym{1} + \rho_{j, k+1} \hbar^{k+1}) \circ 
(\bsym{1} + \rho_{i, k+1} \hbar^{k+1})^{-1} \equiv 
\rho_j^{(k)} \circ (\rho_i^{(k)})^{-1} \circ 
(\rho_{(i, j)})^{-1} \opn{mod} \hbar^{k+2} . \]
\end{proof}

\begin{prop} \label{prop6.1}
Let $\star$ and $\star'$ be two star products on 
$\mcal{O}_{X}[[\hbar]]$. Consider the globally trivialized
deformation quantizations 
$(\mcal{A}, \psi, \bsym{\tau})$ and 
$(\mcal{A}', \psi, \bsym{\tau})$, 
where
$\mcal{A} := (\mcal{O}_{X}[[\hbar]], \star)$,
$\mcal{A}' := (\mcal{O}_{X}[[\hbar]], \star')$,
$\psi := \bsym{1}_{\mcal{O}_{X}}$ and
$\bsym{\tau} := \{ \bsym{1}_{\mcal{O}_{X}[[\hbar]]} \}$. 
Then the deformation quantizations 
$(\mcal{A}, \psi, \bsym{\tau})$ and 
$(\mcal{A}', \psi, \bsym{\tau})$ 
are gauge equivalent, in the sense of Definition 
\tup{\ref{dfn4.6}}, iff there exists a gauge equivalence
$\gamma$ of $\mcal{O}_{X}[[\hbar]]$, in the sense of Definition 
\tup{\ref{dfn4.2}}, such that
\[ f \star' g = \gamma^{-1} \big( \gamma(f) \star \gamma(g)
\big) \]
for all local sections $f, g \in \mcal{O}_{X}$. 
\end{prop}

We leave out the easy proof.

\section{Review of Dir-Inv Modules}
\label{sec2}

In this section we review the concept of dir-inv 
structure, which was introduced in \cite[Section 1]{Ye2}. 
A dir-inv structure is a generalization of adic topology, and it
will turn out to be extremely useful in several places in the paper. 

Let $C$ be a commutative $\K$-algebra. 
We denote by $\cat{Mod} C$ the category of $C$-modules. 

\begin{dfn}
\begin{enumerate}
\item Let $M \in \cat{Mod} C$. An {\em inv module 
structure} on $M$ is an inverse system 
$\{ \mrm{F}^i M \}_{i \in \mbb{N}}$ of $C$-submodules
of $M$. The pair 
$(M, \{ \mrm{F}^i M \}_{i \in \mbb{N}})$
is called an {\em inv $C$-module}.
\item Let $(M, \{ \mrm{F}^i M \}_{i \in \mbb{N}})$ and
$(N, \{ \mrm{F}^i N \}_{i \in \mbb{N}})$ be two inv 
$C$-modules.
A function $\phi : M \to N$ 
($C$-linear or not) is said to be {\em continuous}
if for every $i \in \mbb{N}$ there exists $i' \in \mbb{N}$ such 
that $\phi(\mrm{F}^{i'} M) \subset \mrm{F}^i N$.
\item Define $\cat{Inv} \cat{Mod} C$
to be the category whose objects are the inv $C$-modules, 
and whose morphisms are the continuous $C$-linear 
homomorphisms.  
\end{enumerate}
\end{dfn}

There is a full and faithful embedding of categories
$\cat{Mod} C \inj \cat{Inv} \cat{Mod} C$,
$M \mapsto (M, \{ \ldots, 0, 0 \})$.

Recall that a directed set is a partially ordered set $J$ 
with the property that for any $j_1, j_2 \in J$ there exists 
$j_3 \in J$ such that $j_1, j_2 \leq j_3$.

\begin{dfn}
\begin{enumerate}
\item Let $M \in \cat{Mod} C$. A {\em 
dir-inv module structure} on $M$ is a direct system 
$\{ \mrm{F}_j M \}_{j \in J}$ of $C$-submodules 
of $M$, indexed by a nonempty directed set $J$, 
together with an inv module structure on 
each $\mrm{F}_j M$, such that for every $j_1 \leq j_2$ the 
inclusion $\mrm{F}_{j_1} M \inj \mrm{F}_{j_2} M$
is continuous. The pair 
$(M, \{ \mrm{F}_j M \}_{j \in J})$
is called a {\em dir-inv $C$-module}.
\item Let $(M, \{ \mrm{F}_j M \})_{j \in J}$ and
$(N, \{ \mrm{F}_k N \}_{k \in K})$ 
be two dir-inv $C$-modules.
A function $\phi : M \to N$ 
($C$-linear or not)  is said to be {\em continuous}
if for every $j \in J$ there exists $k \in K$ 
such that $\phi(\mrm{F}_j M) \subset \mrm{F}_k N$, and
$\phi : \mrm{F}_j M \to \mrm{F}_k N$ is a continuous 
function between these two inv $C$-modules.
\item Define $\cat{Dir} \cat{Inv} \cat{Mod} C$
to be the category whose objects are the dir-inv $C$-modules, 
and whose morphisms are the continuous $C$-linear 
homomorphisms.  
\end{enumerate}
\end{dfn}

An inv $C$-module $M$ can be endowed with a
dir-inv module structure $\{ \mrm{F}_j M \}_{j \in J}$, where 
$J := \{ 0 \}$ and $\mrm{F}_0 M := M$.
Thus we get a full and faithful embedding 
$\cat{Inv} \cat{Mod} C \inj 
\cat{Dir} \cat{Inv} \cat{Mod} C$. 

Inv modules and dir-inv modules come in a few ``flavors'': 
trivial, discrete and complete. 
A {\em discrete inv module} is one which is isomorphic, in
$\cat{Inv} \cat{Mod} C$, to an object of 
$\cat{Mod} C$ (via the canonical embedding above). 
A {\em complete inv module} is an inv module 
$(M, \{ \mrm{F}^i M \}_{i \in \mbb{N}})$
such that the canonical map 
$M \to \lim_{\leftarrow i} M / \mrm{F}^i M$
is bijective. A {\em discrete} (resp.\ {\em complete})
{\em dir-inv module} 
is one which is isomorphic, in 
$\cat{Dir} \cat{Inv} \cat{Mod} C$, to a dir-inv module
$(M, \{ \mrm{F}_j M \}_{j \in J})$, where all the inv 
modules $\mrm{F}_j M$ are discrete (resp.\ complete), 
and the canonical map
$\lim_{j \to} \mrm{F}_j M \to M$
in $\cat{Mod} C$ is bijective. A {\em trivial dir-inv 
module} is one which is isomorphic to an object of 
$\cat{Mod} C$. Discrete dir-inv modules are complete, but 
there are also other complete modules, as the next example shows.

\begin{exa} \label{exa2.4}
Assume $C$ is noetherian and $\mfrak{c}$-adically complete 
for some ideal $\mfrak{c}$. Let $M$ be a finitely 
generated $C$-module, and define 
$\mrm{F}^i M := \mfrak{c}^{i+1} M$. Then 
$\{ \mrm{F}^i M \}_{i \in \mbb{N}}$ is called the 
{\em $\mfrak{c}$-adic inv structure}, and of course 
$(M, \{ \mrm{F}^i M \}_{i \in \mbb{N}})$ is a complete inv module.
Next consider an arbitrary $C$-module $M$. We take
$\{ \mrm{F}_j M \}_{j \in J}$ to be the collection of 
finitely generated $C$-submodules of $M$. 
This dir-inv module structure on $M$ 
is called the {\em $\mfrak{c}$-adic dir-inv structure}. 
Again $(M, \{ \mrm{F}_j M \}_{j \in J})$ is a complete 
dir-inv $C$-module. Note that a finitely generated
$C$-module $M$ is discrete as inv module iff 
$\mfrak{c}^{i} M = 0$ for $i \gg 0$; and a $C$-module
is discrete as dir-inv module iff it is a direct 
limit of discrete finitely generated modules.
\end{exa}

The category $\cat{Dir} \cat{Inv} \cat{Mod} C$ is additive.
Given a collection $\{ M_k \}_{k \in K}$ of dir-inv modules, the 
direct sum $\boplus_{k \in K} M_k$ has a dir-inv module structure, 
making it into the coproduct of $\{ M_k \}_{k \in K}$ in
$\cat{Dir} \cat{Inv} \cat{Mod} C$.
Note that if the index set
$K$ is infinite and each $M_k$ is a nonzero 
discrete inv module, then $\boplus_{k \in K} M_k$ is a 
discrete dir-inv module which is not trivial.
The tensor product 
$M \otimes_{C} N$
of two dir-inv modules is again a dir-inv module.
There is a completion functor
$M \mapsto \what{M}$. 
(Warning: if $M$ is complete then 
$\what{M} = M$, but
it is not known if $\what{M}$ is complete for arbitrary 
$M$.) The completed tensor product is  
$M \what{\otimes}_{C} N :=
\what{M \otimes_{C} N}$. 
Completion commutes with direct sums: if 
$M \cong \boplus_{k \in K} M_k$ then 
$\what{M} \cong  \boplus_{k \in K} \what{M}_k$.

A {\em graded dir-inv module} (or graded object in 
$\cat{Dir} \cat{Inv} \cat{Mod} C$)
is a direct sum $M = \boplus_{k \in \mbb{Z}} M_k$, 
where each $M_k$ is a dir-inv module. A {\em DG algebra}
in $\cat{Dir} \cat{Inv} \cat{Mod} C$ is a graded dir-inv 
module $A = \boplus_{k \in \mbb{Z}} A^k$, together 
with continuous $C$-(bi)linear functions 
$\mu : A \times A \to A$ and $\d : A \to A$, 
which make $A$ into a DG $C$-algebra.
If $A$ is a super-commutative associative unital DG algebra in 
$\cat{Dir} \cat{Inv} \cat{Mod} C$, and $\mfrak{g}$ is a DG Lie 
Algebra in $\cat{Dir} \cat{Inv} \cat{Mod} C$, then 
$A \hatotimes{C} \mfrak{g}$ is a DG Lie 
Algebra in $\cat{Dir} \cat{Inv} \cat{Mod} C$. 

Let $A$ be a super-commutative associative unital DG algebra in 
$\cat{Dir} \cat{Inv} \cat{Mod} C$. A {\em DG $A$-module} 
in $\cat{Dir} \cat{Inv} \cat{Mod} C$ is a graded object
$M$ in $\cat{Dir} \cat{Inv} \cat{Mod} C$,
together with continuous $C$-(bi)linear functions 
$\mu : A \times M \to M$ and $\d : M \to M$, 
which make $M$ into a DG $A$-module in the usual sense. 
A {\em DG $A$-module Lie algebra} in
$\cat{Dir} \cat{Inv} \cat{Mod} C$ is a DG Lie algebra
$\mfrak{g}$ in $\cat{Dir} \cat{Inv} \cat{Mod} C$,
together with a continuous $C$-bilinear function
$\mu : A \times \mfrak{g} \to \mfrak{g}$, 
such that $\mfrak{g}$ becomes a DG $A$-module, and 
\[ [a_1 \gamma_1, a_2 \gamma_2] = (-1)^{i_2 j_1}
a_1 a_2 \,[\gamma_1, \gamma_2] \]
for all $a_k \in A^{i_k}$ and 
$\gamma_k \in \mfrak{g}^{j_k}$.

All the constructions above can be geometrized. 
Let $(Y, \mcal{O})$ be a commutative ringed space
over $\K$, i.e.\ $Y$ is a topological space, and $\mcal{O}$ is 
a sheaf of commutative $\K$-algebras on $Y$. 
We denote by $\cat{Mod} \mcal{O}$ the category of 
$\mcal{O}$-modules on $Y$. Then we can talk about the category 
$\cat{Dir} \cat{Inv} \cat{Mod} \mcal{O}$ of dir-inv 
$\mcal{O}$-modules. 

\begin{exa} \label{exa2.1}
Geometrizing Example \ref{exa2.4}, let $\mfrak{X}$ be a noetherian 
formal scheme, with defining ideal $\mcal{I}$. Then any coherent 
$\mcal{O}_{\mfrak{X}}$-module $\mcal{M}$ is an 
inv $\mcal{O}_{\mfrak{X}}$-module, with system of submodules
$\{ \mcal{I}^{i+1} \mcal{M} \}_{i \in \mbb{N}}$,
and $\mcal{M} \cong \what{\mcal{M}}$; cf.\ 
\cite{EGA-I}. We call an $\mcal{O}_{\mfrak{X}}$-module
{\em dir-coherent} if it is the direct limit of coherent 
$\mcal{O}_{\mfrak{X}}$-modules. Any dir-coherent module 
is quasi-coherent, but it is not known if the converse is true.
At any rate, a dir-coherent $\mcal{O}_{\mfrak{X}}$-module 
$\mcal{M}$ is a dir-inv $\mcal{O}_{\mfrak{X}}$-module, where we 
take $\{ \mrm{F}_j \mcal{M} \}_{j \in J}$ to be the collection of 
coherent submodules of $\mcal{M}$. Any 
dir-coherent $\mcal{O}_{\mfrak{X}}$-module is then a complete 
dir-inv module. This dir-inv module structure on $\mcal{M}$ 
is called the {\em $\mcal{I}$-adic dir-inv structure}. 
\end{exa}

If $f : (Y', \mcal{O}') \to (Y, \mcal{O})$ is a morphism of 
ringed spaces and
$\mcal{M} \in \cat{Dir} \cat{Inv} \cat{Mod} \mcal{O}$,
then there is an obvious structure of dir-inv 
$\mcal{O}'$-module on $f^{*} \mcal{M}$, and we define
$f^{\what{*}} \mcal{M} := \what{f^* \mcal{M}}$. 
If $\mcal{M}$ is a graded object in
$\cat{Dir} \cat{Inv} \cat{Mod} \mcal{O}$,
then the inverse images $f^* \mcal{M}$
and $f^{\what{*}} \mcal{M}$ are graded objects in 
$\cat{Dir} \cat{Inv} \cat{Mod} \mcal{O}'$.
If $\mcal{G}$ is an algebra (resp.\ a DG algebra) in 
$\cat{Dir} \cat{Inv} \cat{Mod} \mcal{O}$, then 
$f^* \mcal{G}$ and  $f^{\what{*}}\, \mcal{G}$ are
algebras (resp.\ DG algebras) in 
$\cat{Dir} \cat{Inv} \cat{Mod} \mcal{O}'$. 
Given $\mcal{N} \in \cat{Dir} \cat{Inv} \cat{Mod} \mcal{O}'$
there is an obvious 
dir-inv $\mcal{O}$-module structure on $f_* \mcal{N}$.

\section{Universal Formulas for Deformation Quantization}
\label{sec3}

In this section, as before, $\K$ is a field of characteristic $0$.

 From here to Corollary \ref{cor7.1} we consider the following 
data. Let $\mfrak{g} = \boplus_{j \in \mbb{Z}} \mfrak{g}^j$ 
be a DG Lie algebra over $\K$.
We put on each $\mfrak{g}^j$ the discrete inv $\K$-module structure,
and $\mfrak{g}$ is given the $\boplus$ dir-inv structure; so
$\mfrak{g}$ is a discrete, but possibly nontrivial, DG Lie algebra 
in $\cat{Dir} \cat{Inv} \cat{Mod} \K$. Let 
$A$ be noetherian commutative complete local  
$\K$-algebra with maximal ideal $\mfrak{m}$. 
We put on $A$ and $\mfrak{m}$ the $\mfrak{m}$-adic inv structures. 
For $i \geq 0$ let $A_i := A / \mfrak{m}^{i+1}$, 
which is an artinian local algebra with maximal ideal
$\mfrak{m}_i := \mfrak{m} / \mfrak{m}^{i+1}$; so 
$A_i$ and $\mfrak{m}_i$ are discrete inv modules.
We obtain a new DG Lie algebra
$A \hatotimes{\K} \mfrak{g} = 
A \hatotimes{} \mfrak{g} = \boplus_{j \in \mbb{Z}}\, 
A \hatotimes{} \mfrak{g}^j$, 
and there are related DG Lie algebras 
$\mfrak{m} \hatotimes{} \mfrak{g} \subset
A \hatotimes{} \mfrak{g}$
and
$\mfrak{m}_i \otimes \mfrak{g} \subset
A_i \otimes \mfrak{g}$.
Note that for every $j$ one has 
$A \hatotimes{} \mfrak{g}^j \cong 
\lim_{\leftarrow i} (A_i \otimes \mfrak{g}^j)$
in $\cat{Inv} \cat{Mod} \K$.
In case $A = \K[[\hbar]]$ we shall also use the notation
$\mfrak{g} [[\hbar]]^+ := \mfrak{m} \, \what{\otimes} \, 
\mfrak{g}$,
namely 
$\mfrak{g} [[\hbar]]^+ = 
\boplus_j \, \hbar \mfrak{g}^j [[\hbar]]$.

Recall the correspondence between finite dimensional nilpotent Lie 
algebras and unipotent algebraic groups over $\K$ 
(see \cite[Theorem XVI.4.2]{Ho}). Given a nilpotent Lie algebra 
$\mfrak{h}$ we denote by $\opn{exp}(\mfrak{h})$ the corresponding 
group. This group has the same underlying scheme structure as 
$\mfrak{h}$, and the product is according to the Campbell-Hausdorff 
formula. The assignment $\mfrak{h} \mapsto \opn{exp}(\mfrak{h})$
is functorial. 

For any $i$ the Lie algebra 
$\mfrak{m}_i \otimes \mfrak{g}^0$
is a nilpotent, and in fact it is a direct limit of finite 
dimensional nilpotent Lie algebras. Therefore we obtain a group
$\opn{exp}(\mfrak{m}_i \otimes \mfrak{g}^0)$, 
which is a direct limit of unipotent groups.
Passing to the inverse limit in $i$ we get a group
$\opn{exp}(\mfrak{m} \otimes \mfrak{g}^0) := 
\lim_{\leftarrow i} \opn{exp}(\mfrak{m}_i \otimes \mfrak{g}^0)$.

Given a vector space $V$ over $\K$ let 
$\opn{Aff}(V) := \opn{GL}(V) \ltimes V$,
the group of affine transformations. Its Lie algebra is
$\mfrak{aff}(V) := \mfrak{gl}(V) \ltimes V$. 
If $V$ is finite dimensional then of course
$\opn{Aff}(V)$ is an algebraic group; but we will be interested in 
$V := \mfrak{m} \, \what{\otimes} \, \mfrak{g}^1$.

For 
$\gamma \in \mfrak{m} \, \what{\otimes} \, \mfrak{g}^0$ and
$\omega \in \mfrak{m} \, \what{\otimes} \, \mfrak{g}^1$
define
\[ \opn{af}(\gamma)(\omega) := [\gamma, \omega] - \mrm{d}(\gamma) 
 = (\opn{ad}(\gamma) -\d)(\omega) \in
\mfrak{m} \hatotimes{} \mfrak{g}^1  , \]
where $\d$ and $[-,-]$ are the operations of the DG Lie algebra 
$\mfrak{m} \, \what{\otimes} \, \mfrak{g}$.
A calculation shows that this is a homomorphism of Lie 
algebras
\[ \opn{af} : \mfrak{m} \, \what{\otimes} \, \mfrak{g}^0
\to 
\mfrak{aff}(\mfrak{m} \, \what{\otimes} \, \mfrak{g}^1)
. \]

Recall that the {\em Maurer-Cartan equation} in 
$A \hatotimes{} \mfrak{g}$ is
\begin{equation} \label{eqn3.6}
\d(\omega) + \smfrac{1}{2} [\omega, \omega] = 0 
\end{equation}
for $\omega \in A \hatotimes{} \mfrak{g}^1$.

\begin{lem} \label{lem7.3}
\begin{enumerate}
\item The Lie algebra homomorphism $\opn{af}$ integrates to a 
group homomorphism
\[ \opn{exp}(\opn{af}) : \opn{exp}(\mfrak{m} 
\, \what{\otimes} \, \mfrak{g}^0)
\to \opn{Aff}(\mfrak{m} \, \what{\otimes} \, \mfrak{g}^1) . \]
\item Assume 
$\omega \in \mfrak{m} \, \what{\otimes} \, \mfrak{g}^1$ is a 
solution of the MC equation in 
$\mfrak{m} \, \what{\otimes} \, \mfrak{g}$,
and let $\gamma \in \mfrak{m} \, \what{\otimes} \, \mfrak{g}^0$. 
Then 
$\opn{exp}(\opn{af})(\opn{exp}(\gamma))(\omega)$ 
is also a solution of the MC equation. 
\end{enumerate}
\end{lem}

\begin{proof}
We may assume that 
$\mfrak{g} = \boplus_{j \geq 0} \, \mfrak{g}^{j}$.
First consider the nilpotent case. The DG Lie algebra 
$\mfrak{m}_i \otimes \mfrak{g}$
is the direct limit of sub DG Lie algebras 
$\mfrak{h} = \boplus_{j \geq 0} \, \mfrak{h}^{j}$, 
which are nilpotent, and each $\mfrak{h}^{j}$ is a finite 
dimensional vector space. The arguments of \cite[Section 1.3]{GM}  
apply here, so we obtain a homomorphism of algebraic groups
$\opn{exp}(\opn{af}) : \opn{exp}(\mfrak{h}^0)
\to \opn{Aff}(\mfrak{h}^1)$,
and $\opn{exp}(\mfrak{h}^0)$ preserves the set of solutions of the 
MC equation in $\mfrak{h}^1$. Passing to the direct limit over 
these subalgebras we get a homomorphism of groups
$\opn{exp}(\opn{af}) : \opn{exp}(\mfrak{m}_i \otimes \mfrak{g}^0)
\to \opn{Aff}(\mfrak{m}_i \otimes \mfrak{g}^1)$,
and $\opn{exp}(\mfrak{m}_i \otimes \mfrak{g}^0)$ 
preserves the set of solutions of the 
MC equation in $\mfrak{m}_i \otimes \mfrak{g}^1$. 
Finally we pass to the inverse limit in $i$.
\end{proof}

The formula for  
$\opn{exp}(\opn{af})(\opn{exp}(\gamma))(\omega)$ is, according to 
\cite{GM}:
\begin{equation} \label{eqn8.2}
\opn{exp}(\opn{af})(\opn{exp}(\gamma))(\omega) = 
\opn{exp}(\opn{ad}(\gamma))(\omega) +
\frac{1 - \opn{exp}(\opn{ad}(\gamma))}{\opn{ad}(\gamma)} 
(\d(\gamma)) .
\end{equation}
On the right side of the equation
``$\opn{exp}$'' stands for the usual exponential 
power series
$\opn{exp}(t) = \sum_{k=0}^{\infty} \smfrac{1}{k!} t^k$,
and this makes sense because 
$\lim_{k \to \infty} \opn{ad}(\gamma)^k(\omega) = 0$
in the $\mfrak{m}$-adic inv structure on 
$\mfrak{m} \hatotimes{} \mfrak{g}^1$.

\begin{dfn} \label{dfn7.1}
Elements of the group 
$\opn{exp}(\mfrak{m} \, \what{\otimes} \, \mfrak{g}^0)$
are called {\em gauge equivalences}. We write
\[ \opn{MC}(\mfrak{m} \, \what{\otimes} \, \mfrak{g}) :=
\frac{ \{ \text{solutions of the MC equation in }
\mfrak{m} \, \what{\otimes} \, \mfrak{g} \} }{
\{ \text{gauge equivalences} \} } .  \]
\end{dfn}

\begin{lem} \label{lem7.1}
The canonical projection
\[ \opn{MC}(\mfrak{m} \, \what{\otimes} \, \mfrak{g}) \to
\lim_{\leftarrow i} 
\opn{MC}(\mfrak{m}_i \otimes \mfrak{g}) \]
is bijective.
\end{lem}

The easy proof is omitted. 

\begin{rem}
Consider the super-commutative DG algebra
$\Omega_{\K[t]} = \Omega^0_{\K[t]} \oplus \Omega^1_{\K[t]}$,
where $\K[t]$ is the polynomial algebra in the variable
$t$. There is an induced DG Lie algebra 
$\Omega_{\K[t]} \otimes \mfrak{g}$. For any $\lambda \in \K$ there 
is a DG Lie algebra homomorphism
$\Omega_{\K[t]} \otimes \mfrak{g} \to \mfrak{g}$,
$t \mapsto \lambda$.
Assume $A$ is artinian, and let $\omega_0$ and $\omega_1$ be two 
solutions of the MC equation in $\mfrak{m} \otimes \mfrak{g}$.
According to \cite[4.5.2(3)]{Ko1} the following conditions are 
equivalent:
\begin{enumerate}
\rmitem{i} $\omega_0$ and $\omega_1$ are gauge equivalent, in the 
sense of Definition \ref{dfn7.1}.
\rmitem{ii} There is a solution $\omega(t)$
of the MC equation in the DG Lie algebra 
$\Omega_{\K[t]} \otimes \mfrak{m} \otimes \mfrak{g}$,
such that for $i \in \{ 0, 1 \} \subset \K$,
the specialization homomorphisms
$\Omega_{\K[t]} \otimes \mfrak{m} \otimes \mfrak{g} \to
\mfrak{m} \otimes \mfrak{g}$,
$t \mapsto i$, send $\omega(t) \mapsto \omega_i$. 
\end{enumerate}
See also \cite{Fu} and \cite{Hi1}. We will not need these facts in 
our paper.
\end{rem}
 
For a graded $\K$-module $M$ the expression 
$\bwedge^i M$ denotes the $i$-th super-exterior power. 

\begin{dfn} \label{dfn6.2}
Let $\mfrak{g}$ and $\mfrak{h}$ be two DG Lie algebras. An 
{\em $\mrm{L}_{\infty}$ morphism} 
$\Psi : \mfrak{g} \to \mfrak{h}$
is a collection $\Psi = \{ \psi_i \}_{i \geq 1}$
of $\K$-linear homomorphisms
$\psi_i : \bwedge^i \mfrak{g} \to \mfrak{h}$,
each of them homogeneous of degree $1-i$, satisfying
\[ \begin{aligned}
& \mrm{d} \bigl( \psi_i(\gamma_1 \wedge \ldots 
\wedge \gamma_i) \bigr) - \sum_{k = 1}^i \pm 
\psi_i \bigl( \gamma_1 \wedge \ldots \wedge 
\mrm{d}(\gamma_k) \wedge \ldots \wedge \gamma_i \bigr) = \\
& \quad \smfrac{1}{2} \sum_{\substack{k, l \geq 1 \\ k + l = i}}
\smfrac{1}{k! \, l!} \sum_{\sigma \in \Sigma_i} \pm
\bigl[ \psi_k (\gamma_{\sigma(1)} \wedge \ldots \wedge 
\gamma_{\sigma(k)}),
\psi_l (\gamma_{\sigma(k + 1)} \wedge \ldots \wedge 
\gamma_{\sigma(k + l)}) \bigr] \\
& \quad + \sum_{k < l} \pm
\psi_{i-1} \bigl( [\gamma_k, \gamma_l] \wedge
\gamma_{1} \wedge \ldots \wedge \gamma_{i} \bigr) .
\end{aligned} \]
Here $\gamma_k \in \mfrak{g}$ are homogeneous elements,
$\Sigma_i$ is the permutation group of $\{ 1, \ldots, i \}$,
and the signs depend only on the indices, the permutations and the 
degrees of the elements $\gamma_k$. See \cite[Section 6]{Ke} or
\cite[Theorem 3.1]{CFT} for the explicit signs.
\end{dfn}

An $\mrm{L}_{\infty}$ morphism is a generalization of a DG Lie 
algebra homomorphism. Indeed, 
$\psi_1 : \mfrak{g} \to \mfrak{h}$ is a homomorphism of complexes 
of $\K$-modules, and 
$\mrm{H}(\psi_1) : \mrm{H} \mfrak{g} \to \mrm{H} \mfrak{h}$ 
is a homomorphism of graded Lie algebras.

Suppose 
$\Psi = \{ \psi_i \}_{i \geq 1} : \mfrak{g} \to \mfrak{h}$
is an $\mrm{L}_{\infty}$ morphism.
For every $i$ we can extend the $\K$-multilinear function
$\psi_i : \prod^i \mfrak{g} \to \mfrak{h}$ uniquely
to a continuous $A$-multilinear function
$\psi_{A, i} : 
\prod^i (A \hatotimes{} \mfrak{g})
\to A \hatotimes{} \mfrak{h}$.
These restrict to functions 
$\psi_{A, i} : 
\prod^i (\mfrak{m} \hatotimes{} \mfrak{g})
\to \mfrak{m} \hatotimes{} \mfrak{h}$.
Clearly 
$\Psi_A = \{ \psi_{A, i} \} : 
\mfrak{m} \, \what{\otimes} \, \mfrak{g} 
\to \mfrak{m} \, \what{\otimes} \, \mfrak{h}$
is an $\mrm{L}_{\infty}$ morphism; we call it the continuous 
$A$-multilinear extension of $\Psi$.

\begin{thm}[{\cite[Section 4.4]{Ko1}, \cite[Theorem 2.2.2]{Fu}}] 
\label{thm7.1}
Assume $A$ is artinian. Let $\Psi : \mfrak{g} \to \mfrak{h}$ be 
an $\mrm{L}_{\infty}$ quasi-isomorphism. Then the function
\begin{equation} \label{eqn7.1}
\omega \mapsto \sum_{j \geq 1} \smfrac{1}{j!}
\psi_{A, j} (\omega^j) 
\end{equation}
induces a bijection
\[ \opn{MC}(\Psi_A) :
\opn{MC}(\mfrak{m} \otimes \mfrak{g}) \iso
\opn{MC}(\mfrak{m} \otimes \mfrak{h}) . \]
\end{thm}

\begin{cor} \label{cor7.1}
Let $(A, \mfrak{m})$ be a complete noetherian local 
$\K$-algebra, and let $\Psi : \mfrak{g} \to \mfrak{h}$ be 
an $\mrm{L}_{\infty}$ quasi-isomorphism between two discrete DG 
Lie algebras. 
Then the function \tup{(\ref{eqn7.1})} induces a bijection
\[ \opn{MC}(\Psi_A) :
\opn{MC}(\mfrak{m} \, \what{\otimes} \, \mfrak{g}) \iso
\opn{MC}(\mfrak{m} \, \what{\otimes} \, \mfrak{h}) . \]
\end{cor}

\begin{proof}
Use Theorem \ref{thm7.1} and Lemma \ref{lem7.1}.
\end{proof}

Let $C$ be a commutative $\K$-algebra. The module of derivations 
of $C$ relative to $\K$ is denoted by 
$\mcal{T}_{C / \K} = \mcal{T}_{C}$. For $p \geq -1$ let
$\mcal{T}^p_{\mrm{poly}}(C) := \bwedge^{p+1}_{C} \mcal{T}_{C}$,
the $p$-th exterior power. The direct sum
$\mcal{T}_{\mrm{poly}}(C) := \boplus_{p \geq -1} 
\mcal{T}^p_{\mrm{poly}}(C)$
is a DG Lie algebra over $\K$ 
with trivial differential and with the 
Schouten-Nijenhuis Lie bracket (see \cite{Ko2} for details). 

For any $p \geq 0$ and $m \geq 0$ let 
$\mrm{F}_m \mcal{D}^p_{\mrm{poly}}(C)$ be the set of 
$\K$-multilinear functions $\phi : C^{p+1} \to C$ that are 
differential operators of order $\leq m$ in each argument
(in the sense of \cite{EGA-IV}). For 
$p = -1$ let 
$\mrm{F}_m \mcal{D}^{-1}_{\mrm{poly}}(C) := C$. 
Define
$\mcal{D}^p_{\mrm{poly}}(C) := \bigcup_{m \geq 0}
\mrm{F}_m \mcal{D}^p_{\mrm{poly}}(C)$
and
$\mcal{D}_{\mrm{poly}}(C) := \boplus_{p \geq -1} 
\mcal{D}^p_{\mrm{poly}}(C)$.
This is a sub DG Lie algebra of the shifted Hochschild cochain 
complex of $C$, with shifted Hochschild differential and 
Gerstenhaber Lie bracket (see \cite{Ko1}). 
We view $\mcal{D}_{\mrm{poly}}(C)$ as a 
left $C$-module by the rule
$(c \cdot \phi)(c_1, \ldots, c_{p+1}) :=
c \cdot \phi(c_1, \ldots, c_{p+1})$.

For $p \geq 0$ define
$\mcal{D}^{\mrm{nor}, p}_{\mrm{poly}}(C)$
to be the subset of $\mcal{D}^{p}_{\mrm{poly}}(C)$
consisting of the poly differential operators
$\phi : C^{p+1} \to C$ such that 
$\phi(c_1, \ldots, c_{p+1}) = 0$ if $c_i = 1$ for some $i$. 
For $p = -1$ let
$\mcal{D}^{\mrm{nor}, -1}_{\mrm{poly}}(C) := C$.
Then 
$\mcal{D}^{\mrm{nor}}_{\mrm{poly}}(C) := 
\boplus_{p \geq -1} \mcal{D}^{\mrm{nor}, p}_{\mrm{poly}}(C)$
is a sub DG Lie algebra of 
$\mcal{D}_{\mrm{poly}}(C)$.

For any integer $p \geq 0$ 
there is a $C$-linear homomorphism
\[ \mcal{U}_1 : 
\mcal{T}^{p}_{\mrm{poly}}(C) \to 
\mcal{D}^{\mrm{nor}, p}_{\mrm{poly}}(C) \]
with formula 
\begin{equation} \label{eqn6.11}
\begin{aligned}
& \mcal{U}_1(\xi_1 \wedge \cdots \wedge \xi_{p+1})
(c_1, \ldots, c_{p+1}) := \\
& \qquad {\smfrac{1}{(p+1)!}} \sum_{\sigma \in \Sigma_{p+1}} 
\opn{sgn}(\sigma)
\xi_{\sigma(1)}(c_1) \cdots \xi_{\sigma({p+1})}(c_{p+1})
\end{aligned} 
\end{equation}
for elements
$\xi_1, \ldots, \xi_{p+1} \in \mcal{T}_C$
and
$c_1, \ldots, c_{p+1} \in C$. For $p = -1$ the map
$\mcal{U}_1 : \mcal{T}^{-1}_{\mrm{poly}}(C) \to 
\mcal{D}^{\mrm{nor}, -1}_{\mrm{poly}}(C)$
is the identity (of $C$).

The next result is a variant of the Hochschild-Kostant-Rosenberg 
Theorem. A slightly weaker result appeared in \cite{Ye1}. See 
\cite{Ko1} for the $\mrm{C}^{\infty}$ version. 

\begin{thm}[{\cite[Corollary 4.12]{Ye2}}] \label{thm3.6}
Suppose $C$ is a smooth $\K$-algebra. Then the homomorphism
$\mcal{U}_1 : \mcal{T}_{\mrm{poly}}(C) \to 
\mcal{D}^{\mrm{nor}}_{\mrm{poly}}(C)$
and the inclusion 
$\mcal{D}^{\mrm{nor}}_{\mrm{poly}}(C) \to \mcal{D}_{\mrm{poly}}(C)$
are both quasi-isomorphisms of complexes of $C$-modules. 
\end{thm}

Here is a slight modification of the celebrated result 
of Kontsevich, known as the {\em Kontsevich Formality Theorem}
\cite[Theorem 6.4]{Ko1}. In the form below it is 
is proved in \cite[Theorem 4.13]{Ye2}.

\begin{thm} \label{thm2.0}
Let $\K[\bsym{t}] = \K[t_1, \ldots, t_n]$ 
be the polynomial algebra in $n$ variables, and 
assume that $\mbb{R} \subset \K$. There is a collection of 
$\K$-linear homomorphisms 
\[ \mcal{U}_j :\, \bwedge^{j}
\mcal{T}_{\mrm{poly}}(\K[\bsym{t}]) \to
\mcal{D}_{\mrm{poly}}(\K[\bsym{t}]) , \]
indexed by $j \in \{ 1, 2, \ldots \}$, satisfying the 
following conditions.
\begin{enumerate}
\rmitem{i} The sequence $\mcal{U} = \{ \mcal{U}_j \}$ 
is an $\mrm{L}_{\infty}$-morphism
$\mcal{T}_{\mrm{poly}}(\K[\bsym{t}]) \to
\mcal{D}_{\mrm{poly}}(\K[\bsym{t}])$.
\rmitem{ii} Each $\mcal{U}_j$ is a poly differential operator 
of $\K[\bsym{t}]$-modules.
\rmitem{iii} Each $\mcal{U}_j$ is equivariant for the standard 
action of $\mrm{GL}_n(\K)$ on $\K[\bsym{t}]$.
\rmitem{iv} The homomorphism $\mcal{U}_1$ is given by equation 
\tup{(\ref{eqn6.11})}.
\rmitem{v} For any $j \geq 2$ and
$\alpha_1, \ldots, \alpha_j \in 
\mcal{T}^{0}_{\mrm{poly}}(\K[\bsym{t}])$
one has
$\mcal{U}_j(\alpha_1 \wedge \cdots \wedge \alpha_j) = 0$.
\rmitem{vi}  For any $j \geq 2$,
$\alpha_1 \in \mfrak{gl}_n(\K) \subset 
\mcal{T}^{0}_{\mrm{poly}}(\K[\bsym{t}])$
and
$\alpha_2, \ldots, \alpha_j \in 
\mcal{T}_{\mrm{poly}}(\K[\bsym{t}])$
one has
$\mcal{U}_j(\alpha_1 \wedge \cdots \wedge \alpha_j) = 0$.
\end{enumerate}
\end{thm}

\begin{rem}
Presumably the image of $\mcal{U}_j$ is inside 
$\mcal{D}^{\mrm{nor}}_{\mrm{poly}}(\K[\bsym{t}])$
for all $j$. However we did not verify this.
\end{rem}

\begin{rem}
The methods of Tamarkin \cite{Ta, Hi2}, or suitable arithmetic 
considerations \cite{Ko2}, should make it possible to 
extend Theorem \ref{thm2.0}, and hence all results of our paper, 
to any field $\K$ of characteristic $0$.
\end{rem}

Consider the power series algebra
$\K[[\bsym{t}]] = \K[[t_1, \ldots, t_n]]$.
As in Example \ref{exa2.1}, the $\K[[\bsym{t}]]$-modules
$\mcal{T}_{\mrm{poly}}(\K[[\bsym{t}]])$
and
$\mcal{D}_{\mrm{poly}}(\K[[\bsym{t}]])$
have the $\bsym{t}$-adic dir-inv structures. 
These are DG Lie algebras in $\cat{Dir} \cat{Inv} \cat{Mod} \K$.
Because $\K[\bsym{t}] \to \K[[\bsym{t}]]$ is flat and 
$\bsym{t}$-adically formally \'etale, it follows that 
there is an induced $\mrm{L}_{\infty}$ morphism
$\mcal{U} : \mcal{T}_{\mrm{poly}}(\K[[\bsym{t}]]) \to
\mcal{D}_{\mrm{poly}}(\K[[\bsym{t}]])$.
Since each 
\[ \mcal{U}_j : \sprod^j \, 
\mcal{T}_{\mrm{poly}}(\K[[\bsym{t}]]) \to
\mcal{D}_{\mrm{poly}}(\K[[\bsym{t}]]) \]
is a poly differential operator over 
$\K[[\bsym{t}]]$, it is continuous for the dir-inv 
structures. See \cite[Proposition 4.6]{Ye2} for 
details and proofs. 

Now suppose we are given a complete
super-commutative associative unital DG  
algebra $A = \boplus_{i \geq 0} A^i$ in 
$\cat{Dir} \cat{Inv} \cat{Mod} \K$. 
Let
\[ \mcal{U}_{A; j} : \sprod^j
\bigl( A \hatotimes{} \mcal{T}_{\mrm{poly}}(\K[[\bsym{t}]])
\bigr) \to
A \hatotimes{} \mcal{D}_{\mrm{poly}}(\K[[\bsym{t}]]) \]
be the continuous $A$-multilinear extension of 
$\mcal{U}_{j}$. It almost immediate from Definition \ref{dfn6.2} 
that  $\mcal{U}_{A} = \{ \mcal{U}_{A;j} \}_{j \geq 1}$
is an $\mrm{L}_{\infty}$ morphism; see 
\cite[Proposition 3.25]{Ye2}.

Let's recall the notion of twisting for a DG Lie algebra 
$\mfrak{g}$. Suppose $\omega \in \mfrak{g}^1$ is a solution of the 
MC equation (\ref{eqn3.6}). 
The twisted DG Lie algebra $\mfrak{g}_{\omega}$ is 
the same graded Lie algebra, but the new differential is
$\d_{\omega} := \d + \opn{ad}(\omega)$; i.e.\
$\d_{\omega}(\alpha) = \d(\alpha) + [\omega, \alpha]$. 

\begin{thm}[{\cite[Theorem 0.1]{Ye2}}] \label{thm3.3}
Assume $\mbb{R} \subset \K$. Let $A = \boplus_{i \geq 0} A^i$
be a complete super-commutative associative unital DG algebra in 
$\cat{Dir} \cat{Inv} \cat{Mod} \K$, and let
$\omega \in A^1 \hatotimes{}
\mcal{T}^0_{\mrm{poly}}(\K[[\bsym{t}]])$
be a solution of the Maurer-Cartan equation in 
$A \hatotimes{} \mcal{T}_{\mrm{poly}}(\K[[\bsym{t}]])$.
For any element
$\alpha \in \bwedge^j
\bigl( A \hatotimes{} \mcal{T}_{\mrm{poly}}(\K[[\bsym{t}]])
\bigr)$
define
\[ \mcal{U}_{A, \omega; j}(\alpha) :=
\sum_{k \geq 0} \smfrac{1}{(j+k)!} 
(\mcal{U}_{A; j+k})(\omega^k \wedge \alpha)
\in A \hatotimes{} \mcal{D}_{\mrm{poly}}(\K[[\bsym{t}]]) . \]
Let $\omega' := \mcal{U}_{A; 1}(\omega)$.
Then $\omega'$ is a solution of the Maurer-Cartan equation in
\linebreak
$A \hatotimes{} \mcal{D}_{\mrm{poly}}(\K[[\bsym{t}]])$,
and the sequence 
$\{ \mcal{U}_{A, \omega; j} \}_{j \geq 1}$
is a continuous $A$-multilinear
$\mrm{L}_{\infty}$ quasi-iso\-mor\-phism 
\[ \mcal{U}_{A, \omega} :
\bigl( A \hatotimes{} \mcal{T}_{\mrm{poly}}(\K[[\bsym{t}]])
\bigr)_{\omega} \to 
\bigl( A \hatotimes{} \mcal{D}_{\mrm{poly}}(\K[[\bsym{t}]])
\bigr)_{\omega'} . \]
\end{thm}

The sum occurring in the definition of 
$\mcal{U}_{A, \omega; j}(\alpha)$ is always
finite (but the number of 
nonzero terms depends on the argument $\alpha$). 

The group $\mrm{GL}_n(\K)$ acts on $\K[[\bsym{t}]]$ by linear 
change of coordinates. This is an action by $\K$-algebra 
automorphisms, and hence $\mrm{GL}_n(\K)$ acts on 
$\mcal{T}_{\mrm{poly}}(\K[[\bsym{t}]])$ and
$\mcal{D}_{\mrm{poly}}(\K[[\bsym{t}]])$
by continuous DG Lie algebra automorphisms. Suppose we are given 
an action of $\mrm{GL}_n(\K)$ on $A$ by continuous unital DG 
algebra automorphisms. Then we obtain an action of 
$\mrm{GL}_n(\K)$ on 
$A \hatotimes{} \mcal{T}_{\mrm{poly}}(\K[[\bsym{t}]])$ and
$A \hatotimes{} \mcal{D}_{\mrm{poly}}(\K[[\bsym{t}]])$
by continuous DG Lie algebra automorphisms. 

\begin{prop} \label{prop3.2}
Each operator $\mcal{U}_{A; j}$ is $\mrm{GL}_n(\K)$-equivariant, 
i.e.\ 
$\mcal{U}_{A; j}(g(\alpha))= \linebreak 
g(\mcal{U}_{A; j}(\alpha))$
for any $g \in \mrm{GL}_n(\K)$ and 
$\alpha \in \bwedge^j
\bigl( A \hatotimes{} \mcal{T}_{\mrm{poly}}(\K[[\bsym{t}]])
\bigr)$.
Moreover, if  $\omega$ is $\mrm{GL}_n(\K)$-invariant, then 
each operator $\mcal{U}_{A, \omega; j}$ is 
$\mrm{GL}_n(\K)$-equivariant. 
\end{prop}

\begin{proof}
Using continuity and multilinearity we may assume that 
\[ \alpha = (a_1 \otimes \alpha_1) \wedge \cdots \wedge
(a_j \otimes \alpha_j) ,
\]
with $a_k \in A$ and 
$\alpha_k \in \mcal{T}_{\mrm{poly}}(\K[[\bsym{t}]])$.
Then
\[ \begin{aligned}
g \bigl( \mcal{U}_{A; j}(\alpha) \bigr) & = 
\pm g \bigl( a_1 \cdots a_j \cdot \mcal{U}_{j}(\alpha_1 \wedge
\cdots \wedge \alpha_j) \bigr) \\
& = \pm g(a_1 \cdots a_j) \cdot 
g \bigl( \mcal{U}_{j}(\alpha_1 \wedge
\cdots \wedge \alpha_j) \bigr) \\
& =^{\diamondsuit} \pm g(a_1 \cdots a_j) \cdot 
\mcal{U}_{j} \bigl( g(\alpha_1 \wedge
\cdots \wedge \alpha_j) \bigr) \\
& = \mcal{U}_{A; j} \bigl( g(\alpha) \bigr) ,
\end{aligned} \]
where the equality marked $\diamondsuit$ is due to condition (iii) 
of Theorem \ref{thm2.0}.

We see that 
$g \bigl( \mcal{U}_{A, \omega; j}(\alpha) \bigr) = 
\mcal{U}_{A, g(\omega); j} \bigl( g(\alpha) \bigr)$.
Hence the second assertion. 
\end{proof}

Let $X$ be a smooth irreducible separated $n$-dimensional
$\K$-scheme. 

\begin{prop} \label{prop3.5}
There are sheaves of DG Lie algebras
$\mcal{T}_{\mrm{poly}, X}$,
$\mcal{D}_{\mrm{poly}, X}$
and
$\mcal{D}^{\mrm{nor}}_{\mrm{poly}, X}$
on $X$. As left $\mcal{O}_X$-modules all three are quasi-coherent. 
The sheaves $\mcal{D}_{\mrm{poly}, X}$
and
$\mcal{D}^{\mrm{nor}}_{\mrm{poly}, X}$
are quasi-coherent left $\mcal{D}_{X}$-modules.
For any affine open set $U = \opn{Spec} C \subset X$ one has
$\Gamma(U, \mcal{T}^{}_{\mrm{poly}, X}) = 
\mcal{T}^{}_{\mrm{poly}}(C)$,
$\Gamma(U, \mcal{D}^{}_{\mrm{poly}, X}) = 
\mcal{D}^{}_{\mrm{poly}}(C)$
and
$\Gamma(U, \mcal{D}^{\mrm{nor}}_{\mrm{poly}, X}) = 
\mcal{D}^{\mrm{nor}}_{\mrm{poly}}(C)$
as DG Lie algebras and as $C$-modules.
\end{prop}

\begin{proof}
Let $U' = \opn{Spec} C' \subset U$ be an open subset. 
Then $C \to C'$ is an \'etale ring homomorphism. According to 
\cite[Proposition 4.6]{Ye2} 
there are functorial DG Lie algebra homomorphisms
$\mcal{T}^{}_{\mrm{poly}}(C) \to \mcal{T}^{}_{\mrm{poly}}(C')$
and
$\mcal{D}^{}_{\mrm{poly}}(C) \to \mcal{D}^{}_{\mrm{poly}}(C')$
such that
$C' \otimes_C \mcal{T}^{}_{\mrm{poly}}(C) \cong
\mcal{T}^{}_{\mrm{poly}}(C')$
and
$C' \otimes_C \mcal{D}^{}_{\mrm{poly}}(C) \cong
\mcal{D}^{}_{\mrm{poly}}(C')$. 
Therefore we get quasi-coherent sheaves 
$\mcal{T}^{}_{\mrm{poly}, X}$ and $\mcal{D}^{}_{\mrm{poly}, X}$.

For any
$i \in \set{1, \ldots, p+1}$ let 
$\epsilon_i : \mcal{D}^{p}_{\mrm{poly}, X} \to
\mcal{D}^{p-1}_{\mrm{poly}, X}$
be the map
$\epsilon(\phi)(f_1, \ldots, f_p) \linebreak := 
\phi(f_1, \ldots, 1, \ldots, f_p)$,
with $1$ inserted at the $i$-th position. This is an 
$\mcal{O}_{X}$-linear homomorphism, and 
$\mcal{D}^{\mrm{nor}, p}_{\mrm{poly}, X} = \bigcap
\opn{Ker}(\epsilon_i)$.
Thus $\mcal{D}^{\mrm{nor}, p}_{\mrm{poly}, X}$ is quasi-coherent.

The left $\mcal{D}_X$-module structures on 
$\mcal{D}_{\mrm{poly}, X}$
and
$\mcal{D}^{\mrm{nor}}_{\mrm{poly}, X}$
are by composition of operators.  
\end{proof}

Following Kontsevich we call $\mcal{T}_{\mrm{poly}, X}$ the 
algebra of {\em poly vector fields} on $X$, and
$\mcal{D}_{\mrm{poly}, X}$ is called the 
algebra of {\em poly differential operators}.
The subalgebra $\mcal{D}^{\mrm{nor}}_{\mrm{poly}, X}$ is called 
the algebra of {\em normalized poly differential operators}.

Let us write  
$\mcal{T}_{\mrm{poly}}(X) = \Gamma(X, \mcal{T}_{\mrm{poly}, X})$,
the DG Lie algebra of global poly vector fields on $X$.
We consider each $\mcal{T}^p_{\mrm{poly}}(X)$ as a discrete inv 
module, and 
$\mcal{T}_{\mrm{poly}}(X) = \boplus_p 
\mcal{T}^p_{\mrm{poly}}(X)$
gets the $\boplus$ dir-inv structure, so it is a discrete DG Lie 
algebra in $\cat{Dir} \cat{Inv} \cat{Mod} \K$.
Likewise we define
$\mcal{D}_{\mrm{poly}}(X)$
and 
$\mcal{D}^{\mrm{nor}}_{\mrm{poly}}(X)$.

A series 
$\alpha = \sum_{k = 1}^{\infty} \alpha_k \hbar^k
\in \mcal{T}^1_{\mrm{poly}}(X)[[\hbar]]^+$
satisfying $[\alpha, \alpha] = 0$ is called a {\em formal 
Poisson structure} on $X$. 
Two formal Poisson structure $\alpha$ and
$\alpha'$ are called gauge equivalent if there is some
$\gamma =  \sum_{k = 1}^{\infty} \gamma_k \hbar^k \in 
\mcal{T}^0_{\mrm{poly}}(X)[[\hbar]]^+$
such that 
$\alpha' = \opn{exp}(\opn{ad}(\gamma))(\alpha)$. 
Thus the set 
$\mrm{MC} \big( \mcal{T}_{\mrm{poly}}(X) [[\hbar]]^+ \big)$
is the set of gauge equivalence classes of formal Poisson 
structures on $X$.

\begin{exa}
Let
$\alpha_1 \in \Gamma(X, \bwedge^2_{\mcal{O}_{X}} \mcal{T}_{X})$
be a Poisson structure on $X$ (Definition \ref{dfn7.2}). 
Then $\alpha := \alpha_1 \hbar$ is a formal Poisson structure.
\end{exa}

\begin{prop} \label{prop7.3}
An element  
\[ \beta = \sum_{j=1}^{\infty} \beta_j \hbar^j
\in \mcal{D}^{\mrm{nor}, 1}_{\mrm{poly}}(X) [[\hbar]]^+ \]
is a solution of the Maurer-Cartan equation in 
$\mcal{D}^{\mrm{nor}}_{\mrm{poly}}(X) [[\hbar]]^+$ 
iff the pairing
\[ (f, g) \mapsto f \star_{\beta} g 
:= f g + \sum_{j=1}^{\infty} \beta_j(f, g) \hbar^j , \]
for local sections $f, g \in \mcal{O}_X$, 
is a star product on $\mcal{O}_X[[\hbar]]$ 
\tup{(}see Definition \tup{\ref{dfn6.3})}.
\end{prop}

\begin{proof}
The assertion is actually local: it is enough to prove it for an 
affine open set $U = \opn{Spec} C \subset X$. Take
$\beta \in \mcal{D}^{\mrm{nor}, 1}_{\mrm{poly}}(C) [[\hbar]]^+$. 
We have to prove that $\star_{\beta}$ is an associative product on 
$C[[\hbar]]$ iff $\beta$ is a solution of the MC 
equation in $\mcal{D}^{\mrm{nor}}_{\mrm{poly}}(C) [[\hbar]]^+$. 
This assertion is made in \cite[Corollary 4.5]{Ke}. See also  
\cite[Section 4.6.2]{Ko1}. (For a non-differential 
star product  
this is the original discovery of Gerstenhaber, see \cite{Ge}.) 
\end{proof}

\begin{prop} \label{prop3.10}
Under the identification, in Proposition \tup{\ref{prop7.3}}, 
of solutions of the MC equation in
$\mcal{D}^{\mrm{nor}, 1}_{\mrm{poly}}(X)[\hbar]]^+$
with star products on $\mcal{O}_X[[\hbar]]$, 
the notion of gauge equivalence in Definition \tup{\ref{dfn7.1}} 
coincides with that in Proposition \tup{\ref{prop6.1}}.
\end{prop}

\begin{proof}
Let $\beta$ and $\beta'$ be two solutions of the MC equation in 
$\mcal{D}_{\mrm{poly}}^{\mrm{nor}}(X)[[\hbar]]^+$,
and let $\star$ and $\star'$ be the corresponding star products on
$\mcal{O}_X[[\hbar]]$. 
Given 
$\gamma \in \mcal{D}_{\mrm{poly}}^{\mrm{nor}, 0}(X)[[\hbar]]^+$,
let
$\opn{exp}(\gamma) := 1 + \gamma +
\smfrac{1}{2} \gamma^2 + \cdots$
be the corresponding gauge equivalence of 
$\mcal{O}_{X}[[\hbar]]$. As stated implicitly in 
\cite[Section 4.6.2]{Ko1} and 
\cite[Ch.\ 2, Lemma 4.2 and Section 5.1]{Ke},
one has
$\beta' = \opn{exp}(\opn{af})(\opn{exp}(\pm \gamma))(\beta)$
iff for all local sections $f, g \in \mcal{O}_{X}$ one has
\begin{equation} \label{eqn7.2}
f \star' g = \opn{exp}(\gamma)^{-1}
\big( \opn{exp}(\gamma)(f) \star
\opn{exp}(\gamma)(g) \big) . 
\end{equation}
(The reason for the sign ambiguity is that the references 
\cite{Ko1}, \cite{Ke} and \cite{GM} are inconsistent with each 
other regarding signs, and we did not carry out this calculation 
ourselves.)
\end{proof}

An immediate consequence is:

\begin{cor} \label{cor7.2}
The assignment $\beta \mapsto \star_{\beta}$ of Proposition 
\tup{\ref{prop7.3}} gives rise to a bijection from 
$\mrm{MC} 
\big( \mcal{D}^{\mrm{nor}}_{\mrm{poly}}(X) [[\hbar]]^+ \big)$
to set of gauge equivalence classes of globally trivialized 
deformation quantizations of $\mcal{O}_X$.
\end{cor}

Here is a first approximation of Theorem \ref{thm0.2}.

\begin{cor} \label{cor6.1}
Let $X$ be an $n$-dimensional affine scheme admitting an 
\'etale morphism $X \to \mbf{A}^n_{\K}$. 
Then there is a bijection $Q$ as in 
Theorem \tup{\ref{thm0.2}}.
\end{cor}

\begin{proof}
Write $X = \opn{Spec} C$
and $\mbf{A}^n_{\K} = \opn{Spec} \K[\bsym{t}]$.
Because $\K[\bsym{t}] \to C$ is an \'etale ring homomorphism, 
according to \cite[Proposition 4.6]{Ye2} we have 
$\mcal{T}_{\mrm{poly}}(C) = C \otimes_{\K[\bsym{t}]}
\mcal{T}_{\mrm{poly}}(\K[\bsym{t}])$
and
$\mcal{D}_{\mrm{poly}}(C) = C \otimes_{\K[\bsym{t}]}
\mcal{D}_{\mrm{poly}}(\K[\bsym{t}])$.
By condition (ii) in Theorem \ref{thm2.0} the universal operators 
$\mcal{U}_j$ are poly differential operators over $\K[\bsym{t}]$, 
and hence according to \cite[Proposition 2.6]{Ye2} they extend
to $C$-multilinear operators, giving an $\mrm{L}_{\infty}$ morphism
$\mcal{U} : \mcal{T}_{\mrm{poly}}(C) \to
\mcal{D}_{\mrm{poly}}(C)$;
and by \cite[Corollary 4.12]{Ye2} this is  an $\mrm{L}_{\infty}$
quasi-isomorphism. The inclusion 
$\mcal{D}^{\mrm{nor}}_{\mrm{poly}}(C) 
\to \mcal{D}_{\mrm{poly}}(C)$
is a DG Lie algebra quasi-isomorphism.
Now use Corollaries \ref{cor7.1} and \ref{cor7.2}.
\end{proof}

\begin{rem}
The method of $\mrm{L}_{\infty}$ morphisms is suitable only for 
characteristic $0$. For an approach in positive characteristic see 
\cite{BK2}. 
\end{rem}

\section{Formal Geometry -- Coordinate Bundles etc.}
\label{sec4} 

In this section we translate the notions of formal geometry (in 
the sense of Gelfand-Kazhdan \cite{GK}; 
cf.\ \cite[Section 7]{Ko1}) to the language of 
algebraic geometry (schemes and shea\-ves). As before $\K$ is a 
field of characteristic $0$, and $X$ is a smooth separated 
irreducible scheme over $\K$ of dimension $n$. 

For a 
closed point $x \in X$ the residue field $\bsym{k}(x)$ lifts 
uniquely into the complete local ring 
$\what{\mcal{O}}_{X, x}$, and any choice of system of coordinates 
$\bsym{t} = (t_1, \ldots, t_n)$ gives rise to an isomorphism of 
$\K$-algebras
\[ \what{\mcal{O}}_{X, x} \cong \bsym{k}(x)[[\bsym{t}]] = 
\bsym{k}(x)[[t_1, \ldots, t_n]] . \]
Of course the condition that an $n$-tuple of elements $\bsym{t}$
in the maximal ideal $\mfrak{m}_x$ is a system of coordinates is 
that their residue classes form a basis of the $\bsym{k}(x)$-module
$\mfrak{m}_x / \mfrak{m}_x^2$, the Zariski cotangent space.

Suppose $U$ is an open neighborhood of $x$ and 
$f \in \Gamma(U, \mcal{O}_X)$. The Taylor expansion of $f$ at $x$ 
w.r.t.\ $\bsym{t}$ is
\[ f = \sum_{\bsym{i} \in \mbb{N}^n} a_{\bsym{i}} \,
\bsym{t}^{\bsym{i}} \in \what{\mcal{O}}_{X, x} , \]
where $a_{\bsym{i}} \in \bsym{k}(x)$ and
$\bsym{t}^{\bsym{i}} := t_1^{i_1} \ldots t_n^{i_n}$. The 
coefficients are given by the usual formula
\[ a_{\bsym{i}} = \frac{1}{\bsym{i} !}
\Bigl( {\textstyle
\bigl( \frac{\partial}{\partial t_1} \bigr)^{i_1}
\cdots \bigl( \frac{\partial}{\partial t_n} \bigr)^{i_n} f 
} \Bigr) (x) , \]
where for any $g \in \what{\mcal{O}}_{X, x}$ we write 
$g(x) \in \bsym{k}(x)$ for its residue class.

The {\em jet bundle} of $X$ is an infinite dimensional scheme 
$\opn{Jet} X$, that comes with a projection
$\pi_{\mrm{jet}} : \opn{Jet} X \to X$. 
Given a closed point $x \in X$ the $\bsym{k}(x)$-rational points 
of the fiber $\pi_{\mrm{jet}}^{-1}(x)$ correspond to ``jets of 
functions at $x$'', namely to elements of the complete local ring 
$\what{\mcal{O}}_{X, x}$. Here is a way to visualize such a fiber: 
choose a coordinate system $\bsym{t}$. Then a jet is just the data
$\{ a_{\bsym{i}} \}_{\bsym{i} \in \mbb{N}^n}$
of its Taylor coefficients. So set theoretically the fiber 
$\pi_{\mrm{jet}}^{-1}(x)$ is just the set 
$\bsym{k}(x)^{\mbb{N}^n}$.

The naive description above does not make $\opn{Jet} X$ into a 
scheme. So let us try another approach. Consider the diagonal 
embedding $\Delta : X \to X^2 = X \times X$. 
Let $\mcal{I}_{X, \mrm{alg}}$ be the ideal sheaf
$\opn{Ker}(\Delta^* : \mcal{O}_{X^2} \to \mcal{O}_X)$.
The {\em sheaf of principal parts} of $X$ is 
\[ \mcal{P}_{X} = \mcal{P}_{X / \K} := \lim_{\leftarrow d}
\mcal{O}_{X^2} / \mcal{I}_{X, \mrm{alg}}^d  \]
(cf.\ \cite{EGA-IV}). 
It is a sheaf of commutative rings, equipped with 
two ring homomorphisms
$\mrm{p}^*_1, \mrm{p}^*_2 : \mcal{O}_X \to \mcal{P}_{X}$, 
namely $\mrm{p}^*_1 (f) := f \otimes 1$ and 
$\mrm{p}^*_2 (f) := 1 \otimes f$. We consider $\mcal{P}_{X}$
as a left $\mcal{O}_X$-module via $\mrm{p}^*_1$ and as a right 
$\mcal{O}_X$-module via $\mrm{p}^*_2$. 

The sheaf of rings $\mcal{P}_{X}$ can be thought of the 
the structure sheaf $\mcal{O}_{\mfrak{X}}$ of the formal scheme 
$\mfrak{X}$ which is the formal completion of $X^2$ along 
$\Delta(X)$. We denote by $\mcal{I}_{X}$ the ideal
$\opn{Ker}(\mcal{P}_{X} \to \mcal{O}_{X})$; it is just the 
completion of the ideal $\mcal{I}_{X, \mrm{alg}}$. 
By default we shall consider $\mcal{P}_X$ as an 
$\mcal{O}_X$-algebra via $\mrm{p}_1^*$.

\begin{prop}[{\cite[Lemma 2.6]{Ye1}}] \label{prop1.1}
Let $U \subset X$ be an open set admitting an \'etale morphism
$U \to \mbf{A}^n_{\K} = \opn{Spec} \K[s_1, \ldots, s_n]$. 
For $i = 1, \ldots, n$ define
\[ \til{s}_i := 1 \otimes s_i - s_i \otimes 1 \in
\Gamma(U, \mcal{I}) . \]
Then
\[ \mcal{P}_X|_U \cong \mcal{O}_U [[ \til{s}_1,
\ldots, \til{s}_n ]] \]
as sheaves of $\mcal{O}_U$-algebras, either via 
$\mrm{p}^*_1$ or via $\mrm{p}^*_2$.
\end{prop}

\begin{dfn}
Let $U \subset X$ be an open set. 
\begin{enumerate}
\item A {\em system of \'etale coordinates} on $U$ is a sequence 
$\bsym{s} = (s_1, \ldots, s_n)$ of elements in
$\Gamma(U, \mcal{O}_X)$ s.t.\ the morphism 
$U \to \mbf{A}^n_{\K}$ it determines is \'etale.
\item A {\em system of formal coordinates} on $U$ is a sequence 
$\bsym{t} = (t_1, \ldots, t_n)$ of elements in
$\Gamma(U, \mcal{I}_{X})$ s.t.\ the homomorphism of sheaves of 
rings $\mcal{O}_U [[\bsym{t}]] \to \mcal{P}_{X} |_U$
extending $\mrm{p}_1^*$ is an isomorphism.
\end{enumerate}
\end{dfn}

\begin{prop} \label{prop1.6}
Given a closed point $x \in X$ one has
a canonical isomorphism of $\mcal{O}_X$-algebras \tup{(}via 
$\mrm{p}^*_2$\tup{)}:
\[ \bsym{k}(x) \otimes_{\mcal{O}_X} \mcal{P}_{X} \cong
\what{\mcal{O}}_{X, x} . \]
If $\bsym{t} = (t_1, \ldots, t_n)$ 
is a system of formal coordinates on some 
neighborhood $U$ of $x$, and we let
$t_i(x) := 1 \otimes t_i \in \what{\mcal{O}}_{X, x}$
under the above isomorphism, then the sequence
$\bsym{t}(x) := (t_1(x), \ldots, t_n(x))$
is a system of coordinates in $\what{\mcal{O}}_{X, x}$.
\end{prop}

The easy proof is left out. 

\begin{exa} \label{exa1.1}
Assume $\bsym{s}$ is a system of \'etale coordinates on $U$, and 
let $t_i := \til{s}_i$. 
By Proposition \ref{prop1.1} the sequence
$\bsym{t} := (t_1, \ldots, t_n)$
is a system of formal coordinates on $U$. 
Given a closed point $x \in U$ we have
$t_i (x) = s_i - s_i(x)$, 
where $s_i(x) \in \bsym{k}(x) \subset \what{\mcal{O}}_{X, x}$.
The sequence
$\bsym{t}(x) = \big( t_1(x), \ldots, t_n(x) \big)$
is a system of coordinates in $\what{\mcal{O}}_{X, x}$.
\end{exa}

\begin{cor} \label{cor1.3}
Let $U \subset X$ be an open set admitting a formal
system of coordinates 
$\bsym{t} \in \Gamma(U, \mcal{I}_{X})^{\times n}$. 
For $f \in \Gamma(U, \mcal{O}_X)$ let us write
\[ \mrm{p}^*_2(f) = \sum_{\bsym{i}} 
\mrm{p}^*_1(a_{\bsym{i}})\, \bsym{t}^{\bsym{i}}
\in \Gamma(U, \mcal{P}_{X}) \]
with $a_{\bsym{i}} \in \Gamma(U, \mcal{O}_X)$. 
Then under the isomorphism of Proposition \tup{\ref{prop1.6}} 
we recover the Taylor expansion at any closed point $x \in U$:
\[ f = \sum_{\bsym{i}} 
a_{\bsym{i}}(x)\, \bsym{t}(x)^{\bsym{i}}
\in \what{\mcal{O}}_{X, x} . \]
\end{cor}

Again the easy proof is omitted.

The conclusion from Proposition \ref{prop1.6} is that the sheaf of  
sections of the bundle $\opn{Jet} X$ should be 
$\mcal{P}_{X}$. By the standard schematic formalism we deduce the 
defining formula
\[ \opn{Jet} X := \opn{Spec}_X \mrm{S}_{\mcal{O}_X}
\mcal{D}_X , \]
where 
$\mcal{D}_X = 
\mcal{H}\mathit{om}^{\mrm{cont}}_{\mcal{O}_X}
(\mcal{P}_X, \mcal{O}_X)$,
the sheaf of differential operators, is considered as a locally 
free left $\mcal{O}_X$-module; 
$\mrm{S}_{\mcal{O}_X} \mcal{D}_X$ is the symmetric algebra 
of the $\mcal{O}_X$-module $\mcal{D}_X$; and
$\opn{Spec}_X$ refers to the relative spectrum over $X$ of a 
quasi-coherent $\mcal{O}_X$-algebra.

If $U$ is a sufficiently small affine open set in $X$
admitting an \'etale coordinate system $\bsym{s}$, then 
$\opn{Jet} U$ can be made more explicit. We know that 
\[ \Gamma(U, \mcal{D}_X) = \boplus_{\bsym{i} \in \mbb{N}^n} 
\Gamma(U, \mcal{O}_X)
{\textstyle (\frac{\partial}{\partial s_1})^{i_1}
\cdots (\frac{\partial}{\partial s_n})^{i_n}} , \]
so letting $\xi_{\bsym{i}}$ be a commutative indeterminate we have
\[ \pi_{\mrm{jet}}^{-1}(U) \cong 
\opn{Jet} U \cong \opn{Spec} \Gamma(U, \mcal{O}_X)
[ \{ \xi_{\bsym{i}} \}_{\bsym{i} \in \mbb{N}^n} ] . \]

The next geometric object we need is the {\em bundle of formal 
coordinate systems} of $X$, which we denote by
$\opn{Coor} X$. (In \cite{Ko1} the notation is $X^{\mrm{coor}}$.
However we feel that stylistically $\opn{Coor} X$ is better, since 
it resembles the usual bundle notation $\mrm{T} X$ and 
$\mrm{T}^* X$.) The scheme $\opn{Coor} X$ comes with a projection 
$\pi_{\mrm{coor}} : \opn{Coor} X \to X$, and the fiber over a 
closed point $x \in X$ corresponds to the set of 
$\K$-algebra isomorphisms
$\bsym{k}(x) [[\bsym{t}]] \iso \what{\mcal{O}}_{X, x}$.
We are going to make this more precise below. 

Consider a morphism of schemes $f : Y \to X$. Unless $f$ is 
quasi-finite the sheaf $f^* \mcal{P}_X$ is not particularly 
interesting, but its completion is. To make this work nicely 
we are going to use inv structures (see Section \ref{sec2}).
Since $(X, \mcal{P}_X)$ is a noetherian formal scheme, any 
coherent $\mcal{P}_X$-module $\mcal{M}$ has 
the $\mcal{I}_X$-adic inv structure (see Example \ref{exa2.1}). 
Using the ring homomorphism 
$\mrm{p}_1^* : \mcal{O}_X \to \mcal{P}_{X}$
the module $\mcal{M}$ becomes an inv 
$\mcal{O}_{X}$-module, and so the complete inverse image 
$f^{\what{*}}\, \mcal{M}$ is defined. 
Taking $\mcal{M} := \mcal{P}_X$ we get an $\mcal{O}_{Y}$-algebra
$f^{\what{*}}\, \mcal{P}_X$. By Proposition \ref{prop1.1} 
we see that 
$f^{\what{*}}\, \mcal{I}_X$ is a sheaf of ideals in 
$f^{\what{*}}\, \mcal{P}_X$,
and $f^{\what{*}}\, \mcal{P}_X$ is 
$f^{\what{*}}\, \mcal{I}_X$ -adically complete in the usual sense, 
namely
$f^{\what{*}}\, \mcal{P}_X \cong \lim_{\leftarrow m}
f^{\what{*}}\, \mcal{P}_X / (f^{\what{*}}\, \mcal{I}_X)^m$.

Now we can state the geometric property that should characterize 
$\opn{Coor} X$. There should be a sequence 
$\bsym{t} = (t_1, \ldots, t_n)$ of elements in 
$\Gamma(\opn{Coor} X, \pi_{\mrm{coor}}^{\what{*}}\, \mcal{I}_{X})$
such that for any open set $U \subset X$, the assignment
$\sigma \mapsto \sigma^*(\bsym{t}) 
:= (\sigma^*(t_1), \ldots, \sigma^*(t_n))$
shall be a bijection from the set of sections 
$\sigma : U \to \opn{Coor} X$
to the set of formal coordinate systems on $U$.
Thus the sheaf of sections of 
$\opn{Coor} X$, let us call it $\mcal{Q}$, 
has to be a subsheaf of 
$\mcal{I}_{X}^{\times n} = \mcal{I}_{X} \times \cdots \times 
\mcal{I}_{X}$.
The condition for an $n$-tuple $\bsym{t}$ to be in $\mcal{Q}$ is 
that under the composed map
\begin{equation} \label{eqn1.3}
\mcal{I}_{X}^{\times n} \to 
(\mcal{I}_{X} / \mcal{I}_{X}^2)^{\times n} =
(\Omega^1_{X})^{\times n} \xar{\wedge} \Omega^n_{X} 
\end{equation}
one has
$t_1 \wedge \cdots \wedge t_n \neq 0$. 
Note that $\mcal{Q} = \lim_{\leftarrow m} \mcal{Q}^m$,
where $\mcal{Q}^1$ is the sheaf of frames of $\Omega^1_{X}$. 

We conclude that $\opn{Coor} X$ is a subscheme of 
$(\opn{Jet} X)^{\times n}$.
Specifically, $\opn{Coor} X$ is an open subscheme of 
$\opn{Spec}_X \mrm{S}_{\mcal{O}_X} 
\bigl( (\mcal{D}_X / \mcal{O}_X)^{\oplus n} \bigr)$,
where $\mcal{D}_X / \mcal{O}_X$ is viewed as a locally free left 
$\mcal{O}_X$-module. 
Over any affine open set $U \subset X$ admitting an \'etale 
coordinate system, say $\bsym{s}$, one has
\begin{equation} \label{eqn1.2}
\pi_{\mrm{coor}}^{-1}(U) \cong 
\opn{Coor} U \cong \opn{Spec} \Gamma(U, \mcal{O}_X)
[ \{ \xi_{\bsym{i}, j} \} , d^{-1}] . 
\end{equation}
In this formula $\bsym{i} = (i_1, \ldots, i_n)$
runs over $\mbb{N}^n - \{ (0, \ldots, 0) \}$
and $j$ runs over $\{ 1, \ldots, n \}$. The indeterminate
$\xi_{\bsym{i}, j}$ corresponds to the DO
$\frac{1}{\bsym{i}!} (\frac{\partial}{\partial s_1})^{i_1}
\cdots (\frac{\partial}{\partial s_n})^{i_n}$
in the $j$-th copy of $\mcal{D}_X$. 
The symbol $\bsym{e}_i$ denotes the row whose only nonzero entry 
is $1$ in the $i$-th place, so
$[\xi_{\bsym{e}_i, j}]$ is the matrix whose $(i, j)$ entry is 
the indeterminate $\xi_{\bsym{e}_i, j}$ 
corresponding to the DO $\frac{\partial}{\partial s_i}$
in the $j$-th copy of $\mcal{D}_X$. Finally 
$d := \opn{det}([\xi_{\bsym{e}_i, j}])$.

The next results justify the heuristic considerations above.

\begin{thm} \label{thm4.1}
Consider the functor
$F : (\msf{Sch} / X)^{\mrm{op}} \to \cat{Sets}$ 
defined as follows. For any $X$-scheme $Y$, with structural 
morphism $g : Y \to X$, we let $F Y$ be the
set of $\mcal{O}_Y$-algebra isomorphisms
$\phi : \mcal{O}_Y[[\bsym{t}]] \iso g^{\what{*}} \,  \mcal{P}_X$
such that
$\phi( \mcal{O}_Y[[\bsym{t}]] \cdot \bsym{t}) = 
g^{\what{*}} \,  \mcal{I}_X$. 
Then $\opn{Coor} X$ is a fine moduli space for $F$, namely 
$F \cong \opn{Hom}_{\msf{Sch} / X}(-, \opn{Coor} X)$.
\end{thm}

\begin{proof}
Suppose we are given $g : Y \to X$ and 
$\phi : \mcal{O}_Y[[\bsym{t}]] \iso g^{\what{*}} \,  \mcal{P}_X$. 
Define 
$a_i := \phi(t_i) \in \Gamma(Y, g^{\what{*}} \,  \mcal{I}_X)$. 
These elements satisfy
$a_1 \wedge \cdots \wedge a_n \neq 0$ like in equation
(\ref{eqn1.3}). Each $a_i$ gives rise to an 
$\mcal{O}_Y$-linear sheaf homomorphism
$\mcal{O}_Y \to g^{\what{*}} \,  \mcal{I}_X$. Since 
$\mcal{I}_X / \mcal{I}_X^m$ is a coherent 
locally free $\mcal{O}_X$-module 
for every $m \geq 1$ we see that
\[ \mcal{H}\mathit{om}^{\mrm{cont}}_{\mcal{O}_Y}
(g^{\what{*}} \,  \mcal{I}_X, \mcal{O}_Y) \cong
\mcal{O}_Y \otimes_{g^{-1} \mcal{O}_X}
g^{-1} \mcal{H}\mathit{om}^{\mrm{cont}}_{\mcal{O}_X}
(\mcal{I}_X, \mcal{O}_X) =
g^* (\mcal{D}_X / \mcal{O}_X) . \]
So after dualization, i.e.\ applying the functor
$\mcal{H}om^{\mrm{cont}}_{\mcal{O}_Y}(-, \mcal{O}_Y)$,
each $a_i$ gives  a homomorphism of $\mcal{O}_Y$-modules
$g^* (\mcal{D}_X / \mcal{O}_X) \to \mcal{O}_Y$.
By adjunction we get $\mcal{O}_X$-linear homomorphisms
$\mcal{D}_X / \mcal{O}_X \to g_* \mcal{O}_Y$, and therefore an 
$\mcal{O}_X$-algebra homomorphism
$\mrm{S}_{\mcal{O}_X} 
\bigl( (\mcal{D}_X / \mcal{O}_X)^{\oplus n} \bigr)
\to g_* \mcal{O}_Y$. 
Passing to schemes we obtain a morphism of $X$-schemes
\[ \til{\phi} : Y \to 
\opn{Spec}_X \mrm{S}_{\mcal{O}_X} 
\bigl( (\mcal{D}_X / \mcal{O}_X)^{\oplus n} \bigr) . \]
Because $a_1 \wedge \cdots \wedge a_n \neq 0$ this is actually a 
morphism $\til{\phi} : Y \to \opn{Coor} X$. The process we have 
described is reversible, and hence 
$F Y \cong \opn{Hom}_{\msf{Sch} / X}(Y, \opn{Coor} X)$. 
\end{proof}

\begin{cor} \label{cor4.1}  
There is a canonical isomorphism of 
$\mcal{O}_{\opn{Coor} X}$-algebras
\[ \mcal{O}_{\opn{Coor} X}[[\bsym{t}]] \cong
\pi_{\mrm{coor}}^{\what{*}} \, \mcal{P}_{X} . \]
This isomorphism has 
the following universal property: for any open set $U \subset X$ 
the assignment $\sigma \mapsto \sigma^*(\bsym{t})$
is a bijection of sets
\[ \opn{Hom}_{\msf{Sch} / X}(U, \opn{Coor} X) \iso
\{ \textup{formal coordinate systems on $U$} \} . \]
\end{cor}

\begin{proof}
Applying the theorem to $Y := \opn{Coor} X$, 
$g := \pi_{\mrm{coor}}$ and the identity morphism 
$\til{\phi}_0 : \opn{Coor} X \to \opn{Coor} X$,
we obtain a canonical isomorphism 
$\phi_0 : \mcal{O}_{\opn{Coor} X}[[\bsym{t}]] \cong
\pi_{\mrm{coor}}^{\what{*}} \, \mcal{P}_{X}$
with the desired universal property. 
\end{proof}

On $\opn{Coor} X$ we have a universal Taylor expansion: 

\begin{cor}
Suppose $U \subset X$ is open and $f \in \Gamma(U, \mcal{O}_X)$.
Then there are functions
$a_{\bsym{i}} \in 
\Gamma(\pi_{\mrm{coor}}^{-1}(U), \mcal{O}_{\opn{Coor} X})$
s.t.\
\[ \pi_{\mrm{coor}}^{\what{*}}(\mrm{p}_2^*(f)) = 
\sum_{\bsym{i} \in \mbb{N}^n} a_{\bsym{i}}\, 
\bsym{t}^{\bsym{i}} \in \Gamma(\pi_{\mrm{coor}}^{-1}(U), 
\pi_{\mrm{coor}}^{\what{*}} \, \mcal{P}_X) , \]
where $\bsym{t}$ is the universal coordinate system in
$\Gamma(\opn{Coor} X, \pi_{\mrm{coor}}^{\what{*}} \, \mcal{P}_X)$.
Given a section $\sigma : U \to \opn{Coor} X$ we obtain a Taylor 
expansion
\[ \mrm{p}_2^*(f) = 
\sum_{\bsym{i} \in \mbb{N}^n} \sigma^*(a_{\bsym{i}})\, 
\sigma^*(\bsym{t})^{\bsym{i}} \in
\Gamma(U, \mcal{P}_{X}) \]
as in Corollary \tup{\ref{cor1.3}}.
\end{cor}

The proof is left to the reader. 

Suppose $\bsym{s} = (s_1, \ldots, s_n)$ is an \'etale coordinate 
system on an open set $U \subset X$. As before let 
$\til{s}_i := 1 \otimes s_i - s_i \otimes 1 \in
\Gamma(U, \mcal{P}_X)$, and define
$\bsym{s} := (\til{s}_1, \ldots, \til{s}_n)$, 
which is a formal coordinate system on $U$. 
Then on 
$\opn{Coor} U = \pi_{\mrm{coor}}^{-1}(U)$
we have isomorphisms of $\mcal{O}_{\opn{Coor} U}$-algebras  
\[ \mcal{O}_{\opn{Coor} U}[[\bsym{t}]] \cong
(\pi_{\mrm{coor}}^{\what{*}} \mcal{P}_X)|_{\opn{Coor} U} \cong
\mcal{O}_{\opn{Coor} U}[[\bsym{s}]] . \]
Using the coordinate functions 
$\xi_{\bsym{i}, j} \in 
\Gamma(\opn{Coor} U, \mcal{O}_{\opn{Coor} U})$
from formula (\ref{eqn1.2}) we then have
\begin{equation}
t_j = \sum_{\bsym{i}} \xi_{\bsym{i}, j} \bsym{s}^{\bsym{i}} , 
\end{equation}
where the sum is on 
$\bsym{i} \in \mbb{N}^n - \{ (0, \ldots, 0) \}$.

For $i \geq 1$ let $\opn{Coor}^i X$ be the bundle over $X$ 
parameterizing coordinate systems up to order $i$ (i.e.\ modulo 
order $\geq i+1$). There are projections 
$\opn{Coor} X \to \opn{Coor}^i X \to \opn{Coor}^{i-1} X \to X$.
The next theorem describes the geometry of these bundles.

Let $G(\K)$ be the group of $\K$-algebra automorphisms of 
$\K[[\bsym{t}]]$. Then $G(\K)$ is the group of $\K$-rational 
points of a pro-algebraic group 
$G = \opn{GL}_{n, \K} \ltimes N$, where $N$ is a pro-unipotent 
group. The action of $\opn{GL}_{n}(\K)$ on $\K[[\bsym{t}]]$ is by 
linear change of coordinates; and $N(\K)$ is the subgroup of 
$G(\K)$ consisting of automorphisms that act trivially 
modulo $(\bsym{t})^2$. 

According to Corollary \ref{cor4.1} there is a canonical embedding 
of $\K$-algebras
\begin{equation} \label{eqn4.1}
\K[[\bsym{t}]] \inj 
\Gamma(\opn{Coor} X, \pi_{\mrm{coor}}^{\what{*}} \mcal{P}_X) .
\end{equation}
 
\begin{thm} \label{thm1.1}
\begin{enumerate}
\item $\opn{Coor} X \cong \lim_{\leftarrow i} \opn{Coor}^i X $
as schemes over $X$. 
\item $\opn{Coor} X$ is a $G$-torsor over $X$.
The action of $G$ on $\opn{Coor} X$ is characterized by the fact 
that the embedding \tup{(\ref{eqn4.1})} is $G(\K)$-equivariant. 
\item $\opn{Coor}^1 X$ is a $\opn{GL}_{n, \K}$-torsor over
$X$, and $\opn{Coor} X$ is a $\opn{GL}_{n, \K}$-equivariant 
$N$-torsor over $\opn{Coor}^{1} X$. 
\item The geometric quotient \tup{(}cf.\ \cite{GIT}\tup{)}
\[ \opn{LCC} X := \opn{Coor} X / \opn{GL}_{n, \K} \]
exists, with projection
$\pi_{\mrm{gl}} : \opn{Coor} X \to \opn{LCC} X$, 
and $\opn{Coor} X$ is a $\opn{GL}_{n, \K}$-torsor over 
$\opn{LCC} X$.
\item Let $U \subset X$ be an affine open set admitting an \'etale 
coordinate system. Then all the torsors in parts \tup{(2-4)}
are trivial over $U$ \tup{(}i.e.\ they admit sections\tup{)}.  
\end{enumerate}
\end{thm}

``LCC'' stands for ``linear coordinate classes''. 
In \cite{Ko1} the notation for $\opn{LCC} X$ is
$X^{\mrm{aff}}$. Note that the bundle $\opn{LCC} X$ has no group 
action; but locally, for $U$ as in part (5), 
there's a non-canonical isomorphism of schemes
$\opn{LCC} U \cong N \times U$. $\opn{Coor}^1 X$ is the frame 
bundle of $\Omega^1_X$. 
The various bundles and projections are depicted in
Figure \ref{fig0}. 

\begin{figure}
\[ \UseTips \xymatrix{
& \opn{Coor} X
\ar[dl]_{} 
\ar[dr]^{\pi_{\mrm{gl}}} 
\ar[dd]^{\pi_{\mrm{coor}}} 
\\
\opn{Coor}^i X
\ar[dr]
& &
\opn{LCC} X
\ar[dl]^{\pi_{\mrm{lcc}}}
\\
& X 
} \]
\caption{} \label{fig0}
\end{figure}

\begin{proof}
(1) This is an immediate consequence of the moduli property of 
$\opn{Coor} X$ (see Theorem \ref{thm4.1}), and an analogous 
property of $\opn{Coor}^i X$.

\medskip \noindent
(2) Given $g \in G(\K)$ let us denote by $g(\bsym{t})$ the 
sequence $(g(t_1), \ldots, g(t_n))$ in $\K[[\bsym{t}]]$.
By Theorem \ref{thm4.1} there exists a unique $X$-morphism
$\til{g} : \opn{Coor} X \to \opn{Coor} X$
such that the algebra homomorphism
$\til{g}^* : 
\Gamma(\opn{Coor},  \pi_{\mrm{coor}}^{\what{*}} \mcal{P}_X) \to
\Gamma(\opn{Coor},  \pi_{\mrm{coor}}^{\what{*}} \mcal{P}_X)$
sends $\bsym{t}$ to $g(\bsym{t})$. 
We have to prove that $\til{g}$ is an automorphism, and that
$g \mapsto \til{g}$ is a group homomorphism from $G(\K)$ to
$\opn{Aut}_{\msf{Sch} / X}(\opn{Coor} X)$.

Now via the embedding (\ref{eqn4.1}), the homomorphism
$\til{g}^*$ restricts to the automorphism $g$ on 
$\K[[\bsym{t}]]$. If $g$ is the identity automorphism of 
$\K[[\bsym{t}]]$, then by uniqueness $\til{g}$ has to be the 
identity automorphism of $\opn{Coor} X$. Next take two elements
$g_1, g_2 \in G(\K)$. Then
\[ \begin{aligned}
& \wtil{g_2 \circ g_1}^*(\bsym{t}) = (g_2 \circ g_1)(\bsym{t})
= g_2 (g_1(\bsym{t})) = \til{g}_2^* (g_1(\bsym{t})) =
g_1( \til{g}_2^*(\bsym{t})) \\
& \qquad = g_1(g_2(\bsym{t})) =
(\til{g}_1^* \circ \til{g}_2^*)(\bsym{t})
= (\til{g}_2 \circ \til{g}_1)^*(\bsym{t}) .
\end{aligned} \] 
Thus indeed we have a group action. 

Due to the moduli property this action becomes geometric, i.e.\ it 
is a morphism of schemes 
$G \times \opn{Coor} X \to \opn{Coor} X$.
The explicit local description (\ref{eqn1.2}) 
shows that $\opn{Coor} X$ is in fact a $G$-torsor over $X$.

\medskip \noindent
(3, 4) These are consequence of (2). 

\medskip \noindent
(5) Clear from formula (\ref{eqn1.2}).
\end{proof}

\section{Formal Differential Calculus}
\label{sec5}

As before $\K$ is a field of characteristic $0$, and $X$ is a 
smooth separated irreducible $n$-dimensional $\K$-scheme.

Recall the algebra homomorphism 
$\mrm{p}_1^* : \mcal{O}_{X} \to \mcal{P}_{X}$.
We define 
$\mcal{T}(\mcal{P}_X / \mcal{O}_{X}; \mrm{p}_1^*)$
to be the sheaf of derivations of $\mcal{P}_{X}$ relative to 
$\mcal{O}_{X}$. Thus for any affine open set
$U = \opn{Spec} C \subset X$, writing
$\what{A} := \Gamma(U, \mcal{P}_X)$, we have
\[ \Gamma \bigl( U, \mcal{T}(\mcal{P}_X / \mcal{O}_{X}; \mrm{p}_1^*)
\bigr) =
\mcal{T}_{\what{A} / C} = \opn{Der}_C(\what{A}) . \]
Similarly we define
$\mcal{T}^{i}_{\mrm{poly}}
(\mcal{P}_X / \mcal{O}_{X}; \mrm{p}_1^*)$
and
$\mcal{D}^{i}_{\mrm{poly}}
(\mcal{P}_X / \mcal{O}_{X}; \mrm{p}_1^*)$.

\begin{lem} \label{lem5.1} 
Let $\mcal{G}$ stand either for $\mcal{T}_{\mrm{poly}}$ or
$\mcal{D}_{\mrm{poly}}$, so that 
$\mcal{G}_X = \mcal{T}_{\mrm{poly}, X}$ etc. 
\begin{enumerate}
\item The graded left $\mcal{O}_{X}$-module 
$\mcal{P}_{X} \otimes_{\mcal{O}_X} \mcal{G}_{X}$
is a DG Lie algebra in \linebreak
$\cat{Dir} \cat{Inv} \cat{Mod} \mcal{O}_{X}$. The homomorphism
$\mcal{G}_{X} \to \mcal{P}_{X} \otimes_{\mcal{O}_X}
\mcal{G}_{X}$
given by $\gamma \mapsto 1 \otimes \gamma$
is a DG Lie algebra homomorphism.
\item There is a canonical isomorphism 
\[ \mcal{P}_{X} \otimes_{\mcal{O}_X} \mcal{G}_{X}
\cong
\mcal{G}(\mcal{P}_{X}  / \mcal{O}_X; \mrm{p}_1^*) \]
of sheaves of DG Lie algebras in 
$\cat{Dir} \cat{Inv} \cat{Mod} \mcal{O}_{X}$.
\item Suppose $f : Y \to X$ is a morphism of schemes. Then
\[ f^{\what{*}} \,
\mcal{G}(\mcal{P}_{X}  / \mcal{O}_X; \mrm{p}_1^*)
\cong
\mcal{G}(f^{\what{*}}  \mcal{P}_{X}  / \mcal{O}_Y)  \]
as DG Lie algebras in $\cat{Dir} \cat{Inv} \cat{Mod} \mcal{O}_{Y}$.
\end{enumerate}
\end{lem}

\begin{proof}
(1) Let $U \subset X$ be an affine open set, and define
$C := \Gamma(U, \mcal{O}_X)$ and
$A := C \otimes C$. Let $\mfrak{a} := \opn{Ker}(A \to C)$,
and let $\what{A}$ be the $\mfrak{a}$-adic completion of $A$.
The left $C$-module
$C \otimes \mcal{G}(C)$ is a DG Lie algebra over $C$.
When we consider $C \otimes \mcal{G}(C)$ 
as an $A$-module, the bracket
\[ [-,-] : (C \otimes \mcal{G}(C)) \times 
(C \otimes \mcal{G}(C)) \to
C \otimes \mcal{G}(C) \]
and the differential
\[ \d : C \otimes \mcal{G}(C) \to
C \otimes \mcal{G}(C) \]
are poly differential operators, and hence they are
continuous for the $\mfrak{a}$-adic dir-inv structure 
(see \cite[Example 1.8]{Ye2}).
So according to \cite[Proposition 2.3]{Ye2},
$\what{C \otimes \mcal{G}(C)}$
is a DG Lie algebra in $\cat{Dir} \cat{Inv} \cat{Mod} C$.
But
\[ \Gamma(U, \mcal{P}_{X} \otimes_{\mcal{O}_X} 
\mcal{G}_X) \cong \what{A} \otimes_C \mcal{G}(C)
\cong \what{C \otimes \mcal{G}(C)} . \]

\medskip \noindent
(2, 3) By definition $\mcal{G}(C) = \mcal{G}(C / \K)$. 
By base change there is an isomorphism
$B \otimes \mcal{G}(C / \K) \cong 
\mcal{G} \bigl( (C \otimes B) / B \bigr)$
for any $\K$-algebra $B$.
\end{proof}

\begin{dfn} \label{dfn5.2}
Consider the de Rham differential
$\d : \mcal{O}_{X^2} \to \Omega^1_{X^2 / X} = 
\mrm{p}^*_1\, \Omega^1_X$
relative to the projection  
$\mrm{p}_2 : X^2 \to X$. Passing to the completion along the 
diagonal we obtain the {\em Grothendieck connection}
\[ \nabla_{\mcal{P}} : \mcal{P}_X \to
\mcal{P}_X \otimes_{\mcal{O}_{X}} \Omega^1_X . \]
\end{dfn}

Let $\mcal{M}$ be an $\mcal{O}_{X}$-module. Then the connection
$\nabla_{\mcal{P}}$ extends uniquely
to a degree $1$ endomorphism of the graded sheaf
\[ \Omega^{}_{X} \otimes_{\mcal{O}_{X}} 
\mcal{P}_X \otimes_{\mcal{O}_{X}} \mcal{M} = \boplus_{p \geq 0}
\Omega^{p}_{X} \otimes_{\mcal{O}_{X}} 
\mcal{P}_X \otimes_{\mcal{O}_{X}} \mcal{M} . \]
The formula is
\[ \nabla_{\mcal{P}}(\alpha \otimes a \otimes m) :=
\d(\alpha) \otimes a \otimes m + (-1)^p \alpha \wedge 
\nabla_{\mcal{P}}(a) \otimes m \]
for local sections
$\alpha \in \Omega^{p}_{X}$, $a \in \mcal{P}_X$ and 
$m \in \mcal{M}$. The connection is integrable, i.e.\
$\nabla_{\mcal{P}} \circ \nabla_{\mcal{P}} = 0$, and it makes 
$\Omega^{}_{X} \otimes_{\mcal{O}_{X}} 
\mcal{P}_X \otimes_{\mcal{O}_{X}} \mcal{M}$
into a DG $\Omega_X$-module. 

\begin{thm}[{\cite[Theorem 4.4]{Ye3}}]
Let $X$ be a smooth $\K$-scheme and let $\mcal{M}$ be an 
$\mcal{O}_{X}$-module. Then the map
\[ \mcal{M} \to \Omega^{}_{X} \otimes_{\mcal{O}_{X}} 
\mcal{P}_X \otimes_{\mcal{O}_{X}} \mcal{M}, \
m \mapsto 1 \otimes 1 \otimes m \]
is a quasi-isomorphism.
\end{thm}

Given any $\K$-scheme $Y$ let $\K_Y$ be the constant sheaf $\K$ on 
$Y$. We consider $\Omega^p_Y = \Omega^p_{Y / \K}$ as a discrete inv 
$\mcal{O}_Y$-module, and 
$\Omega_Y = \boplus_{p \geq 0} \Omega^p_Y$
gets the $\boplus$ dir-inv structure. Thus $\Omega_Y$
is a discrete DG algebra in $\cat{Dir} \cat{Inv} \cat{Mod} \K_Y$.
Note that if $Y$ is infinite dimensional then 
$\Omega_Y$ will be unbounded.

Suppose $\mcal{M}$ is a quasi-coherent $\mcal{O}_{X}$-module. 
Then for any $p$ the sheaf
$\Omega^{p}_{X} \otimes_{\mcal{O}_{X}} 
\mcal{P}_X \otimes_{\mcal{O}_{X}} \mcal{M}$
is a dir-coherent $\mcal{P}_X$-module (see Example \ref{exa2.1}),
so it has the $\mcal{I}_X$-adic dir-inv module structure.
The connection 
\[ \nabla_{\mcal{P}} : 
\Omega^{p}_{X} \otimes_{\mcal{O}_{X}} 
\mcal{P}_X \otimes_{\mcal{O}_{X}} \mcal{M} \to
\Omega^{p+1}_{X} \otimes_{\mcal{O}_{X}} 
\mcal{P}_X \otimes_{\mcal{O}_{X}} \mcal{M} \]
is a differential operator of $\mcal{P}_X$-modules (of order 
$\leq 1$), and therefore it is continuous for the dir-inv 
structures (see \cite[Proposition 2.3]{Ye2}). So in fact 
$\Omega^{}_{X} \otimes_{\mcal{O}_{X}} \mcal{P}_X 
\otimes_{\mcal{O}_{X}} \mcal{M}$
is a DG $\Omega^{}_X$-module in 
$\cat{Dir} \cat{Inv} \cat{Mod} \K_X$.

Suppose $f : Y \to X$ is some morphism of schemes. 
The complete pullback 
$f^{\what{*}} (\mcal{P}_X \otimes_{\mcal{O}_{X}}
\mcal{M})$
is a dir-inv $\mcal{O}_{Y}$-module. Moreover 
$\Omega^{}_{Y} \hatotimes{\mcal{O}_{Y}}
f^{\what{*}} (\mcal{P}_X \otimes_{\mcal{O}_{X}} \mcal{M})$
is a DG $\Omega^{}_Y$-module in 
$\cat{Dir} \cat{Inv} \cat{Mod} \K_Y$.
Its differential is also denoted by $\nabla_{\mcal{P}}$. 
In particular, when $\mcal{M} = \mcal{O}_{X}$, we obtain a 
super-commutative associative unital DG algebra
$\Omega_{Y} \hatotimes{\mcal{O}_{Y}}
f^{\what{*}} \mcal{P}_X$
in $\cat{Dir} \cat{Inv} \cat{Mod} \K_Y$.
Its degree $0$ component is $f^{\what{*}} \mcal{P}_X$,
which is a complete commutative algebra in 
$\cat{Inv} \cat{Mod} \mcal{O}_{Y}$. 
For details and proofs see \cite[Section 1]{Ye2}. 

\begin{prop}  \label{prop5.1}
Let $\mcal{G}$ denote either $\mcal{T}^{}_{\mrm{poly}}$ or
$\mcal{D}^{}_{\mrm{poly}}$, so that 
$\mcal{G}_X = \mcal{T}^{}_{\mrm{poly}, X}$ etc.
Also let $\d_{\mcal{G}}$ and $[-,-]_{\mcal{G}}$ denote the 
differential and the bracket of $\mcal{G}$.
\begin{enumerate}
\item The graded sheaf 
$\Omega^{}_X \otimes_{\mcal{O}_{X}} \mcal{P}_X
\otimes_{\mcal{O}_{X}} \mcal{G}_X$
is a DG $\Omega^{}_X$-module Lie algebra in 
$\cat{Dir} \cat{Inv} \cat{Mod} \K_X$. The differential is 
$\nabla_{\mcal{P}} + \bsym{1} \otimes \bsym{1} \otimes 
\d_{\mcal{G}}$, 
and the bracket is the continuous $\Omega^{}_X$-bilinear 
extension of $[-,-]_{\mcal{G}}$.
\item The canonical map
$\mcal{G}_X \to \Omega^{}_X \otimes_{\mcal{O}_{X}} \mcal{P}_X
\otimes_{\mcal{O}_{X}} \mcal{G}_X$
is a DG Lie algebra quasi-isomorphism. 
\item Suppose $f : Y \to X$ is a morphism of schemes. Then 
\[ \Omega_{Y} \hatotimes{\mcal{O}_{Y}}
f^{\what{*}} (\mcal{P}_X \otimes_{\mcal{O}_{X}} \mcal{G}_X) \]
is a DG $\Omega_{Y}$-module Lie algebra in 
$\cat{Dir} \cat{Inv} \cat{Mod} \K_Y$. 
The canonical map
\[ f^{-1} (\Omega^{}_X \otimes_{\mcal{O}_{X}} \mcal{P}_X
\otimes_{\mcal{O}_{X}} \mcal{G}_X) \to 
\Omega_{Y} \hatotimes{\mcal{O}_{Y}}
f^{\what{*}} (\mcal{P}_X \otimes_{\mcal{O}_{X}} \mcal{G}_X) \]
is a homomorphism of DG Lie algebras.
\end{enumerate}
\end{prop}

The explicit formulas for part (1) are
\[ \d(\alpha_1 \otimes 1 \otimes \gamma_1) :=
\d(\alpha_1) \otimes 1 \otimes \gamma_1 
+ (-1)^{i_1} \alpha_1 \otimes 1 \otimes \d_{\mcal{G}}(\gamma_1)
\]
and
\[ [\alpha_1 \otimes 1 \otimes \gamma_1, 
\alpha_2 \otimes 1 \otimes \gamma_2] := (-1)^{j_1 i_2}
(\alpha_1 \wedge \alpha_2) \otimes 1 \otimes
[\gamma_1, \gamma_2] \]
for $\alpha_k \in \Omega^{i_k}_X$
and $\gamma_k \in \mcal{G}_X^{j_k}$.

\begin{proof}
(1) Using the notation of the proof of Lemma \ref{lem5.1},
$\Omega_C \otimes \mcal{G}(C)$ is a DG $\Omega_C$-module
Lie algebra in $\cat{Dir} \cat{Inv} \cat{Mod} \K$, with the 
$\mfrak{a}$-adic dir-inv structure. Hence so is its completion
\[ \what{\Omega_C \otimes \mcal{G}(C)} \cong
\Omega_C \otimes_C \what{A} \otimes_C  \mcal{G}(C)
\cong \Gamma(U, \Omega_X \otimes_{\mcal{O}_{X}}
\mcal{P}_X \otimes_{\mcal{O}_{X}} \mcal{G}_{X}) . \]

\medskip \noindent
(2) In the the proof of part (1) the inclusion
$\mcal{G}(C) \subset \Omega_C \otimes \mcal{G}(C)$
is a DG algebra homomorphism. According to \cite[Theorem 3.4]{Ye3}
it is a quasi-isomorphism.

\medskip \noindent
(3) This is by \cite[Proposition 1.22(2)]{Ye2}.
\end{proof}

\begin{dfn}
Let $\mcal{G}$ denote either $\mcal{T}^{}_{\mrm{poly}}$ or
$\mcal{D}^{}_{\mrm{poly}}$, and let $\d$ be the de Rham 
differential on $\Omega_{\opn{Coor} X}$. 
Put on $\mcal{G}(\K[[\bsym{t}]])$ the $\bsym{t}$-adic dir-inv 
structure. Define
\[ \mrm{d}_{\mrm{for}} := \d \otimes \bsym{1} :
\Omega^p_{\opn{Coor} X} \hatotimes{} \mcal{G}(\K[[\bsym{t}]]) \to
\Omega^{p+1}_{\opn{Coor} X} \hatotimes{} \mcal{G}(\K[[\bsym{t}]]) 
. \]
\end{dfn}

According to \cite[Proposition 1.19]{Ye2}, 
\[ \Omega_{\opn{Coor} X} \hatotimes{} \mcal{G}(\K[[\bsym{t}]]) =
\boplus_{p, q}
\Omega^p_{\opn{Coor} X} \hatotimes{} \mcal{G}^q(\K[[\bsym{t}]]) \]
is a DG Lie algebra in 
$\cat{Dir} \cat{Inv} \cat{Mod} \K_{\opn{Coor} X}$,
with differential
$\mrm{d}_{\mrm{for}} + \bsym{1} \otimes \d_{\mcal{G}}$.
The explicit formula is
\[ (\mrm{d}_{\mrm{for}} + \bsym{1} \otimes \d_{\mcal{G}})
(\alpha \otimes \gamma) = 
\d(\alpha) \otimes \gamma + 
(-1)^p \alpha \otimes \d_{\mcal{G}}(\gamma) \]
for $\alpha \in \Omega_{\opn{Coor} X}^p$ and 
$\gamma \in \mcal{G}(\K[[\bsym{t}]])$.

\begin{thm}[Universal Taylor Expansion] \label{thm5.3} 
Let $\mcal{G}$ denote either $\mcal{T}^{}_{\mrm{poly}}$ or
$\mcal{D}^{}_{\mrm{poly}}$. There is a canonical isomorphism 
\[ \Omega_{\opn{Coor} X} \, 
\what{\otimes}_{\mcal{O}_{\opn{Coor} X}} \,
\pi_{\mrm{coor}}^{\what{*}} (\mcal{P}_X \otimes_{\mcal{O}_{X}}
\mcal{G}_{X}) \cong 
\Omega_{\opn{Coor} X} \, 
\what{\otimes} \, 
\mcal{G}^{}(\K [[ \bsym{t} ]])  \]
of graded Lie algebras in
$\cat{Dir} \cat{Inv} \cat{Mod} \mcal{O}_{\opn{Coor} X}$, 
extending the isomorphism of Corollary \tup{\ref{cor4.1}}.
\end{thm}

Warning: the isomorphism in the theorem does 
not respect the differentials; cf.\ Proposition \ref{prop2.3} 
below.

\begin{proof}
By Corollary \ref{cor4.1} we know that
$\pi_{\mrm{coor}}^{\what{*}}\, \mcal{P}_{X} \cong 
\mcal{O}_{\opn{Coor} X}[[\bsym{t}]]$
canonically as inv $\mcal{O}_{\opn{Coor} X}$-algebras.
Using Lemma \ref{lem5.1} we then obtain isomorphisms of graded Lie 
algebras over $\mcal{O}_{\opn{Coor} X}$
\[ \begin{aligned}
\pi_{\mrm{coor}}^{\what{*}} \,
(\mcal{P}_{X} \otimes_{\mcal{O}_X} 
\mcal{G}_{X}) & \cong 
\pi_{\mrm{coor}}^{\what{*}} \, \mcal{G}
(\mcal{P}_X / \mcal{O}_X; \mrm{p}_1^*) \\
& \cong 
\mcal{G}(\pi_{\mrm{coor}}^{\what{*}} 
\mcal{P}_{X} / \mcal{O}_{\opn{Coor} X}) \\
& \cong
\mcal{G}(\mcal{O}_{\opn{Coor} X}[[\bsym{t}]] /
\mcal{O}_{\opn{Coor} X}) \\
& \cong
\mcal{O}_{\opn{Coor} X} \, 
\what{\otimes} \, \mcal{G}(\K[[\bsym{t}]]) .
\end{aligned} \]
Finally we may apply the functor
$\Omega_{\opn{Coor} X} \, \what{\otimes}_{\mcal{O}_{\opn{Coor} X}} 
\, -$.
\end{proof}

\begin{dfn} \label{dfn2.3}
The {\em Maurer-Cartan form of $X$} is
\[ \begin{aligned}
& \omega_{\mrm{MC}} := \sum_{i = 1}^n
\nabla_{\mcal{P}}(t_i) \cdot \smfrac{\partial}{\partial t_i}
\in \Gamma \big( \opn{Coor} X, \Omega^1_{\opn{Coor} X} \, 
\what{\otimes} \, 
\mcal{T}^{0}_{\mrm{poly}}(\K [[ \bsym{t} ]]) \big) ,
\end{aligned} \]
where $\nabla_{\mcal{P}}(t_i)$ is defined 
using the canonical isomorphism in Theorem \ref{thm5.3}.
\end{dfn}

The Lie algebra 
$\mcal{T}^{0}_{\mrm{poly}}(\K [[\bsym{t}]]) 
= \mcal{T}(\K [[\bsym{t}]])$
is also a Lie subalgebra of 
$\mcal{D}^{0}_{\mrm{poly}}(\K [[\bsym{t}]])$ \linebreak 
$= \mcal{D}(\K [[\bsym{t}]])$. 
Keeping the notation of Theorem \ref{thm5.3}, for any local section
$\alpha \in \Omega_{\opn{Coor} X} \, 
\what{\otimes} \, \mcal{G}(\K [[ \bsym{t} ]])$
let 
\[ \opn{ad}(\omega_{\mrm{MC}})(\alpha) := 
[\omega_{\mrm{MC}}, \alpha] . \]
The operation $\opn{ad}(\omega_{\mrm{MC}})$ is $\K$-linear 
endomorphism of degree $1$ of the graded sheaf
$\Omega_{\opn{Coor} X} \, \what{\otimes} \, 
\mcal{G}(\K [[ \bsym{t} ]])$.

\begin{prop} \label{prop2.3}
Let $\mcal{G}$ denote either $\mcal{T}^{}_{\mrm{poly}}$ or
$\mcal{D}^{}_{\mrm{poly}}$.
Under the isomorphism of Theorem \tup{\ref{thm5.3}}, 
there is equality
\[ \nabla_{\mcal{P}} = \mrm{d}_{\mrm{for}} + 
\opn{ad}(\omega_{\mrm{MC}}) \]
as endomorphisms of 
$\Omega_{\opn{Coor} X} \, \what{\otimes} \, 
\mcal{G}(\K [[ \bsym{t} ]])$.
\end{prop}

\begin{proof}
First we shall consider $\mcal{G} = \mcal{T}_{\mrm{poly}}$.
Let's write
$\omega := \omega_{\mrm{MC}}$,
$\mrm{d} := \mrm{d}_{\mrm{for}}$,
$\nabla := \nabla_{\mcal{P}}$ and
$\partial_i := \smfrac{\partial}{\partial t_i}$.
For any multi-index
$\bsym{j} = (j_1 < \cdots < j_q)$
let's write 
\[ \bsym{\partial}_{\bsym{j}} := \partial_{j_1} \wedge \cdots
\wedge \partial_{j_q} \in 
\mcal{T}^{q-1}_{\mrm{poly}}(\K [[ \bsym{t} ]]) . \] 
Take a local section 
$\alpha \in \Omega^p_{\opn{Coor} X}$ and
a multi-index $\bsym{i} \in \mbb{N}^n$, and consider
$\alpha \otimes \bsym{t}^{\bsym{i}} \bsym{\partial}_{\bsym{j}} \in
\Omega^p_{\opn{Coor} X} \, \what{\otimes} \, 
\mcal{T}^{q-1}_{\mrm{poly}}(\K [[ \bsym{t} ]])$.
Then
\[ \begin{aligned}
\nabla
(\alpha \otimes \bsym{t}^{\bsym{i}} \bsym{\partial}_{\bsym{j}}) 
& = \mrm{d}(\alpha) \cdot \bsym{t}^{\bsym{i}} \cdot
\bsym{\partial}_{\bsym{j}}
\pm \alpha \cdot \nabla(\bsym{t}^{\bsym{i}}) \cdot
\bsym{\partial}_{\bsym{j}} \pm
\alpha \cdot \bsym{t}^{\bsym{i}} \cdot
\nabla(\bsym{\partial}_{\bsym{j}}) \\
\mrm{d} 
(\alpha \otimes \bsym{t}^{\bsym{i}} \bsym{\partial}_{\bsym{j}}) 
& = \mrm{d}(\alpha) \cdot \bsym{t}^{\bsym{i}} \cdot
\bsym{\partial}_{\bsym{j}} \\
\opn{ad}(\omega)
(\alpha \otimes \bsym{t}^{\bsym{i}} \bsym{\partial}_{\bsym{j}}) 
& = \pm \alpha \cdot \opn{ad}(\omega)(\bsym{t}^{\bsym{i}}) \cdot
\bsym{\partial}_{\bsym{j}} \pm
\alpha \cdot \bsym{t}^{\bsym{i}} \cdot
\opn{ad}(\omega)(\bsym{\partial}_{\bsym{j}}) .
\end{aligned} \]
Now
\[ \nabla(\bsym{t}^{\bsym{i}}) = 
\sum_k \nabla(t_k) \cdot \partial_k(\bsym{t}^{\bsym{i}}) = 
\opn{ad}(\omega)(\bsym{t}^{\bsym{i}}) . \]
It remains to show that
\[ \nabla(\bsym{\partial}_{\bsym{j}}) =
\opn{ad}(\omega)(\bsym{\partial}_{\bsym{j}}) . \]

Take any 
$\beta \in \Omega_{\opn{Coor} X} \, \what{\otimes} \, 
\mcal{T}_{\mrm{poly}}^{q-1}(\K [[ \bsym{t} ]])$,
and write it as 
$\beta = \sum_{\bsym{k}}\, \beta_{\bsym{k}} 
\bsym{\partial}_{\bsym{k}}$,
where the sum is over the multi-indices
$\bsym{k} = (k_1 < \cdots < k_q)$, 
and $\beta_{\bsym{k}} \in \Omega_{\opn{Coor} X}
\, \what{\otimes} \, \K [[ \bsym{t} ]]$.
Then
$[\beta, \bsym{t}^{\bsym{k}}] = \beta_{\bsym{k}}$.
We see that $\beta = 0$ iff 
$[\beta, \bsym{t}^{\bsym{k}}] = 0$
for all such $\bsym{k}$. 
Therefore is suffices to prove that
$[\nabla(\bsym{\partial}_{\bsym{j}}), \bsym{t}^{\bsym{k}}] =
[\opn{ad}(\omega)(\bsym{\partial}_{\bsym{j}}),  
\bsym{t}^{\bsym{k}}]$
for any $\bsym{k}$. Now 
$[\bsym{\partial}_{\bsym{j}}, \bsym{t}^{\bsym{k}}] \in \K$
(it is $0$ or $1$). Because $\nabla$ is a $\K$-linear
derivation, we have
\[ 0 = \nabla([\bsym{\partial}_{\bsym{j}}, \bsym{t}^{\bsym{k}}])
= [\nabla(\bsym{\partial}_{\bsym{j}}), \bsym{t}^{\bsym{k}}] \pm
[\bsym{\partial}_{\bsym{j}}, \nabla(\bsym{t}^{\bsym{k}})]  . \]
Likewise
\[ 0 = \opn{ad}(\omega)(
[\bsym{\partial}_{\bsym{j}}, \bsym{t}^{\bsym{k}}])
= [\opn{ad}(\omega)(\bsym{\partial}_{\bsym{j}}), 
\bsym{t}^{\bsym{k}}] \pm
[\bsym{\partial}_{\bsym{j}}, 
\opn{ad}(\omega)(\bsym{t}^{\bsym{k}})]  . \]
And by definition of $\omega$ we have
\[ \opn{ad}(\omega)(\bsym{t}^{\bsym{k}}) = 
\sum_l \, \nabla(t_l) \, \partial_l(\bsym{t}^{\bsym{k}}) = 
\nabla(\bsym{t}^{\bsym{k}}) . \]

The case $\mcal{G} = \mcal{D}^{}_{\mrm{poly}}$ is handled 
similarly, using the basis 
\[ \bsym{\partial}_{\bsym{j}_1, \ldots, \bsym{j}_q} :=
(\smfrac{\partial}{\partial \bsym{t}})^{\bsym{j}_1} \otimes
\cdots \otimes 
(\smfrac{\partial}{\partial \bsym{t}})^{\bsym{j}_{q}}
\in  \mcal{D}^{q-1}_{\mrm{poly}}(\K[[\bsym{t}]]) , \]
see equation \cite[equation 4.3]{Ye2}.
\end{proof}
 
\begin{prop} \label{prop2.2}
The form $\omega_{\mrm{MC}}$ satisfies the identity
\[ \mrm{d}_{\mrm{for}}(\omega_{\mrm{MC}}) + \smfrac{1}{2}
[\omega_{\mrm{MC}}, \omega_{\mrm{MC}}] = 0 ; \]
namely it is a solution of the MC equation in the DG Lie algebra 
\linebreak 
$\Omega_{\opn{Coor} X} \, \what{\otimes} \, 
\mcal{T}_{\mrm{poly}}(\K [[ \bsym{t} ]])$
with differential $\mrm{d}_{\mrm{for}}$.
\end{prop}

\begin{proof}
Let's write
$\omega := \omega_{\mrm{MC}}$,
$\mrm{d} := \mrm{d}_{\mrm{for}}$,
$\nabla := \nabla_{\mcal{P}}$,
$\beta := \mrm{d}(\omega) + \smfrac{1}{2}[\omega, \omega]$
and
$\mfrak{g} := \Omega_{\opn{Coor} X} \, \what{\otimes} \, 
\mcal{T}_{\mrm{poly}}(\K [[ \bsym{t} ]])$.
As explained in the proof of the previous proposition, 
it suffices to show that $[\beta, t_i] = 0$ for all $i$. 

By definition of $\omega$ we have
\begin{equation} \label{eqn2.1}
[\omega, t_i] = \nabla_{\mcal{P}}(t_i) . 
\end{equation}
Next we use the fact that $\mrm{d}$ is an odd
derivation of $\mfrak{g}$ to obtain
\[ \mrm{d}([\omega, t_i]) = [\mrm{d}(\omega), t_i] -  
[\omega, \mrm{d}(t_i)] . \]
But $\mrm{d}(t_i) = 0$, so
\begin{equation} \label{eqn2.2}
[\mrm{d}(\omega), t_i] = 
\mrm{d}(\nabla_{\mcal{P}}(t_i)) . 
\end{equation}

The graded Jacobi identity in $\mfrak{g}$ tells us that
\[ [[\omega, \omega], t_i] + [[t_i, \omega], \omega] -
[[\omega, t_i], \omega] = 0 . \]
Hence 
$[[\omega, \omega], t_i] = 2 [\omega, [\omega, t_i]]$,
and plugging in (\ref{eqn2.1}) we arrive at
\begin{equation} \label{eqn2.3}
\smfrac{1}{2} [[\omega, \omega], t_i] = 
\opn{ad}(\omega)(\nabla_{\mcal{P}}(t_i)) . 
\end{equation}
Finally, combining (\ref{eqn2.2}), (\ref{eqn2.3}) and Proposition
\ref{prop2.3} we get
\[ \begin{aligned}
& [\beta, t_i] = [\mrm{d}(\omega), t_i] + 
\smfrac{1}{2} [[\omega, \omega], t_i] =
\mrm{d}(\nabla_{\mcal{P}}(t_i)) + 
\opn{ad}(\omega)(\nabla_{\mcal{P}}(t_i)) \\
& \qquad = (\mrm{d} + \opn{ad}(\omega))(\nabla_{\mcal{P}}(t_i))
= \nabla_{\mcal{P}}(\nabla_{\mcal{P}}(t_i)) = 0 . 
\end{aligned} \]
\end{proof}

According to Theorem \ref{thm1.1}(2) the group $G(\K)$ of 
$\K$-algebra automorphisms of 
$\K[[\bsym{t}]]$ acts on the bundle $\opn{Coor} X$. Therefore for 
any open set $U \subset X$ this group acts on the algebra
$\Gamma \bigl( \pi_{\mrm{coor}}^{-1}(U), 
\pi_{\mrm{coor}}^{\what{*}} \mcal{P}_X \bigr)$.
More generally, let $\mcal{G}$ denote either 
$\mcal{T}^{}_{\mrm{poly}}$ or $\mcal{D}^{}_{\mrm{poly}}$.
Let's introduce the temporary notation
\[ \mfrak{h}(U, \mcal{G}) :=
\Gamma \bigl( \pi_{\mrm{coor}}^{-1}(U), 
\Omega_{\opn{Coor} X} \, 
\what{\otimes}_{\mcal{O}_{\opn{Coor} X}} \,
\pi_{\mrm{coor}}^{\what{*}} (\mcal{P}_X \otimes_{\mcal{O}_{X}}
\mcal{G}) \bigr)  \]
and
\[ \mfrak{h}'(U, \mcal{G}) :=
\Gamma \bigl( \pi_{\mrm{coor}}^{-1}(U), 
\Omega_{\opn{Coor} X} \hatotimes{}
\mcal{G}(\K[[\bsym{t}]]) \bigr) . \]
These are graded Lie algebras. The group $G(\K)$ acts on 
$\mfrak{h}(U, \mcal{G})$ via its geometric action on 
$\opn{Coor} X$. On the other hand there is 
an action of $G(\K)$ on $\mfrak{h}'(U, \mcal{G})$
via its action on $\Gamma \bigl( \pi_{\mrm{coor}}^{-1}(U), 
\Omega_{\opn{Coor} X} \bigr)$
and on $\K[[\bsym{t}]]$. 

\begin{prop} \label{prop5.3}
The canonical isomorphism 
$\mfrak{h}(U, \mcal{G}) \cong \mfrak{h}'(U, \mcal{G})$ 
of Theorem \tup{\ref{thm5.3}} is $G(\K)$-equivariant.
\end{prop}

\begin{proof} 
By Theorem \ref{thm1.1}(2) the algebra  isomorphism
\[ \Gamma \bigl(\pi_{\mrm{coor}}^{-1}(U), 
\mcal{O}_{\opn{Coor} X}[[\bsym{t}]] \bigr) \cong
\Gamma \bigl(\pi_{\mrm{coor}}^{-1}(U), 
\pi_{\mrm{coor}}^{\what{*}} \mcal{P}_X \bigr) \]
of Corollary \ref{cor4.1} is $G(\K)$-equivariant. 
Tracing the isomorphisms used in the proof of 
Theorem \ref{thm5.3} we deduce the same for the isomorphism
$\mfrak{h}(U, \mcal{G}) \cong \mfrak{h}'(U, \mcal{G})$.
\end{proof}

We view $\omega_{\mrm{MC}}$ as an element of
$\mfrak{h}(X, \mcal{T}^{}_{\mrm{poly}}) 
\cong \mfrak{h}'(X, \mcal{T}^{}_{\mrm{poly}})$.
Due to 
Proposition \ref{prop5.3} we can talk about {\em the} action of 
$G(\K)$ on $\omega_{\mrm{MC}}$. 
Recall that $\mrm{GL}_n(\K)$ sits inside $G(\K)$ as the group of 
linear changes of coordinates. 

\begin{prop} \label{prop5.2}
The element $\omega_{\mrm{MC}}$ is $\mrm{GL}_n(\K)$-invariant.
\end{prop}

\begin{proof}
Since the Grothendieck connection $\nabla_{\mcal{P}}$ on
$\mfrak{h}(X, \mcal{T}^{}_{\mrm{poly}})$ 
is induced from $X$, it commutes with 
the action of $G(\K)$. Hence in particular 
$g(\nabla_{\mcal{P}}(t_i)) = \nabla_{\mcal{P}}(g(t_i))$
for any $g \in \mrm{GL}_n(\K)$ and $i \in \{1, \ldots, n\}$.

Fix such a matrix $g = [g_{i,j}]$. So $g_{i,j} \in \K$ and
$g(t_i) = \sum_{j=1}^n g_{i,j} t_j$. 
Let $h = [h_{i,j}] := (g^{-1})^{\mrm{t}}$, 
the transpose inverse matrix. Then in the induced action of 
$\mrm{GL}_n(\K)$ on 
$\mcal{T}^{0}_{\mrm{poly}}(\K[[\bsym{t}]]) = 
\mcal{T}(\K[[\bsym{t}]])$
we have 
$g(\frac{\partial}{\partial t_i}) =
\sum_{j=1}^n h_{i,j} \frac{\partial}{\partial t_j}$.
Thus
\[ \begin{aligned}
g(\omega_{\mrm{MC}}) & = 
g \Bigl( \sum_i
\nabla_{\mcal{P}}(t_i) \cdot \smfrac{\partial}{\partial t_i}
\Bigr) =
\sum_i g \bigl( \nabla_{\mcal{P}}(t_i) \bigr) 
\cdot g (\smfrac{\partial}{\partial t_i}) \\
& = \sum_i \nabla_{\mcal{P}} \bigl( g(t_i) \bigr) 
\cdot g(\smfrac{\partial}{\partial t_i}) =
\sum_{i,j,k} g_{i,j} h_{i,k} \bigl( \nabla_{\mcal{P}}(t_j) \bigr) 
\cdot \smfrac{\partial}{\partial t_k} \\
& = \sum_{j} \nabla_{\mcal{P}}(t_j) 
\cdot \smfrac{\partial}{\partial t_j} =  \omega_{\mrm{MC}} . 
\end{aligned} \]
\end{proof}

\begin{rem} \label{rem4.1}
The adjoint of $\omega_{\mrm{MC}}$ is an element of
$\Gamma \big(\opn{Coor} X, \mcal{T}_{\opn{Coor} X} 
\, \what{\otimes}\, \what{\Omega}^1_{\K[[\bsym{t}]] / \K}
\big),$
and it gives rise to a Lie algebra homomorphism
$\mcal{T}_{\K[[\bsym{t}]]} \to 
\Gamma(\opn{Coor} X, \mcal{T}_{\opn{Coor} X})$. 
In this way $\mcal{T}_{\K[[\bsym{t}]]}$ acts infinitesimally 
on $\opn{Coor} X$. Now inside 
$\mcal{T}_{\K[[\bsym{t}]]} = \boplus_{i = 1}^n 
\K[[\bsym{t}]] \frac{\partial}{\partial t_i}$
there is a subalgebra 
$\mfrak{g} := \boplus_{i,j} 
\K[[\bsym{t}]] \, t_i \frac{\partial}{\partial t_j}$.
The Lie algebra $\mfrak{g}$ is the Lie algebra of the 
pro-algebraic group 
$G = \opn{Aut}(\K[[\bsym{t}]])$, the group of $\K$-algebra 
automorphisms of $\K[[\bsym{t}]]$. The infinitesimal action of
$\mfrak{g}$ on $\opn{Coor} X$ is the differential 
of the action of $G$ on $\opn{Coor} X$
(cf.\ Theorem \ref{thm1.1}). 
The action of $\mcal{T}_{\K[[\bsym{t}]]}$ on $\opn{Coor} X$
is the main feature of the Gelfand-Kazhdan formal geometry. 
However we do not use this action 
(at least not directly) in our paper.
\end{rem}

\section{Review of Mixed Resolutions}
\label{sec6}

As always $\K$ is a field of characteristic $0$. In this section 
we review the constructions and results of the paper \cite{Ye3}. 

Let $\bsym{\Delta}$ denote the category with set of
objects the natural numbers. For any $p, q \in \mbb{N}$ the 
set of morphisms in $\bsym{\Delta}$ from $p$ to $q$ is the set
$\bsym{\Delta}_p^q$ of order preserving functions
$\alpha : \{ 0, \ldots, p \} \to \{ 0, \ldots, q \}$.
Recall that a cosimplicial object in some category
$\cat{C}$ is a functor $C : \bsym{\Delta} \to \cat{C}$.
Usually one writes $C^p$ instead of $C(p)$, and refers to the 
sequence $\{ C^p \}_{p \geq 0}$ as a cosimplicial object (the 
morphisms remaining implicit). The category of 
cosimplicial objects in $\cat{C}$ is denoted by
$\bsym{\Delta} \msf{C}$. 

We are interested in cosimplicial dir-inv $\K$-modules, i.e.\ in 
objects $M = \{ M^p \}_{p \geq 0}$ in 
$\bsym{\Delta} \cat{Dir} \cat{Inv} \cat{Mod} \K$.
As explained in \cite{Ye3}, there is a functor
\[ \what{\til{\mrm{N}}} :
\bsym{\Delta} \cat{Dir} \cat{Inv} \cat{Mod} \K \to
\cat{DGMod} \K , \]
the latter being the category of complexes of $\K$-modules. This 
is the {\em complete Thom-Sullivan normalization} functor, which 
is a generalization of constructions in \cite{HS} and \cite{HY}. 
By definition there is an embedding
\[ \what{\til{\mrm{N}}}{}^{q} M \subset
\prod_{i=0}^{\infty} \left( 
\Omega^{q}(\bsym{\Delta}^{i}_{\K}) \, \what{\otimes} \, 
M^{i} \right) . \]
Here
$\bsym{\Delta}^i_{\K} := \opn{Spec} \K[t_0, \ldots, t_i] / 
(t_0 + \cdots + t_i - 1)$ 
is the $i$-dimensional geometric simplex, and
$\Omega^{q}(\bsym{\Delta}^{i}_{\K}) := 
\Gamma(\bsym{\Delta}^i_{\K}, 
\Omega^q_{\bsym{\Delta}^i_{\K}})$
is a discrete inv $\K$-module. The differential 
$\partial : \what{\til{\mrm{N}}}{}^q M \to
\what{\til{\mrm{N}}}{}^{q+1} M$
is induced from the de Rham differentials
$\d : \Omega^{q}(\bsym{\Delta}^{i}_{\K}) \to
\Omega^{q+1}(\bsym{\Delta}^{i}_{\K})$.

Let $X$ be a separated smooth irreducible $n$-dimensional 
$\K$-scheme. Choose an affine open covering
$\bsym{U} = \{ U_{(0)}, \ldots, U_{(m)} \}$
of $X$. Given 
$\bsym{i} = (i_0, \ldots, i_q) \in \bsym{\Delta}_q^m$
let
$U_{\bsym{i}} := U_{(i_0)} \cap \cdots \cap U_{(i_m)}$,
and let
$g_{\bsym{i}} : U_{\bsym{i}} \to X$
be the inclusion. For a sheaf $\mcal{M}$ on $X$ we write
\[ \mrm{C}^q(\bsym{U}, \mcal{M}) :=
\prod_{\bsym{i} \in \bsym{\Delta}_q^m} 
g_{\bsym{i} *}\, g_{\bsym{i}}^{-1}\, \mcal{M} . \]
The sequence
$\{ \mrm{C}^q(\bsym{U}, \mcal{M}) \}_{q \geq 0}$
is then a cosimplicial sheaf on $X$. This is a variant of the 
\v{C}ech resolution of $\mcal{M}$.

Suppose $\mcal{M}$ is a dir-inv $\K_X$-module, i.e.\ a sheaf of 
$\K$-modules on $X$ with a dir-inv structure. For any open set
$V \subset X$ we then have a cosimplicial dir-inv $\K$-module 
$\bigl\{ \Gamma \bigl(V, \mrm{C}^q(\bsym{U}, \mcal{M}) \bigr) 
\bigr\}_{q \geq 0}$. Applying the functor 
$\what{\til{\mrm{N}}}{}^q$ to it we obtain a $\K$-module
$\what{\til{\mrm{N}}}{}^q
\Gamma \bigl(V, \mrm{C}(\bsym{U}, \mcal{M}) \bigr)$.
It turns out that the presheaf
$V \mapsto \what{\til{\mrm{N}}}{}^{q}
\Gamma \bigl(V, \mrm{C}(\bsym{U}, \mcal{M}) \bigr)$
is a sheaf, and we denote it by
$\what{\til{\mrm{N}}}{}^{q} \mrm{C}(\bsym{U}, \mcal{M})$.
So there is a functor
\[ \what{\til{\mrm{N}}} \mrm{C}(\bsym{U}, -) :
\cat{Dir} \cat{Inv} \cat{Mod} \K_X \to
\cat{DGMod} \K_X , \]
and there is a functorial homomorphism 
$\mcal{M} \to 
\what{\til{\mrm{N}}} \mrm{C}(\bsym{U}, \mcal{M})$.
If $\mcal{M}$ is a complete dir-inv module then according to 
\cite[Theorem 3.7]{Ye3} the homomorphism 
$\mcal{M} \to 
\what{\til{\mrm{N}}} \mrm{C}(\bsym{U}, \mcal{M})$
is in fact a quasi-isomorphism. We call 
$\what{\til{\mrm{N}}}{}^{} \mrm{C}(\bsym{U}, \mcal{M})$
the {\em commutative \v{C}ech resolution} of $\mcal{M}$, since
$\what{\til{\mrm{N}}}{}^{} \mrm{C}(\bsym{U}, \mcal{O}_X)$
is a super-commutative DG algebra.

Now suppose $\mcal{M}$ is a quasi-coherent $\mcal{O}_X$-module. 
Let $p$ be some natural number. Then
$\Omega^p_X \otimes_{\mcal{O}_{X}} \mcal{P}_X 
\otimes_{\mcal{O}_{X}} \mcal{M}$
is a complete dir-inv $\mcal{P}_X$-module with the 
$\mcal{I}_X$-adic dir-inv structure. 
Define
\[ \opn{Mix}^{p,q}_{\bsym{U}}(\mcal{M}) :=
\what{\til{\mrm{N}}}{}^{q} \mrm{C}(\bsym{U}, 
\Omega^p_X \otimes_{\mcal{O}_{X}} \mcal{P}_X 
\otimes_{\mcal{O}_{X}} \mcal{M}) . \]
This is a sheaf on $X$, and there is an embedding of sheaves
\begin{equation} \label{eqn5.3}
\opn{Mix}^{p,q}_{\bsym{U}}(\mcal{M})
\subset 
{\displaystyle \prod_{j \in \mbb{N}}}\ \
{\displaystyle \prod_{\bsym{i} \in 
\bsym{\Delta}_j^{m}}} \,
g_{\bsym{i} *}\, g_{\bsym{i}}^{-1}\,
\big( \Omega^{q}(\bsym{\Delta}^j_{\K})
\, \what{\otimes} \, (\Omega^{p}_X \otimes_{\mcal{O}_{X}}
\mcal{P}_X \otimes_{\mcal{O}_{X}} \mcal{M}) \big) .
\end{equation}
In addition to the differential 
$\partial : \opn{Mix}^{p,q}_{\bsym{U}}(\mcal{M}) \to
\opn{Mix}^{p,q+1}_{\bsym{U}}(\mcal{M})$
there is a second differential
$\nabla_{\mcal{P}} : \opn{Mix}^{p,q}_{\bsym{U}}(\mcal{M}) \to
\opn{Mix}^{p+1,q}_{\bsym{U}}(\mcal{M})$
coming from the connection $\nabla_{\mcal{P}}$ of Definition
\ref{dfn5.2}. We now totalize 
\[ \opn{Mix}^{i}_{\bsym{U}}(\mcal{M}) := 
\boplus_{p+q = i} \opn{Mix}^{p,q}_{\bsym{U}}(\mcal{M})  \]
and let
$\d_{\mrm{mix}} := \partial + (-1)^q \nabla_{\mcal{P}}$. 
This is the {\em mixed resolution} of $\mcal{M}$, which is a 
functor
\[ \opn{Mix}^{}_{\bsym{U}} :
\cat{QCoh} \mcal{O}_{X} \to \cat{DGMod} \K_X . \]

\begin{thm}[{\cite[Theorem 4.14]{Ye3}}] \label{thm6.6}
Let $\mcal{M}$ be a quasi-coherent $\mcal{O}_{X}$-module.
\begin{enumerate}
\item There is a functorial quasi-isomorphism
$\mcal{M} \to \opn{Mix}^{}_{\bsym{U}}(\mcal{M})$.
\item There is a functorial isomorphism 
$\Gamma \big( X, \opn{Mix}^{}_{\bsym{U}}(\mcal{M}) \big) \cong
\mrm{R} \Gamma(X, \mcal{M})$ 
in \linebreak $\msf{D}(\cat{Mod} \K)$.
\end{enumerate}
\end{thm}

Of course the functor 
$\opn{Mix}^{}_{\bsym{U}}$ can be extended to bounded below
complexes of quasi-coherent $\mcal{O}_{X}$-modules, by totalizing.

The sheaves of DG Lie algebras 
$\mcal{T}^{}_{\mrm{poly}, X}$ and
$\mcal{D}^{}_{\mrm{poly}, X}$ are bounded below complexes of 
quasi-coherent $\mcal{O}_{X}$-modules, so the above theorem 
applies to them. In addition we have:

\begin{prop} \label{prop6.2}
Let $\mcal{G}_X$ stand for either $\mcal{T}^{}_{\mrm{poly}, X}$ or 
$\mcal{D}^{}_{\mrm{poly}, X}$. Then 
$\opn{Mix}^{}_{\bsym{U}}(\mcal{G}_X)$
is a sheaf of DG Lie algebras, with differential 
\[ \d_{\mrm{mix}} + (-1)^i \d_{\mcal{G}} :
\opn{Mix}^{i}_{\bsym{U}}(\mcal{G}^j_X) \to
\opn{Mix}^{i+1}_{\bsym{U}}(\mcal{G}^j_X) \oplus
\opn{Mix}^{i}_{\bsym{U}}(\mcal{G}^{j+1} _X) . \]
The quasi-isomorphism 
$\mcal{G}_X \to \opn{Mix}^{}_{\bsym{U}}(\mcal{G}_X)$
of Theorem \tup{\ref{thm6.6}(1)} 
is a homomorphism of DG Lie algebras.
\end{prop}

\begin{proof}
By Proposition \ref{prop5.1} the sheaf
$\Omega^{}_X \otimes_{\mcal{O}_{X}} \mcal{P}_X
\otimes_{\mcal{O}_{X}} \mcal{G}_X$
is a DG $\Omega_X$-module Lie algebra in 
$\cat{Dir} \cat{Inv} \cat{Mod} \K_X$. Now use 
\cite[Proposition 5.5]{Ye3}.
\end{proof}

Suppose $\pi : Z \to X$ is some morphism of schemes (possibly of
infinite type). A {\em simplicial section} of $\pi$ based on the 
covering $\bsym{U}$ is a collection of $X$-morphisms
\[ \bsym{\sigma} = \{ \sigma_{\bsym{i}} :
\bsym{\Delta}^{q}_{\K} \times U_{\bsym{i}} \to Z \} \]
indexed by $\bsym{i} \in \bsym{\Delta}^{m}_{q}$,
$q \in \mbb{N}$, which satisfies the simplicial relations (see 
\cite[Definition 6.1]{Ye3}). 

The sheaf $\Omega^p_Z$ is considered as a discrete inv $\K_Z$-module,
and $\Omega_Z = \boplus_p \Omega^p_Z$
has the $\boplus$ dir-inv structure. Given a quasi-coherent 
$\mcal{O}_{X}$-module $\mcal{M}$ the graded sheaf 
$\Omega^{}_{Z} \, 
\what{\otimes}_{\mcal{O}_{Z}} \, \pi^{\what{*}}\, 
(\mcal{P}_X \otimes_{\mcal{O}_X} \mcal{M})$
is then a DG $\Omega_Z$-module in 
$\cat{Dir} \cat{Inv} \cat{Mod} \K_Z$, 
with differential $\nabla_{\mcal{P}}$. See Section \ref{sec2}.

Let $A$ be an associative unital super-commutative DG $\K$-algebra. 
Consider a homogeneous $A$-multilinear function
$\phi : M_1 \times \cdots \times M_r \to N$,
where $M_1, \ldots, M_r, N$ are DG $A$-modules. There is an 
operation of composition for such functions: given  
functions $\psi_i : \prod_j L_{i,j} \to M_i$ 
the composition is
$\phi \circ (\psi_1 \times \cdots \times \psi_r)
: \prod_{i,j} L_{i,j} \to N$.
There is also a summation operation: if 
$\phi_j : \prod_i M_i \to N$ are homogeneous of equal degree then 
so is their sum $\sum_j \phi_j$. 
Finally let
$\phi \circ \d : \prod_i M_i \to N$ be the homogeneous 
$\K$-multilinear function
\[ (\phi \circ \d)(m_1, \ldots, m_r) := \sum_{i=1}^r
\pm \phi(m_1, \ldots, \d(m_i), \ldots, m_r) \]
with Koszul signs.

\begin{thm}[{\cite[Theorem 6.3]{Ye3}}] \label{thm6.8}
Suppose $\bsym{\sigma}$ is simplicial section of $\pi : Z \to X$ 
based on $\bsym{U}$. Let 
$\mcal{M}_1, \ldots, \mcal{M}_r, \mcal{N}$ be quasi-coherent 
$\mcal{O}_{X}$-modules, and let
\[ \phi : \prod_{i = 1}^r 
\bigl( \Omega_Z \hatotimes{\mcal{O}_Z} \pi^{\what{*}}
(\mcal{P}_X \otimes_{\mcal{O}_{X}} \mcal{M}_i) \big) 
\to \Omega_Z \hatotimes{\mcal{O}_Z} \pi^{\what{*}}
(\mcal{P}_X \otimes_{\mcal{O}_{X}} \mcal{N}) \]
be a continuous $\Omega_Z$-multilinear sheaf 
morphism on $Z$ of degree $k$. Then there is an induced 
$\K$-multilinear sheaf morphism of degree $k$
\[ \bsym{\sigma}^*(\phi) :
\prod_{i = 1}^r \opn{Mix}_{\bsym{U}}(\mcal{M}_i) \to 
\opn{Mix}_{\bsym{U}}(\mcal{N})  \]
on $X$ with the following properties.
\begin{enumerate}
\rmitem{i} The assignment $\phi \mapsto \bsym{\sigma}^*(\phi)$ 
respects the operations of composition and summation.
\rmitem{ii} If $\phi = \pi^{\what{*}}(\phi_0)$ for some continuous 
$\Omega_X$-multilinear morphism
\[ \phi_0  : \prod_{i=1}^r
\big( \Omega^{}_{X} \otimes_{\mcal{O}_{X}} 
\mcal{P}_X \otimes_{\mcal{O}_X} \mcal{M}_i \big) 
\to \Omega^{}_{X} \otimes_{\mcal{O}_{X}}
\mcal{P}_X \otimes_{\mcal{O}_X} \mcal{N} \] 
then 
$\bsym{\sigma}^*(\phi) = 
\what{\til{\mrm{N}}} \mrm{C}(\bsym{U}, \phi_0)$.
\rmitem{iii} Suppose $\psi$ is another such 
$\Omega_Z$-multilinear 
sheaf morphism of degree $k+1$, and the equation
$\nabla_{\mcal{P}} \circ \phi - (-1)^k \phi \circ \nabla_{\mcal{P}}
= \psi$
holds. Then
\[ \d_{\mrm{mix}} \circ \bsym{\sigma}^*(\phi) - (-1)^k
\bsym{\sigma}^*(\phi) \circ \d_{\mrm{mix}} =
\bsym{\sigma}^*(\psi) . \]
\end{enumerate}
\end{thm}

We are interested in the bundle 
$\pi_{\mrm{lcc}} : \opn{LCC} X \to X$.

\begin{thm} \label{thm6.7}
Assume each affine open set $U_{(i)}$ admits an \'etale morphism 
to $\mbf{A}^n_{\K}$. Then there exist sections 
$\sigma_{(i)} : U_{(i)} \to \opn{LCC} X$, and furthermore 
they extend to a simplicial section 
$\bsym{\sigma}$ of $\pi_{\mrm{lcc}} : \opn{LCC} X \to X$.
\end{thm}

\begin{proof}
By Theorem \ref{thm1.1}, $\pi_{\mrm{coor}} : \opn{Coor} X \to X$
is a locally trivial $G$-torsor over
$X$, and $G = \opn{GL}_{n, \K} \ltimes N$, 
where $N$ is a pro-unipotent group. By definition
$\opn{LCC} X = \opn{Coor} X / \opn{GL}_{n, \K}$.
According to Example \ref{exa1.1} and Corollary \ref{cor4.1}, 
for any $i$ there is a section 
$\sigma_{(i)} : U_{(i)} \to \opn{Coor} X$. 
Now use \cite[Theorem 2.2]{Ye4}.
\end{proof}

\begin{figure}
\includegraphics{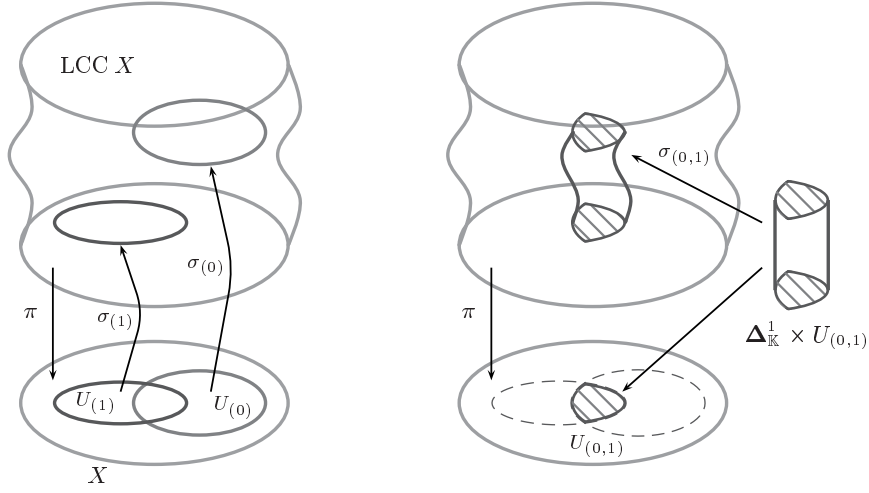}
\caption{Simplicial sections, $q = 1$. 
We start with sections over two open sets $U_{(0)}$ and $U_{(1)}$ 
in the left diagram; and we pass to a simplicial 
section $\sigma_{(0, 1)}$ on the right.} 
\label{fig1}
\end{figure}

Here is the idea behind the proof of \cite[Theorem 2.2]{Ye4}.
There is an averaging process for sections of torsors under 
unipotent groups. The bundle $\opn{LCC} X$ is ``almost'' a torsor 
under the pro-unipotent group $N$. Given a multi-index
$\bsym{i} = (i_0, \ldots, i_q)$ the morphism
$\sigma_{\bsym{i}} :
\bsym{\Delta}^{q}_{\K} \times U_{\bsym{i}} \to \opn{LCC} X$
is then a family of weighted averages of the sections
$\sigma_{(i_0)}, \ldots, \sigma_{(i_q)} : 
U_{\bsym{i}} \to \opn{LCC} X$, 
parameterized by the simplex $\bsym{\Delta}^{q}_{\K}$.
See Figure \ref{fig1} for an illustration. 

\section{The Global $\mrm{L}_{\infty}$ Quasi-isomorphism}

In this section we prove the main result of the paper.
Here again $X$ is a smooth irreducible
separated $n$-dimensional 
scheme over the field $\K$, and also $\mbb{R} \subset \K$. 

Fix an open covering 
$\bsym{U} = \{ U_{(0)}, \ldots, U_{(m)} \}$ 
of the scheme $X$ consisting of affine open sets, each admitting 
an \'etale morphism $U_{(i)} \to \mbf{A}^n_{\K}$. For every $i$ let
$\sigma_{(i)} : U_{(i)} \to \opn{LCC} X$ be the corresponding 
section of $\pi_{\mrm{lcc}} : \opn{LCC} X \to X$, 
and let $\bsym{\sigma}$ be the 
resulting simplicial section (see Theorem \ref{thm6.7}).

Let $\mcal{M}$ be a bounded below complex of 
quasi-coherent $\mcal{O}_{X}$-modules.
The mixed resolution $\opn{Mix}^{}_{\bsym{U}}(\mcal{M})$ 
was defined in Section \ref{sec6}. For any integer $i$ let
$\mrm{G}^i \opn{Mix}^{}_{\bsym{U}}(\mcal{M}) :=
\boplus_{j = i}^{\infty} \opn{Mix}^{j}_{\bsym{U}}(\mcal{M})$,
so 
$\{ \mrm{G}^i \opn{Mix}^{}_{\bsym{U}}(\mcal{M}) \}_{i \in \mbb{Z}}$
is a descending filtration of 
$\opn{Mix}^{}_{\bsym{U}}(\mcal{M})$
by subcomplexes, with 
$\mrm{G}^i \opn{Mix}^{}_{\bsym{U}}(\mcal{M}) = 
\opn{Mix}^{}_{\bsym{U}}(\mcal{M})$
for $i \ll 0$ and
$\bigcap_{i}\, \mrm{G}^i \opn{Mix}^{}_{\bsym{U}}(\mcal{M}) = 0$.
Let 
\[ \opn{gr}_{\mrm{G}}^i \opn{Mix}^{}_{\bsym{U}}(\mcal{M}) :=
\mrm{G}^i \opn{Mix}^{}_{\bsym{U}}(\mcal{M}) \ / \
\mrm{G}^{i+1} \opn{Mix}^{}_{\bsym{U}}(\mcal{M}) \]
and
$\opn{gr}_{\mrm{G}} \opn{Mix}^{}_{\bsym{U}}(\mcal{M}) :=
\boplus_{i}\, \opn{gr}_{\mrm{G}}^i 
\opn{Mix}^{}_{\bsym{U}}(\mcal{M})$.

By Proposition \ref{prop6.2}, if $\mcal{G}_X$ is either 
$\mcal{T}^{}_{\mrm{poly}, X}$ or 
$\mcal{D}^{}_{\mrm{poly}, X}$, then 
$\opn{Mix}^{}_{\bsym{U}}(\mcal{G}_X)$ is a sheaf of DG Lie 
algebras on $X$, and
$\mcal{G}_X \to \opn{Mix}^{}_{\bsym{U}}(\mcal{G}_X)$ 
is a DG Lie algebra quasi-isomorphism. 

Note that if 
$\phi : \opn{Mix}^{}_{\bsym{U}}(\mcal{M}) \to
\opn{Mix}^{}_{\bsym{U}}(\mcal{N})$
is a homomorphism of complexes that respects the filtration
$\{ \mrm{G}^i \opn{Mix}^{}_{\bsym{U}} \}$,
then there exists an induced homomorphism of complexes
\[ \opn{gr}_{\mrm{G}}(\phi) : 
\opn{gr}_{\mrm{G}} \opn{Mix}^{}_{\bsym{U}}(\mcal{M}) \to
\opn{gr}_{\mrm{G}} \opn{Mix}^{}_{\bsym{U}}(\mcal{N}) . \]

Suppose $\mcal{G}$ and $\mcal{H}$ are sheaves of DG Lie algebras 
on a topological space $Y$. An $\mrm{L}_{\infty}$ morphism
$\Psi : \mcal{G} \to \mcal{H}$ is a sequence of sheaf morphisms
$\psi_j : \sprod^j \mcal{G} \to \mcal{H}$, such that for every 
open set $V \subset Y$ the sequence 
$\{ \Gamma(V, \psi_j) \}_{j \geq 1}$ is an $\mrm{L}_{\infty}$ 
morphism
$\Gamma(V, \mcal{G}) \to \Gamma(V, \mcal{H})$. If 
$\psi_1 : \mcal{G} \to \mcal{H}$ is a 
quasi-isomorphism then $\Psi$ is called an $\mrm{L}_{\infty}$ 
quasi-morphism.

\begin{thm} \label{thm5.1}
Let $\bsym{U}$ and $\bsym{\sigma}$ be as above. Then there is an 
induced $\mrm{L}_{\infty}$ quasi-isomorphism
\[ \Psi_{\bsym{\sigma}} = 
\set{\Psi_{\bsym{\sigma}; j}}_{j \geq 1} :
\opn{Mix}^{}_{\bsym{U}}(\mcal{T}_{\mrm{poly}, X}) 
\to \opn{Mix}^{}_{\bsym{U}}(\mcal{D}_{\mrm{poly}, X}) . \]
The homomorphism $\Psi_{\bsym{\sigma}; 1}$ respects the filtration
$\{ \mrm{G}^i \opn{Mix}^{}_{\bsym{U}} \}$,
and 
\[ \opn{gr}_{\mrm{G}}(\Psi_{\bsym{\sigma}; 1}) = 
\opn{gr}_{\mrm{G}}(\opn{Mix}^{}_{\bsym{U}}(\mcal{U}_{1})) . \]
\end{thm}

\begin{proof}
Let $Y$ be some
$\K$-scheme, and denote by $\K_Y$ the constant sheaf. For any $p$ 
we view $\Omega_Y^p$ as a discrete inv $\K_Y$-module, and we put 
on $\Omega_Y = \boplus_{p \in \mbb{N}}\, \Omega_Y^p$ 
direct sum dir-inv structure. So $\Omega_Y$ is a discrete (and 
hence complete)
DG algebra in $\cat{Dir} \cat{Inv} \cat{Mod} \K_Y$.

We shall abbreviate
$\mcal{A} := \Omega_{\opn{Coor} X}$, so that 
$\mcal{A}^0 = \mcal{O}_{\opn{Coor} X}$ etc.
As explained above, $\mcal{A}$ is a DG algebra in 
$\cat{Dir} \cat{Inv} \cat{Mod} \K_{\opn{Coor} X}$,
with discrete (but not trivial) dir-inv module structure.

There are sheaves of DG Lie algebras 
$\mcal{A} \, \what{\otimes} \
\mcal{T}^{}_{\mrm{poly}}(\K[[\bsym{t}]])$
and 
$\mcal{A} \, \what{\otimes} \
\mcal{D}^{}_{\mrm{poly}}(\K[[\bsym{t}]])$
on the scheme $\opn{Coor} X$. The differentials are
$\mrm{d}_{\mrm{for}} = \mrm{d} \otimes \bsym{1}$
and $\mrm{d}_{\mrm{for}} + \bsym{1} \otimes \mrm{d}_{\mcal{D}}$
respectively. As explained just prior to Theorem \ref{thm3.3}, 
there is a continuous $\mcal{A}$-multilinear 
$\mrm{L}_{\infty}$ morphism
\[ \mcal{U}_{\mcal{A}} = 
\{ \mcal{U}_{\mcal{A} ;j} \}_{j \geq 1} :
\mcal{A} \, \what{\otimes} \
\mcal{T}^{}_{\mrm{poly}}(\K[[\bsym{t}]]) \to
\mcal{A} \, \what{\otimes} \
\mcal{D}^{}_{\mrm{poly}}(\K[[\bsym{t}]]) . \]

The MC form $\omega := \omega_{\mrm{MC}}$ is a global section of
$\mcal{A}^1 \, \what{\otimes} \
\mcal{T}^{0}_{\mrm{poly}}(\K[[\bsym{t}]])$
satisfying the MC equation in the DG Lie algebra
$\mcal{A} \, \what{\otimes} \
\mcal{T}^{}_{\mrm{poly}}(\K[[\bsym{t}]])$. 
See Proposition \ref{prop2.2}. 
According to Theorem  \ref{thm3.3} the global section 
$\omega' := \mcal{U}_{\mcal{A} ;1}(\omega) \in 
\mcal{A}^1 \hatotimes{} 
\mcal{D}^{0}_{\mrm{poly}}(\K[[\bsym{t}]])$
is a solution of the MC equation in the DG Lie algebra
$\mcal{A} \, \what{\otimes} \
\mcal{D}^{}_{\mrm{poly}}(\K[[\bsym{t}]])$,
and there is a continuous $\mcal{A}$-multilinear
$\mrm{L}_{\infty}$ morphism
\[ \mcal{U}_{\mcal{A}, \omega} = \{ \mcal{U}_{\mcal{A}, \omega ;j} 
\}_{j \geq 1} : 
\big( \mcal{A} \, \what{\otimes} \,
\mcal{T}^{}_{\mrm{poly}}(\K[[\bsym{t}]]) \big)_{\omega} \to 
\big( \mcal{A} \, \what{\otimes} \,
\mcal{D}^{}_{\mrm{poly}}(\K[[\bsym{t}]]) \big)_{\omega'} \]
between the twisted DG Lie algebras. 
The formula is
\begin{equation} \label{eqn9.6}
\mcal{U}_{\mcal{A}, \omega; j}(\gamma_1 \cdots \gamma_j) = 
\sum_{k \geq 0} \smfrac{1}{(j+k)!}
\mcal{U}_{\mcal{A} ;j + k} (\omega^k \gamma_1 \cdots \gamma_j)
\end{equation} 
for 
$\gamma_1, \ldots, \gamma_j
\in \mcal{A} \, \what{\otimes} \, 
\mcal{T}^{}_{\mrm{poly}}(\K[[\bsym{t}]])$.
The two twisted DG Lie algebras have differentials
$\mrm{d}_{\mrm{for}} + \opn{ad}(\omega)$
and
$\mrm{d}_{\mrm{for}} + \opn{ad}(\omega')
+ \bsym{1} \otimes \mrm{d}_{\mcal{D}}$
respectively. 

By Theorem \ref{thm5.3} (the universal Taylor expansions) 
there are canonical isomorphisms of graded Lie algebras
in $\cat{Dir} \cat{Inv} \cat{Mod} \K_{\opn{Coor} X}$
\[ \mcal{A} \, \what{\otimes} \,
\mcal{T}^{}_{\mrm{poly}}(\K[[\bsym{t}]]) \cong
\mcal{A}
\, \what{\otimes}_{\mcal{A}^0} \,
\pi_{\mrm{coor}}^{\what{*}}(\mcal{P}_X \otimes_{\mcal{O}_{X}}
\mcal{T}^{}_{\mrm{poly}, X}) \]
and
\[ \mcal{A} \, \what{\otimes} \,
\mcal{D}^{}_{\mrm{poly}}(\K[[\bsym{t}]]) \cong
\mcal{A} \, \what{\otimes}_{\mcal{A}^0} \,
\pi_{\mrm{coor}}^{\what{*}}(\mcal{P}_X \otimes_{\mcal{O}_{X}}
\mcal{D}^{}_{\mrm{poly}, X}) . \]
Proposition \ref{prop2.3} tells us that 
\[ \mrm{d}_{\mrm{for}} + 
\opn{ad}(\omega) = \nabla_{\mcal{P}}  \]
under these identifications. Therefore 
\[ \mcal{U}_{\mcal{A}, \omega} : 
\mcal{A} \, \what{\otimes}_{\mcal{A}^0} \,
\pi_{\mrm{coor}}^{\what{*}}(\mcal{P}_X \otimes_{\mcal{O}_{X}}
\mcal{T}^{}_{\mrm{poly}, X}) 
\to \mcal{A} \, \what{\otimes}_{\mcal{A}^0} \,
\pi_{\mrm{coor}}^{\what{*}}(\mcal{P}_X \otimes_{\mcal{O}_{X}}
\mcal{D}^{}_{\mrm{poly}, X}) \]
is a continuous $\mcal{A}$-multilinear
$\mrm{L}_{\infty}$ morphism between these DG Lie algebras, 
whose differentials are
$\nabla_{\mcal{P}}$ and
$\nabla_{\mcal{P}} + \bsym{1} \otimes \d_{\mcal{D}}$ 
respectively.

By Propositions \ref{prop5.3} and \ref{prop5.2} the form 
$\omega$ is $\mrm{GL}_n(\K)$-invariant. 
So according to Proposition \ref{prop3.2} each of the operators
$\mcal{U}_{\mcal{A}; j}$ and 
$\mcal{U}_{\mcal{A}, \omega; j}$ is $\mrm{GL}_n(\K)$-equivariant. 
We conclude that $\omega$ is a global section of
\[ \Omega^1_{\opn{LCC} X} 
\, \what{\otimes}_{\mcal{O}_{\opn{LCC} X}} \,
\pi_{\mrm{lcc}}^{\what{*}}(\mcal{P}_X \otimes_{\mcal{O}_{X}}
\mcal{T}^{0}_{\mrm{poly}, X})  , \]
and the operators $\mcal{U}_{\mcal{A} ;j}$ and
$\mcal{U}_{\mcal{A}, \omega ;j}$
descend to continuous $\Omega_{\opn{LCC} X}$-multilinear operators
\[ \begin{aligned}
& \mcal{U}_{\mcal{A} ;j}, \ \mcal{U}_{\mcal{A}, \omega ;j} : 
\prod\nolimits^j \big( \Omega^{}_{\opn{LCC} X}
\, \what{\otimes}_{\mcal{O}_{\opn{LCC} X}} \,
\pi_{\mrm{lcc}}^{\what{*}}(\mcal{P}_X \otimes_{\mcal{O}_{X}}
\mcal{T}^{}_{\mrm{poly}, X}) \big) \\
& \qquad \qquad \qquad \qquad \to
\Omega^{}_{\opn{LCC} X}
\, \what{\otimes}_{\mcal{O}_{\opn{LCC} X}} \,
\pi_{\mrm{lcc}}^{\what{*}}(\mcal{P}_X \otimes_{\mcal{O}_{X}}
\mcal{D}^{}_{\mrm{poly}, X})
\end{aligned} \]
satisfying formula (\ref{eqn9.6}). The sequence
$\mcal{U}_{\mcal{A}, \omega} = 
\{ \mcal{U}_{\mcal{A}, \omega; j} \}_{j \geq 1}$
is an $\mrm{L}_{\infty}$ morphism.

According to Theorem \ref{thm6.8} there are induced 
operators 
\[ \bsym{\sigma}^*(\mcal{U}_{\mcal{A}; j}), \,
\bsym{\sigma}^*(\mcal{U}_{\mcal{A}, \omega; j}) :
\prod\nolimits^j
\opn{Mix}^{}_{\bsym{U}}
(\mcal{T}^{}_{\mrm{poly}, X}) 
\to \opn{Mix}^{}_{\bsym{U}}
(\mcal{D}^{}_{\mrm{poly}, X}) . \]
The $\mrm{L}_{\infty}$ identities in Definition \ref{dfn6.2},
when applied to the $\mrm{L}_{\infty}$ morphism 
$\mcal{U}_{\mcal{A}, \omega}$, are of the form considered in 
Theorem \ref{thm6.8}(iii). 
Therefore these identities are preserved by $\bsym{\sigma}^*$,
and we conclude that the sequence 
$\set{ \bsym{\sigma}^*(\mcal{U}_{\mcal{A}, \omega; j})}
_{j \geq 1}$
is an $\mrm{L}_{\infty}$ morphism.
There's a global section 
$\bsym{\sigma}^*(\omega) 
\in \opn{Mix}^{1}_{\bsym{U}} (\mcal{T}^{0}_{\mrm{poly}, X})$, 
and the formula
\begin{equation} \label{eqn10.4}
\bsym{\sigma}^*(\mcal{U}_{\mcal{A}, \omega; j})
(\gamma_1 \cdots \gamma_j) 
= \sum_{k \geq 0} \smfrac{1}{(j+k)!}
\bsym{\sigma}^*(\mcal{U}_{\mcal{A}; j+k}) \big( 
\bsym{\sigma}^*(\omega)^k 
\gamma_1 \cdots \gamma_j \big) 
\end{equation}
holds for local sections
$\gamma_1, \ldots, \gamma_j \in 
\opn{Mix}^{}_{\bsym{U}}(\mcal{T}^{}_{\mrm{poly}, X})$.
This sum is finite, the number of nonzero terms in it depending on 
the bidegree of $\gamma_1 \cdots \gamma_j$. 
Indeed, if
$\gamma_1 \cdots \gamma_j \in
\opn{Mix}^{q}_{\bsym{U}}(\mcal{T}^{p}_{\mrm{poly}, X})$
then
\begin{equation} \label{eqn11.1}
\bsym{\sigma}^*(\mcal{U}_{\mcal{A}; j+k})
\big( \bsym{\sigma}^*(\omega)^k 
\gamma_1 \cdots \gamma_j \big) \in 
\opn{Mix}^{q+k}_{\bsym{U}}(\mcal{D}^{p+1-j-k}_{\mrm{poly}, X}) ,
\end{equation}
which is is zero for $k > p-j+2$; see proof of  
\cite[Theorem 3.23]{Ye2}.

Finally we define
$\Psi_{\bsym{\sigma}; j} :=
\bsym{\sigma}^*(\mcal{U}_{\mcal{A}, \omega; j})$. 
The collection
$\Psi_{\bsym{\sigma}}  = 
\{ \Psi_{j, \bsym{\sigma}} \}_{j = 1}^{\infty}$
is then an $\mrm{L}_{\infty}$ morphism.
 From equation (\ref{eqn11.1}) we see that
$\Psi_{\bsym{\sigma}; 1}$ respects the filtration
$\{ \mrm{G}^i \opn{Mix}^{}_{\bsym{U}} \}$, 
and according to equation (\ref{eqn10.4}) we see that  
\[ \opn{gr}_{\mrm{G}}(\Psi_{\bsym{\sigma}; 1}) = 
\opn{gr}_{\mrm{G}}(\bsym{\sigma}^*(\mcal{U}_{\mcal{A}; 1})) =
\opn{gr}_{\mrm{G}}( \opn{Mix}^{}_{\bsym{U}}(\mcal{U}_1)) . \]
According to Theorems \ref{thm3.6} and \ref{thm6.6}
the homomorphism $\opn{Mix}^{}_{\bsym{U}}(\mcal{U}_1)$ is a 
quasi-isomorphism. Since the complexes 
$\opn{Mix}^{}_{\bsym{U}}(\mcal{T}^{}_{\mrm{poly}, X})$
and
$\opn{Mix}^{}_{\bsym{U}}(\mcal{D}^{}_{\mrm{poly}, X})$
are bounded below, and the filtration is nonnegative and exhaustive,
it follows that $\Psi_{\bsym{\sigma}; 1}$ 
is also a quasi-isomorphism.
\end{proof}

\begin{cor} \label{cor11.1}
Taking global sections in Theorem \tup{\ref{thm5.1}} we get an 
$\mrm{L}_{\infty}$ quasi-iso\-mor\-phism
\[ \Gamma(X, \Psi_{\bsym{\sigma}}) = 
\{ \Gamma(X, \Psi_{\bsym{\sigma}; j}) \}_{j \geq 1} :
\Gamma \big( X, \opn{Mix}_{\bsym{U}} 
(\mcal{T}^{}_{\mrm{poly}, X}) \bigl) \to
\Gamma \big( X, \opn{Mix}_{\bsym{U}}
(\mcal{D}^{}_{\mrm{poly}, X}) \bigl) . \]
\end{cor}

\begin{proof}
By Theorem \ref{thm6.6} the homomorphism 
\[ \Gamma(X, \opn{Mix}^{}_{\bsym{U}}(\mcal{U}_1)) :
\Gamma \big( X, \opn{Mix}_{\bsym{U}} 
(\mcal{T}^{}_{\mrm{poly}, X}) \bigl) \to
\Gamma \big( X, \opn{Mix}_{\bsym{U}}
(\mcal{D}^{}_{\mrm{poly}, X}) \bigl) \]
is a quasi-isomorphism, and by the interaction with the filtration 
$\{ \mrm{G}^i \opn{Mix}^{}_{\bsym{U}} \}$ we see that 
$\Gamma(X, \opn{Mix}^{}_{\bsym{U}}(\Psi_{\bsym{\sigma}; 1}))$
is also a quasi-isomorphism. 
\end{proof}

\begin{cor} \label{cor5.5}
The data $(\bsym{U}, \bsym{\sigma})$ induces a bijection
\[ \begin{aligned}
\opn{MC}(\Psi_{\bsym{\sigma}}) & :\
\opn{MC} \Bigl( \Gamma \bigl( X, \opn{Mix}_{\bsym{U}}
(\mcal{T}^{}_{\mrm{poly}, X}) \bigl)[[\hbar]]^+ \Bigr) \\
& \qquad \iso \
\opn{MC} \Bigl( \Gamma \bigl( X, \opn{Mix}_{\bsym{U}}
(\mcal{D}^{}_{\mrm{poly}, X}) \bigl)[[\hbar]]^+ \Bigr) .
\end{aligned} \]
\end{cor}

\begin{proof}
Use Corollaries \ref{cor11.1} and \ref{cor7.1}.
\end{proof}

Recall that
$\mcal{T}_{\mrm{poly}}(X) = 
\Gamma(X, \mcal{T}_{\mrm{poly}, X})$
and
$\mcal{D}^{\mrm{nor}}_{\mrm{poly}}(X) = 
\Gamma(X, \mcal{D}^{\mrm{nor}}_{\mrm{poly}, X})$; 
and the latter is the DG Lie algebra of global poly differential 
operators that vanish if one of their arguments is $1$.

\begin{thm} \label{thm5.2}
Assume $\mrm{H}^q(X, \mcal{D}^{\mrm{nor}, p}_{\mrm{poly}, X}) = 0$
for all $p$ and all $q>0$. Then there is a canonical function 
\[ Q : 
\mrm{MC} \bigl( \mcal{T}_{\mrm{poly}}(X)[[\hbar]]^+ \bigr) 
\iso 
\mrm{MC} \bigl( \mcal{D}^{\mrm{nor}}_{\mrm{poly}}(X)[[\hbar]]^+ 
\bigr) \]
preserving first order terms. If moreover 
$\mrm{H}^q(X, \mcal{T}^{p}_{\mrm{poly}, X}) = 0$
for all $p$ and all $q>0$, then $Q$ is bijective. The function 
$Q$ is called the {\em quantization map}, and it 
is characterized as follows. Choose an open covering
$\bsym{U} = \{ U_{(0)}, \ldots, U_{(m)} \}$ 
of $X$ consisting of affine open sets, each admitting 
an \'etale morphism $U_{(i)} \to \mbf{A}^n_{\K}$. Let 
$\bsym{\sigma}$ be the associated simplicial section of
$\opn{LCC} X \to X$. Then there is a commutative 
diagram 
\begin{equation} \label{eqn7.3}
\begin{CD}
\mrm{MC} \bigl( \mcal{T}^{}_{\mrm{poly}}(X)[[\hbar]]^+ \bigr) 
@> Q >>
\mrm{MC} \bigl( \mcal{D}^{\mrm{nor}}_{\mrm{poly}}(X)[[\hbar]]^+ 
\bigr) \\
@V \opn{MC}(\opn{inc}) VV  @V \opn{MC}(\opn{inc}) VV \\
\mrm{MC} \bigl( 
\Gamma(X, \opn{Mix}_{\bsym{U}}  
(\mcal{T}^{}_{\mrm{poly}, X}))[[\hbar]]^+ \bigr)
@> \opn{MC}(\Psi_{\bsym{\sigma}}) >>
\mrm{MC} \bigl( 
\Gamma(X, \opn{Mix}_{\bsym{U}}
(\mcal{D}^{}_{\mrm{poly}, X}))[[\hbar]]^+ \bigr) 
\end{CD} 
\end{equation}
in which the right vertical arrow is bijective. 
Here $\Psi_{\bsym{\sigma}}$ is the $\mrm{L}_{\infty}$ 
quasi-isomorphism from Theorem \tup{\ref{thm5.1}}, and 
``$\mrm{inc}$'' denotes the various inclusions of DG Lie
algebras.
\end{thm}

Let's elaborate a bit on the statement above. It says that 
to any formal Poisson structure 
$\alpha = \sum_{j = 1}^{\infty} \alpha_j \hbar^j
\in \mcal{T}^{1}_{\mrm{poly}}(X)[[\hbar]]^+$
there corresponds a star product $\star_{\beta}$, with  
$\beta = \sum_{j = 1}^{\infty} \beta_j \hbar^j 
\in \mcal{D}^{\mrm{nor}, 1}_{\mrm{poly}}(X)[[\hbar]]^+$
(cf.\ Proposition \ref{prop7.3}). 
The element $\beta = Q(\alpha)$ 
is uniquely determined up to gauge equivalence by 
$\opn{exp}(\mcal{D}^{\mrm{nor}, 0}_{\mrm{poly}}(X)[[\hbar]]^+)$.
Given any local sections $f, g \in \mcal{O}_{X}$ one has
\[ \beta_1(f, g) - \beta_1(g, f) = 2 \{ f, g \}_{\alpha_1} . \]
The quantization map $Q_{}$ can be calculated (at 
least in theory)
using the collection of sections $\bsym{\sigma}$ and the universal 
formulas for deformation in Theorem \ref{thm2.0}. 

We'll need a lemma before proving the theorem.

\begin{lem} \label{lem11.2}
Let $f, g \in \mcal{O}_{X} = \mcal{D}^{-1}_{\mrm{poly}, X}$
be local sections.
\begin{enumerate}
\item For any  
$\beta \in 
\opn{Mix}^0_{\bsym{U}}(\mcal{D}^{1}_{\mrm{poly}, X})$
one has 
\[ [[ \beta, f], g] = \beta(g, f) - \beta(f, g) 
\in \opn{Mix}^0_{\bsym{U}}(\mcal{O}_{X}) . \]
\item For any 
$\beta \in 
\opn{Mix}^1_{\bsym{U}}(\mcal{D}^{0}_{\mrm{poly}, X})
\oplus \opn{Mix}^2_{\bsym{U}}(\mcal{D}^{-1}_{\mrm{poly}, X})$
one has $[[ \beta, f], g] = 0$.
\item Let 
$\gamma \in 
\opn{Mix}_{\bsym{U}}(\mcal{D}^{}_{\mrm{poly}, X})^0$,
and define
$\beta := (\d_{\mrm{mix}} + \d_{\mcal{D}})(\gamma)$.
Then $[[ \beta, f], g] = 0$.
\end{enumerate}
\end{lem}

\begin{proof}
(1) Proposition \ref{prop6.2} implies that the embedding 
(\ref{eqn5.3}):
\[ \begin{aligned}
& \opn{Mix}_{\bsym{U}}(\mcal{D}^{}_{\mrm{poly}, X}) \\
& \qquad \subset \
\bigoplus_{p,q,r} \
{\displaystyle \prod_{j \in \mbb{N}}}\ \
{\displaystyle \prod_{\bsym{i} \in 
\bsym{\Delta}_j^{m}}} \,
g_{\bsym{i} *}\, g_{\bsym{i}}^{-1}\,
\big( \Omega^{q}(\bsym{\Delta}^j_{\K})
\, \what{\otimes} \, (\Omega^{p}_X \otimes_{\mcal{O}_{X}}
\mcal{P}_X \otimes_{\mcal{O}_{X}} 
\mcal{D}^{r}_{\mrm{poly}, X}) \big)
\end{aligned} \]
is a DG Lie algebra homomorphism.
So by continuity we might as well assume that 
$\beta = a D$ with 
$a \in \Omega^0_X = \mcal{O}_X$ and 
$D \in \mcal{D}^{1}_{\mrm{poly}, X}$.
Moreover, since the Lie bracket of 
$\Omega^{}_X \otimes_{\mcal{O}_{X}}
\mcal{P}_X \otimes_{\mcal{O}_{X}} \mcal{D}^{}_{\mrm{poly}, X}$
is $\Omega_X$-bilinear, we may assume that $a = 1$, 
i.e.\ $\beta = D$. Now the assertion
is clear from the definition 
of the Gerstenhaber Lie bracket, see \cite[Section 3.4.2]{Ko1}.

\medskip \noindent
(2) Applying the same reduction as above, but with 
$D \in \mcal{D}^{r}_{\mrm{poly}, X}$
and $r \in \set{0, -1}$, we get
$[[D, f], g] \in \mcal{D}^{r-2}_{\mrm{poly}, X} = 0$.

\medskip \noindent
(3) By part (2) it suffices to show that 
$[[ \beta, f], g] = 0$ for
$\beta := \d_{\mcal{D}}(\gamma)$ and
$\gamma \in \opn{Mix}^{0}_{\bsym{U}}
(\mcal{D}^{0}_{\mrm{poly}, X})$.
As explained above we may further assume that 
$\gamma = D \in \mcal{D}^{0}_{\mrm{poly}, X}$. 
Now the formulas for $\d_{\mcal{D}}$ and $[-,-]$
in \cite[Section 3.4.2]{Ko1} imply that \linebreak
$[[\d_{\mcal{D}}(D), f], g] = 0$. 
\end{proof}

\begin{proof}[Proof of Theorem \tup{\ref{thm5.2}}]
We are given that 
$\mrm{H}^q(X, \mcal{D}^{\mrm{nor}, p}_{\mrm{poly}, X}) = 0$
for all $p$ and all $q>0$; and therefore 
$\Gamma(X, \mcal{D}^{\mrm{nor}}_{\mrm{poly}, X}) = 
\mrm{R} \Gamma(X, \mcal{D}^{\mrm{nor}}_{\mrm{poly}, X})$
in the derived category $\msf{D}(\cat{Mod} \K)$.
Now by Theorem \ref{thm3.6} the inclusion 
$\mcal{D}^{\mrm{nor}}_{\mrm{poly}, X} \to
\mcal{D}^{}_{\mrm{poly}, X}$
is a quasi-isomorphism, and by Theorem \ref{thm6.6}(1)
the inclusion 
$\mcal{D}^{}_{\mrm{poly}, X} \to 
\opn{Mix}_{\bsym{U}}(\mcal{D}^{}_{\mrm{poly}, X})$
is a quasi-isomorphism. According to Theorem \ref{thm6.6}(2) we 
have
$\Gamma \bigl( X, \opn{Mix}_{\bsym{U}}
(\mcal{D}^{}_{\mrm{poly}, X}) \bigr) = 
\mrm{R} \Gamma \bigl( X, \opn{Mix}_{\bsym{U}}
(\mcal{D}^{}_{\mrm{poly}, X}) \bigr)$.
The conclusion is that
\begin{equation} \label{eqn7.4}
\mcal{D}^{\mrm{nor}}_{\mrm{poly}}(X) = 
\Gamma(X, \mcal{D}^{\mrm{nor}}_{\mrm{poly}, X}) \to
\Gamma \bigl( X, \opn{Mix}_{\bsym{U}}
(\mcal{D}^{}_{\mrm{poly}, X}) \bigr)
\end{equation}
is a quasi-isomorphism of complexes of $\K$-modules. 
But in view of Proposition \ref{prop6.2} 
this is also a homomorphism of DG Lie algebras. 

From (\ref{eqn7.4}) we deduce that 
\[ \mcal{D}^{\mrm{nor}}_{\mrm{poly}}(X)[[\hbar]]^+ \to
\Gamma \bigl( X, \opn{Mix}_{\bsym{U}}
(\mcal{D}^{}_{\mrm{poly}, X}) \bigr)[[\hbar]]^+ \]
is a quasi-isomorphism of DG Lie algebras. Using 
Corollary \ref{cor7.1} we see that the right vertical arrow
in the diagram (\ref{eqn7.3}) is bijective.
Therefore this diagram defines $Q_{}$ uniquely.

According to Corollary \ref{cor5.5} the bottom arrow in 
diagram (\ref{eqn7.3}) is a bijection. The left vertical arrow 
comes from the DG Lie algebra homomorphism
\[ \mcal{T}_{\mrm{poly}}(X)[[\hbar]]^+ \to
\Gamma \bigl( X, \opn{Mix}_{\bsym{U}}
(\mcal{T}^{}_{\mrm{poly}, X}) \bigr)[[\hbar]]^+ , \]
which is a quasi-isomorphism when 
$\mrm{H}^q(X, \mcal{T}^{p}_{\mrm{poly}, X}) = 0$
for all $p$ and all $q>0$. So in case of this further vanishing of 
cohomology the map $Q$ is bijective. 

Now suppose $\bsym{U}' = \{ U'_{(0)}, \ldots, U'_{(m')} \}$ 
is another such covering of $X$, with sections
$\sigma'_{(i)} : U'_{(i)} \to \opn{LCC} X$. 
Without loss of generality we may assume that $m' \geq m$, and 
that $U'_{(i)} = U_{(i)}$ and $\sigma'_{(i)} = \sigma_{(i)}$
for all $i \leq m$. There is a morphism of simplicial schemes
$f : \bsym{U} \to \bsym{U}'$, that is an open and closed embedding.
Correspondingly there is a commutative diagram 
\[ \begin{CD}
\mrm{MC} \bigl( \mcal{T}^{}_{\mrm{poly}}(X)[[\hbar]]^+ \bigr) 
@> Q >>
\mrm{MC} \bigl( \mcal{D}^{\mrm{nor}}_{\mrm{poly}}(X)[[\hbar]]^+ 
\bigr) \\
@V \opn{MC}(\opn{inc}) VV  @V \opn{MC}(\opn{inc}) VV \\
\mrm{MC} \bigl( 
\Gamma(X, \opn{Mix}_{\bsym{U}'}  
(\mcal{T}^{}_{\mrm{poly}, X}))[[\hbar]]^+ \bigr)
@> \opn{MC}(\Psi_{\bsym{\sigma}'}) >>
\mrm{MC} \bigl( 
\Gamma(X, \opn{Mix}_{\bsym{U}'}
(\mcal{D}^{}_{\mrm{poly}, X}))[[\hbar]]^+ \bigr) \\
@V \mrm{MC}(f^*) VV  @V \opn{MC}(f^*) VV \\
\mrm{MC} \bigl( 
\Gamma(X, \opn{Mix}_{\bsym{U}}
(\mcal{T}^{}_{\mrm{poly}, X}))[[\hbar]]^+ \bigr)
@> \opn{MC}(\Psi_{\bsym{\sigma}}) >>
\mrm{MC} \bigl( 
\Gamma(X, \opn{Mix}_{\bsym{U}}
(\mcal{D}^{}_{\mrm{poly}, X}))[[\hbar]]^+ \bigr)  ,
\end{CD} \]
where the vertical arrows on the right
are bijections. We conclude that $Q$ is 
independent of $\bsym{U}$ and $\bsym{\sigma}$. 

Finally we must show that $Q$ preserves first order terms. Let
$\alpha = \sum_{j=1}^{\infty} \alpha_j \hbar^j$ be a formal Poisson 
structure, and let 
$\beta = \sum_{j=1}^{\infty} \beta_j \hbar^j
\in \mcal{D}^{\mrm{nor}}_{\mrm{poly}}(X)^1[[\hbar]]^+$
be a solution of the MC equation, 
such that $\beta = Q(\alpha)$ modulo gauge equivalence.
This means that there exists some 
\[ \gamma = \sum_{k \geq 1} \gamma_k \hbar^k \in
\Gamma(X, \opn{Mix}_{\bsym{U}}
(\mcal{D}^{}_{\mrm{poly}, X}))^0[[\hbar]]^+ \]
such that
\[ \sum_{j \geq 1} \smfrac{1}{j!}
\Psi_{\bsym{\sigma}; j}(\alpha^j) = 
\opn{exp}(\opn{af})(\opn{exp}(\gamma))(\beta) , \]
with notation as in Lemma \ref{lem7.3}.
In the first order term (i.e.\ the coefficient of $\hbar^1$) of 
this equation we have
\begin{equation} \label{eqn10.5}
\Psi_{\bsym{\sigma}; 1}(\alpha_1) = \beta_1 - 
(\d_{\mrm{mix}} + \d_{\mcal{D}})(\gamma_1) ;
\end{equation}
see equation (\ref{eqn8.2}). Now by definition (see proof of 
Theorem \ref{thm5.1})
\[ \Psi_{\bsym{\sigma}; 1}(\alpha_1) = 
\bsym{\sigma}^*(\mcal{U}_{\mcal{A}, \omega; 1}) =
\sum_{k \geq 0} \smfrac{1}{(1+k)!}
\bsym{\sigma}^*(\mcal{U}_{\mcal{A}; 1+k}) \big( 
\bsym{\sigma}^*(\omega_{\mrm{MC}})^k \alpha_1 \big) , \]
and the component in 
$\Gamma(X, \opn{Mix}^0_{\bsym{U}}
(\mcal{D}^{1}_{\mrm{poly}, X}))[[\hbar]]^+$ is the summand with 
$k=0$, namely 
$\bsym{\sigma}^*(\mcal{U}_{\mcal{A}; 1}) (\alpha_1) = 
\mcal{U}_{1}(\alpha_1)$. Using Lemma \ref{lem11.2} we get
\[ [[\Psi_{\bsym{\sigma}; 1}(\alpha_1), f], g] = 
[[ \mcal{U}_{1}(\alpha_1), f], g] =
\mcal{U}_{1}(\alpha_1)(g, f) - \mcal{U}_{1}(\alpha_1)(f, g) = 
- 2 \{ f, g \}_{\alpha_1} ,  \]
\[ [[ \beta_1, f], g] = \beta_1(g, f) - \beta_1(f, g) \]
and
\[ [[ (\d_{\mrm{mix}} + \d_{\mcal{D}})(\gamma_1), f], g] = 0  \]
for every local sections $f, g \in \mcal{O}_{X}$.
Combining these equations with equation (\ref{eqn10.5}) the proof 
is done.
\end{proof}

One says that $X$ is a {\em $\mcal{D}$-affine variety} if 
$\mrm{H}^q(X, \mcal{M}) = 0$ for every quasi-coherent left
$\mcal{D}_X$-module $\mcal{M}$ and every $q > 0$.

\begin{cor} \label{cor9.1}
Assume $X$ is $\mcal{D}$-affine. Then the quantization map $Q$
of Theorem \tup{\ref{thm5.2}} may be interpreted as a 
canonical function
\[ Q : 
\frac{ \{ \tup{formal Poisson structures on $X$} \} }
{\tup{gauge equivalence}} \iso
\frac{ \{ \tup{deformation quantizations of $\mcal{O}_X$} \} }
{\tup{gauge equivalence}} \]
preserving first order terms. If $X$ is affine then $Q$ is 
bijective.
\end{cor}

\begin{proof}
By definition the left side is 
$\mrm{MC} \bigl( \mcal{T}^{}_{\mrm{poly}}(X)[[\hbar]]^+ \bigr)$.
On the other hand, according to Theorem \ref{thm4.5} every 
deformation quantization of $\mcal{O}_{X}$ can be trivialized 
globally, and by Proposition \ref{prop6.1} any gauge equivalence 
between globally trivialized deformation quantizations is a
global gauge equivalence. Hence the right side is 
$\mrm{MC} \bigl( \mcal{D}^{\mrm{nor}}_{\mrm{poly}}(X)[[\hbar]]^+ 
\bigr)$.
Since each $\mcal{D}^{\mrm{nor}, p}_{\mrm{poly}, X}$
is a quasi-coherent left $\mcal{D}_X$-module, and each
$\mcal{T}^{p}_{\mrm{poly}, X}$ is a quasi-coherent  
$\mcal{O}_X$-module, we can apply Theorem \ref{thm5.2}.
\end{proof}

Suppose $f : X' \to X$ is an \'etale morphism. According to 
\cite[Prposition 4.6]{Ye2} there are DG Lie algebra homomorphisms
$f^* : \mcal{T}^{}_{\mrm{poly}}(X) \to 
\mcal{T}^{}_{\mrm{poly}}(X')$
and
$f^* : \mcal{D}^{\mrm{nor}}_{\mrm{poly}}(X) \to 
\mcal{D}^{\mrm{nor}}_{\mrm{poly}}(X')$. Given a formal Poisson 
structure $\alpha$ on $X$ we then obtain a formal Poisson 
structure $f^*(\alpha)$ on $X'$. Similarly a star product
$\star$ on $\mcal{O}_{X}[[\hbar]]$ induces a star product
$f^*(\star)$ on $\mcal{O}_{X'}[[\hbar]]$,
 
\begin{cor} \label{cor7.9}
The quantization map $Q$ respects \'etale morphisms. Namely if $X$ 
and $X'$ are $\mcal{D}$-affine schemes and $f : X' \to X$ is an 
\'etale morphism, then for any formal Poisson structure 
$\alpha$ on $X$ one has
$Q(f^*(\alpha)) = f^*(Q(\alpha))$.
\end{cor}

\begin{proof}
This is clear from the proof of Theorem \ref{thm5.2}.
\end{proof}

\section{Complements and Remarks}

Suppose $C$ is some smooth commutative $\K$-algebra, where $\K$ is 
a field containing $\mbb{R}$. 
It is conceivable to look for a star product on $C[[\hbar]]$ that 
is {\em non-differential}. Namely, a $\K[[\hbar]]$-bilinear, 
associative, unital multiplication $\star$ on $C[[\hbar]]$ of the 
form
\[ f \star g = f g + \sum_{k = 1}^{\infty} 
\beta_{k}(f, g) \hbar^k , \]
where the normalized $\K$-bilinear functions $\beta_k : C^2 \to C$ 
are not necessarily bi-differential operators. Indeed, classically 
this was the type of deformation that had been considered (cf.\ 
\cite{Ge}). There is a corresponding notion of non-differential 
gauge equivalence, via an automorphism 
$\gamma = \bsym{1}_C + \sum_{k = 1}^{\infty} \gamma_{k} \hbar^k$
of $C[[\hbar]]$ with $\gamma_k : C \to C$ normalized $\K$-linear 
functions. 

\begin{prop} \label{prop9.5}
Let $C$ be a smooth $\K$-algebra. Then the obvious function
\[ \frac{ \{ \tup{star products on $C[[\hbar]]$} \} }
{\tup{gauge equivalence}}
\to
\frac{ \{ \tup{non-differential star products 
on $C[[\hbar]]$} \} }
{\tup{non-differential gauge equivalence}} \]
is bijective. 
\end{prop} 

\begin{proof}
Let us denote by $\mcal{G}(C)$ the shifted full Hochschild 
cochain complex of $C$, and let $\mcal{G}^{\mrm{nor}}(C)$
be the subcomplex of normalized cochains. It is a well-known fact 
that the inclusion 
$\mcal{G}^{\mrm{nor}}(C) \inj \mcal{G}(C)$ is a quasi-isomorphism
(it is an immediate consequence of \cite[Corollary X.2.2]{ML}). 
By \cite[Lemma 4.3]{Ye1} the $C$-linear map
$\mcal{U}_1 : \mcal{T}^{}_{\mrm{poly}}(C) \to \mcal{G}(C)$
is a quasi-isomorphism, and by Theorem \ref{thm3.6} the map
$\mcal{U}_1 : \mcal{T}^{}_{\mrm{poly}}(C) \to 
\mcal{D}^{\mrm{nor}}_{\mrm{poly}}(C)$
is a quasi-isomorphism. The upshot is that the inclusion
$\mcal{D}^{\mrm{nor}}_{\mrm{poly}}(C) \inj 
\mcal{G}^{\mrm{nor}}(C)$
is a quasi-isomorphism of DG Lie algebras. 
Now we can use Propositions \ref{prop7.3} and \ref{prop3.10}, as 
well as their ``classical'' non-differential variants 
(see proof of \cite[Corollary 4.5]{Ke}).
\end{proof}

Combining Proposition  \ref{prop9.5} 
with Corollary \ref{cor9.1} (applied to 
$X := \opn{Spec} C$) we obtain:

\begin{cor}
Let $C$ be a smooth $\K$-algebra. 
Then there is a canonical bijection of sets
\[ \begin{aligned}
& Q : \frac{ \{ \tup{formal Poisson structures on $C$} \} }
{\tup{gauge equivalence}}  \\
& \qquad \qquad  \iso
\frac{ \{ \tup{non-differential star products 
on $C[[\hbar]]$} \} }
{\tup{non-differential gauge equivalence}} 
\end{aligned} \]
preserving first order terms. 
\end{cor}

\begin{que}
In case $X$ is affine and admits an \'etale morphism 
$X \to \mbf{A}^n_{\K}$, how are the the deformation quantizations 
of Corollary \ref{cor9.1} Corollary \ref{cor6.1} related?
\end{que}

\begin{rem}
The methods of this paper, combined with the ideas of 
\cite{CFT}, can be used to prove the following 
result. Suppose $\mbb{R} \subset \K$ and 
$\mrm{H}^2(X, \mcal{O}_X) = 0$. Let $\alpha$ be any Poisson 
structure on $X$. Then the Poisson variety $(X, \alpha)$ admits a 
deformation quantization, in the sense of Definition \ref{dfn0.1}.
\end{rem}

\begin{que}
Given a smooth scheme $X$, is it possible to determine which 
Poisson structures on $X$ can be quantized? The papers \cite{NT} 
and \cite{BK1} say that for a symplectic structure to be 
quantizable there are cohomological obstructions. Can anything 
like that be done for a degenerate Poisson structure?
\end{que}

\begin{rem}
Artin worked out a noncommutative deformation theory for 
\linebreak schemes
that goes step by step, from $\K[\hbar] / (\hbar^m)$ to 
$\K[\hbar] / (\hbar^{m+1})$; see \cite{Ar1} and \cite{Ar2}. 
The first order data is a Poisson structure, and at each step 
there are well defined obstructions to the process. Presumably 
Artin's deformations are deformation quantizations  
in the sense of Definition \ref{dfn4.5}, namely they admit 
differential structures; but this requires a proof.

In the case of the projective plane $\mbf{P}^2$ and a nonzero 
Poisson structure $\alpha$, the zero locus of $\alpha$ is a cubic 
divisor $E$. Assume $E$ is smooth. Artin asserts (private 
communication) that a particular deformation of 
$\mcal{O}_{\mbf{P}^2}$ with first order term $\alpha$ lifts to 
a deformation of the homogeneous coordinate ring 
$B := 
\boplus_{i \geq 0} \Gamma(\mbf{P}^2, \mcal{O}_{\mbf{P}^2}(i))$.
Namely there is a graded $\K[[\hbar]]$-algebra structure on 
$\boplus_{i \geq 0} B_i[[\hbar]]$,
say $\star_{}$, such that 
$a \star_{} b \equiv a b \ \opn{mod} \hbar$,
and
$a \star_{} b - b \star_{} a \equiv 
2 \hbar \set{a, b}_{\alpha} \opn{mod} \hbar^2$,
for all $a, b \in B$. Moreover, after tensoring with the field
$\K((\hbar))$ this should be 
a three dimensional Sklyanin algebra, presumably 
with associated elliptic curve
$\K((\hbar)) \times_{\K} E$.

The above should be compared to \cite[Section 3]{Ko3}. 
\end{rem}


\end{document}